\documentclass[11pt,3p]{elsarticle}

\makeatletter
\def\ps@pprintTitle{%
 \let\@oddhead\@empty
 \let\@evenhead\@empty
 \def\@oddfoot{\centerline{\thepage}}%
 \let\@evenfoot\@oddfoot}
\makeatother

\usepackage{graphicx}
\newcommand{\dt}{\mathrm{d}t}

\newcommand{\dd}{\mathrm{d}}
\usepackage{dsfont}
\usepackage{lineno}
\usepackage{amsmath}
\usepackage{amssymb}
\usepackage{wrapfig}
\usepackage{comment}
\usepackage{subfigure}
\usepackage{floatflt}
\usepackage{dsfont}
\usepackage{overpic}
\usepackage{pgfplots}
\pgfplotsset{width=7cm,compat=1.3}
\usepackage{mwe}
\usepackage{calc}
\usepackage{xcolor}
\usepackage{footnote}
\makesavenoteenv{tabular}
\makesavenoteenv{table}
\pgfplotsset{compat=newest} 
\pgfplotsset{plot coordinates/math parser=false} 
\newlength\figureheight 
\newlength\figurewidth 

\newcommand{\ii}{\textbf{i}}

\newcommand{\R}{\mathbb{R}}

\newcommand\scalemath[2]{\scalebox{#1}{\mbox{\ensuremath{\displaystyle #2}}}}

\newcommand{\sus}[1]{{#1}}

\begin{document}

\begin{frontmatter}

\title{On complex dynamics in a Purkinje and a ventricular cardiac cell model}

\author[AE]{Andr\'e H. Erhardt\corref{mycorrespondingauthor}\fnref{contribution}}
\address[AE]{Department of Mathematics, University of Oslo,
           P.O.Box 1053 Blindern, 0316 Oslo, Norway}
\ead{andreerh@math.uio.no}
\author[SS]{Susanne Solem\fnref{contribution}}
\address[SS]{Department of Mathematical Sciences, Norwegian University of Science and Technology, 7491 Trondheim, Norway}
\ead{susanne.solem@ntnu.no}
\cortext[mycorrespondingauthor]{Corresponding author}
\fntext[contribution]{The authors have contributed equally to the content of this manuscript.}

\begin{abstract}
 Cardiac muscle cells can exhibit complex patterns including irregular behaviour such as chaos or (chaotic) early afterdepolarisations (EADs), which can lead to sudden cardiac death. Suitable mathematical models and their analysis help to predict the occurrence of such phenomena and to decode their mechanisms. The focus of this paper is the investigation of dynamics of cardiac muscle cells described by systems of ordinary differential equations. This is generically performed by studying a Purkinje cell model and a modified ventricular cell model. We find chaotic dynamics with respect to the leak current in the Purkinje cell model, and EADs and chaos with respect to a reduced fast potassium current and an enhanced calcium current in the ventricular cell model --- features that have been experimentally observed and are known to exist in some models, but are new to the models under present consideration. We also investigate the related monodomain models of both systems to study synchronisation and the behaviour of the cells on macro-scale in connection with the discovered features. 
 The models show qualitatively the same behaviour to what has been experimentally observed. However, for certain parameter settings the dynamics occur within a non-physiological range. 
\end{abstract}


\end{frontmatter}

\section{Introduction}\label{sec:intro}
\noindent Nowadays, mathematical modelling and numerical simulations are essential and standard approaches to study and analyse real world problems and phenomena in life science. One major aim is the understanding of complex dynamics and behaviour of these systems. For this purpose, bifurcation theory has proven to be a very helpful and powerful tool in order to investigate dynamical systems and their (complex) dynamics, see~\cite{Guckenheimer,Kuznetsov,ShilnikovPartI,ShilnikovPartII,Wiggins} for an overview. Furthermore, numerical bifurcation analysis has become a profitable tool in the study of (for instance) climate, neuronal and cardiac models. 

Cardiac muscle cells can exhibit complex patterns of oscillations like spiking and bursting, which is related to ion current interactions of the considered cell. Aside from normal action potentials of a cardiac muscle cell, certain kinds of cardiac arrhythmia can occur. This includes specific types of abnormal heart rhythms, which can lead to sudden cardiac death. In addition, irregular behaviour, such as (deterministic) chaos or chaotic early afterdepolarisations, has been observed in experimental as well as in computational studies, see~\cite{QU_review,Sato} and the references therein. It is therefore highly interesting and important to understand the complex behaviour and mechanism of such biological phenomena. Moreover, cardiac dynamics or heart rhythms can be quite sensitive to the influence of certain drugs, which has been investigated experimentally but also in computational studies, see e.g.~\cite{drugs1,drugs2,Rodriguez_drug}. In later years, the focus has been to continuously move towards interdisciplinary research, including biology, computer science, and mathematics, to tackle these issues. As a consequence, the number of existing mathematical models based on experimental data is also continuously increasing.

The development of a good and precise mathematical model is essential to design numerical experiments for the study of cardiac dynamics, but also for the investigation of the influence of certain external effects such as drugs or oxidative stress. To this end, mathematical analysis is key to decode occurring phenomena and to validate a derived model in all details. The newly gained information of the considered model, can then be either used to improve the model or to proceed with the original aspiration, e.g. the investigation of optimal properties of drugs~\cite{drugs2}.

To this end, bifurcation theory has been utilised to investigate the dynamics of cardiac muscle cells in recent years, see e.g.~\cite{Kurata,Tran,Tsumoto2017,Xie}. Continuing on this line of research, this paper highlights how useful bifurcation theory can be for the understanding of complex cardiac dynamics and how it can be applied to find hidden features and dynamics of cardiac single cell models described by ordinary differential equations (ODEs) through numerical investigations. 

In the end, it is the synchronisation of a large group of cells that decides whether a cardiac arrhythmia spreads or dies out. For this reason, a brief study of how the micro-scale single cell features of these models affect the behaviour of an ensemble of cells at the macro-scale level ($cm$) is provided. This is done by an up-scaling of the ODE system to a PDE--ODE monodomain model.

All of the above will be done by an in-depth mathematical and numerical investigation of the two cardiac cell models introduced in~\cite{Noble,Bernus}, where one is a model of a Purkinje cell, and the other a model of a human ventricular cell. In particular, we find new features of the considered models, such as chaos and early afterdepolarisations, and investigate how this affects groups of cells at the tissue level.

We find chaotic dynamics in both models considered. Similar chaotic dynamics to what we discover can be observed in experiments, cf.~\cite{QU_review}. However, the dynamics seems to appear in a non-physiological range. This either requires the improvement of the models or the experimental validation. 

In addition, we discover EADs in the human ventricular cell model. This behaviour does seem to be within the physiological range \cite{Vandersickel2}, which is a validation of the model in this case. 

Finally, we show that the (in)validity of the models in terms of being within the physiological range carries over to the synchronisation effects in the corresponding monodomain models. 

These findings clearly highlight the advantages of bifurcation theory in the analysis of cardiac muscle cell dynamics by detecting unexpected or maybe non-physiological behaviour of the model. 

\vspace{0.8em}
\noindent \textbf{Outline of the paper.} In Section~\ref{sec:background} we start with a mathematical and biological description of the models and problems under investigation. A brief background on the modelling of cardiac muscle cells and on up-scaling to a monodomain model at the tissue level (cf.~\cite{KeenerII,CLAYTON2009951,KELDERMANN20091000}) is provided. Furthermore, we perform a stability analysis of the ODE system from \cite{Noble}, and show how to extend this analysis to the macro-scale monodomain model. The approach is explained in detail for the four dimensional model~\cite{Noble}. However, the ansatz can also be used for more complex models, see~\cite{modeloverview,FINK20112}. 

The structure of the systems in~\cite{Noble,Bernus} are similar. Nevertheless, the behaviour and dynamics that they display can be quite disparate due to different complexity and parameter settings. In Section~\ref{sec:bif}, we apply a numerical bifurcation analysis in order to derive a complete understanding of the dynamics of the model from~\cite{Noble} with respect to certain parameters. The analysis is then extended to the corresponding macro-scale monodomain model. Based on the results in Section~\ref{sec:bif}, we then continue our analysis by studying a ten dimensional version of the model from~\cite{Bernus} in Section~\ref{sec:Bernus_bif}. Finally, we close our paper with a discussion in Section \ref{sec:conclusion} and then a conclusion.

\section{Biological and mathematical background}\label{sec:background}
\begin{figure}[h]
\centering
\subfigure[Scheme of a cardiac muscle cell.]{\includegraphics[width=0.5\textwidth]{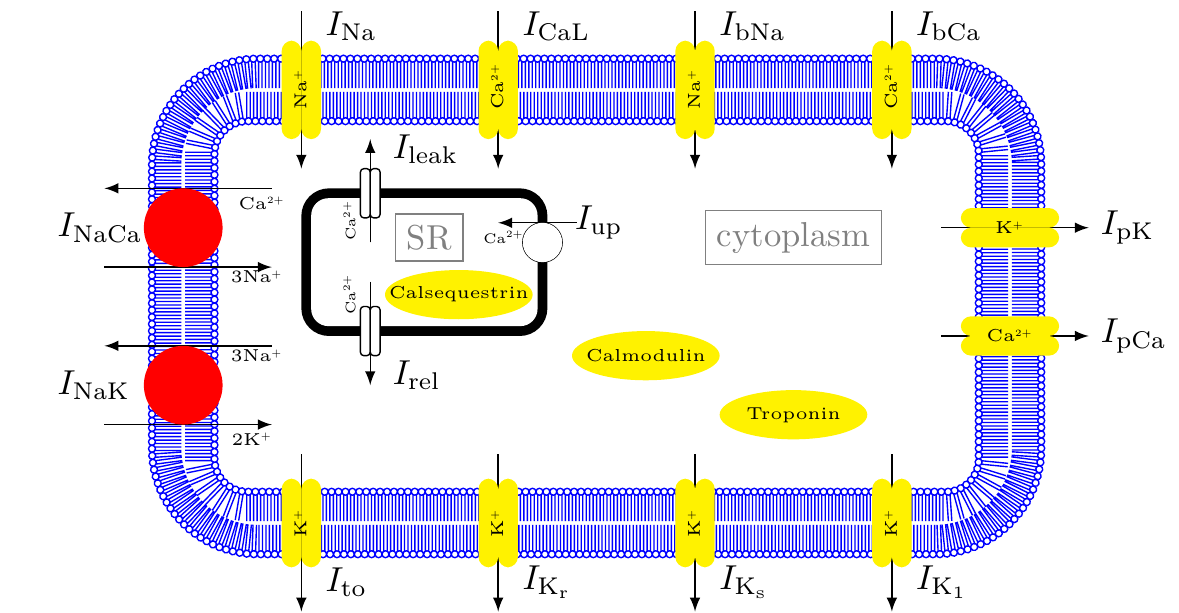}}\subfigure[Physical system of a cardiac muscle cell.]{\includegraphics[width=0.5\textwidth]{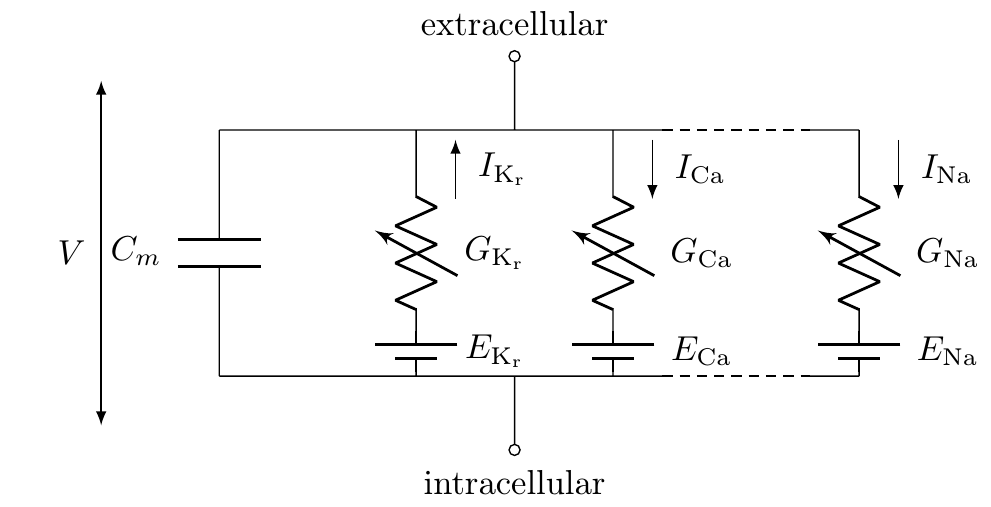}}
\caption{(a) Scheme of a cardiac muscle cell~\cite{TNNP04}, where SR denotes the sarcoplasmic reticulum, $I_\text{NaCa}=\mathrm{Na}^+/\mathrm{Ca}^{2+}$ exchanger current  and  $I_\text{NaK}=\mathrm{Na}^+/\mathrm{K}^+$ pump current. (b) Physical system of a cardiac muscle cell. The dashed lines denote the ion currents, which are not included due to the lack of space.}\label{fig:cell}
\end{figure}

\noindent The history of mathematical modelling of action potentials (APs) of excitable biological cells like neurons and cardiac cells starts with the famous and pioneering Hodgkin--Huxley (HH) model from 1952~\cite{HH}. In~\cite{HH}, the authors established a mathematical approach that can be used to model APs of excitable biological cells by a system of ODEs. The first model of a cardiac cell is the Noble model from 1962~\cite{Noble} of a generic Purkinje cell. In 1991, Luo and Rudy published an ionic model for cardiac action potential in guinea pig ventricular cells~\cite{LR}. In the last decades, there has been an immense development in the modelling of cardiac muscle cells, see e.g.~\cite{Fenton_Karma,TNNP04,TP06,minimalmodel}. These conductance--based models represent a minimal biophysical interpretation of an excitable biological cell in which current flow across the membrane is due to charging of the membrane capacitance and movement of ions across ion channels, cf. Figure~\ref{fig:cell}. Ion channels are selective for particular ionic species. In general, an AP is a temporary, characteristic variance in the membrane potential of an excitable biological cell from its resting potential. The molecular mechanism of an AP is based on the interaction of voltage-sensitive ion channels. The reason for the formation and the special properties of the AP is established in the properties of different groups of ion channels in the plasma membrane. An initial stimulus activates the ion channels as soon as a certain threshold potential is reached. Then, these ion channels break open and/or up allowing an ion current flow, which changes the membrane potential. A normal AP is always uniform and the cardiac muscle cell AP is typically divided into four phases: the resting phase, the upstroke phase, the (long) plateau phase and the repolarisation phase. This mechanism is based on several different currents. One example is the potassium current $I_\text{K}$ which is usually divided into a fast ($I_{\text{K}_r}$) and a slow current ($I_{\text{K}_s}$), cf. scheme in Figure~\ref{fig:cell}(a).

This electrophysiological behaviour can be described by the ordinary differential equation: 

\begin{align*}
	C_m\frac{\dd V}{\dt}=-I_\mathrm{ion}+I_\mathrm{stimulus},
\end{align*}
where $V$ denotes the voltage (in $mV$) and $t$ the time (in $ms$), while $I_\mathrm{ion}$ is the sum of all transmembrane ionic currents. $I_\mathrm{stimulus}$ represents the externally applied stimulus and $C_m$ denotes membrane capacitance. 

The model in~\cite{Noble} contains a single potassium current $I_\text{K}$, while the authors of~\cite{Bernus} merged a fast ($I_{\text{K}_r}$) and a slow current ($I_{\text{K}_s}$) to derive their model, which is based on the system in~\cite{PB}. In this paper, we will slightly modify the model from \cite{Bernus}, i.e. we replace $I_\text{K}$ by the currents $I_{\text{K}_r}$ and $I_{\text{K}_s}$ from~\cite{PB}. Furthermore, the different ion currents may depend on different gating variables, individual ionic conductances $G_i$ and Nernst potentials $E_i$, $i=\mathrm{Na},\mathrm{K},\mathrm{Ca}$, etc., cf. Section~\ref{sec:model}. 

We want to highlight that cardiac cell models usually have different time scales and may exhibit so-called mixed-mode oscillations~\cite{MMOs}, cf.~\cite{AE_MMOs}, and/or chaotic behaviour, cf.~\cite{Tran,GLASS198389,WEISS,DELANGE2012365,KBE}, which can be linked to certain cardiac arrhythmia. 

For instance, if there are depolarising variations of the membrane voltage, we are speaking about afterdepolarisations (ADs). These ADs are divided into early (EADs) and delayed afterdepolarisations (DADs). This division depends on the timing obtaining the AP. EADs occur either in the plateau or in the repolarisation phase of the AP and are benefited by an elongation of the AP, while DADs occur after the repolarisation phase is completed. EADs are resulting, for example, from a reduction of the repolarising $\mathrm{K}^+$ currents or an enhancement in $\mathrm{Ca}^{2+}$ currents, see e.g.~\cite{WEISS}. Triggers for this are congenital disorders of ion channels or the ingestion of medicaments. In general, EADs are additional small amplitude spikes (mathematically speaking mixed-mode oscillations), i.e. pathological voltage oscillations, during the plateau or repolarisation phase. They are caused by ion channel diseases, oxidative stress or drugs. Furthermore, the presence of EADs strongly correlates with the onset of dangerous cardiac arrhythmias, including torsades de pointes (TdP), which is a specific type of abnormal heart rhythm that can lead to sudden cardiac death, see~\cite{Roden,VN,Vandersickel2}. Thus, it is highly important to understand the complex behaviour of such biological phenomena~\cite{PANFILOV20191}. 



Finally, we want to point out that not all of the existing cardiac cell models may include complex dynamics such as EADs or chaos.

\subsection{A Purkinje cardiac cell model}\label{sec:model} 
\noindent First, we focus on the model from~\cite{Noble}, which reads as follows:

\begin{align}
	\frac{\dd V}{\dt}=&-\frac{I_\mathrm{Na}+I_\mathrm{K}+I_\mathrm{L}}{C_m}=:\mathcal{F},\label{model}
\end{align}
with the membrane capacitance $C_\mathrm{m}=12\ \frac{\mu F}{cm^2}$ and the ion currents $I_\mathrm{Na}$ (sodium), $I_\mathrm{K}$ (potasssium), and $I_\mathrm{L}$ (leak current), described by $I_\mathrm{Na}=(G_\mathrm{Na}m^3h+0.14)(V-E_\mathrm{Na}),$

\begin{align*}
 I_\mathrm{K}=\left(G_{\mathrm{K}_1}n^4+G_{\mathrm{K}_2}\exp\left(-\frac{V+90}{50}\right)+\frac{G_{\mathrm{K}_2}}{80}\exp\left(\frac{V+90}{60}\right)\right)(V-E_\mathrm{K}),
\end{align*}
and $I_\mathrm{L}=G_\mathrm{L}(V-E_\mathrm{L})$, respectively. The individual ionic conductances are given by $G_{\mathrm{K}_1}=G_{\mathrm{K}_2}=1.2\ \frac{mS}{cm^2}$, $G_\mathrm{Na}=400\ \frac{mS}{cm^2}$ and $G_\mathrm{L}=0.075\ \frac{mS}{cm^2}$, while the Nernst potentials are given by $E_\mathrm{K} = -100\ mV$, $E_\mathrm{Na} = 40\ mV$ and $E_\mathrm{L}=-60\ mV$. Furthermore, the different gating variables $m$, $h$ and $n$ satisfy the differential equation 

\begin{align}
	\frac{\dd y}{\dt}=a_y(1-y)-b_yy=a_y-(a_y+b_y)y=\frac{y_\infty(V)-y}{\tau_y(V)},\label{diff_gating_variable}
\end{align}
where $y$ represents the gating variables $m$, $h$ and $n$, while $y_\infty:=y_\infty(V)=\frac{a_y}{a_y+b_y}$ 
denotes the equilibrium of the gating variable $y$ and $\tau_y:=\tau_y(V)=\frac{1}{a_y+b_y}$ its relaxation time constant with 


$$
   \begin{alignedat}{3}
     &a_h=0.17\exp\left(-\frac{V+90}{20}\right),  \quad  &b_h=\frac{1}{1+\exp\left(-\frac{V+42}{10}\right)}, \quad &a_m=\frac{0.1(V+48)}{1-\exp\left(-\frac{V+48}{15}\right)},
     \\ &b_m=\frac{0.12(V+8)}{\exp\left(\frac{V+8}{5}\right)-1},&a_n=\frac{0.0001(V+50)}{1-\exp\left(-\frac{V+50}{10}\right)},  \quad  &b_n=0.002\exp\left(-\frac{V+90}{80}\right). & 
   \end{alignedat}
$$
Notice that the gating variables are important for the activation and inactivation of the different ion currents, which is needed for the ion current flows and the resulting action potential. The Noble model~(\ref{model}) describes the long lasting action and pace-maker potentials of the Purkinje fibres of the heart based on the Hodgkin--Huxley formulation~\cite{HH}. While the sodium current is very similar to the one from~\cite{HH}, the potassium current differs in its formulation and the calcium current is missing. Nevertheless, the solutions to system~(\ref{model}) closely resembles the Purkinje fibre action and pace-maker potentials. It is shown that its behaviour in response to applied currents and to changes in ionic permeability corresponds fairly well with that observed experimentally. Furthermore, the Noble model~(\ref{model}) is one of the earliest mathematical models of cardiac APs which is able to produce a certain type of cardiac arrhythmia, so-called alternans in AP duration (APD). These alternans in APD is a result of a period doubling bifurcation with respect to the cycle length, see~\cite{alternansnoble,Fentonalternans}. A period doubling bifurcation is a creation or destruction of a periodic orbit with double the period of the original orbit.  

In~\cite{Noble} the author numerically studied the influence of the conductance $G_\mathrm{L}$ of the leak current on the trajectory of the Noble model~\eqref{model}. We will extend this numerical study by using bifurcation analysis to derive a more detailed insight into the behaviour of the solutions to~(\ref{model}) by varying $G_\mathrm{L}$, see Section~\ref{sec:bif_gl}.

\subsubsection{Stability analysis}
\noindent The steady state or equilibrium of system~(\ref{model}) is determined by $h\equiv h_\infty(V)$, $m\equiv m_\infty(V)$, $n\equiv n_\infty(V)$ and solving the algebraic equation

\begin{align}
\mathcal{F}(V,h_\infty(V),m_\infty(V),n_\infty(V))=0.\label{eq}
\end{align}
This yields $X_\infty =(V_\infty ,h_\infty(V_\infty ),m_\infty(V_\infty ),n_\infty(V_\infty ))^T$, where $X_\infty $ is depending on several system parameters, cf. system~(\ref{model}). Furthermore, the Jacobian of the right hand side of system~(\ref{model}) evaluated at $X_\infty $ is given by

\begin{align}
\mathcal{J}=\left.\left(\scalemath{1}{\begin{array}{cccc} 
\displaystyle\frac{\partial \mathcal{F}}{\partial V}&\displaystyle\frac{\partial \mathcal{F}}{\partial h}&\displaystyle\frac{\partial \mathcal{F}}{\partial m}&\displaystyle\frac{\partial \mathcal{F}}{\partial n}\\[1.5ex]
\displaystyle\frac{1}{\tau_h}\displaystyle\frac{\partial h_\infty}{\partial V}&-\displaystyle\frac{1}{\tau_h}&0&0\\[1.5ex]
\displaystyle\frac{1}{\tau_m}\displaystyle\frac{\partial m_\infty}{\partial V}&0&-\displaystyle\frac{1}{\tau_m}&0\\[1.5ex]
\displaystyle\frac{1}{\tau_n}\displaystyle\frac{\partial n_\infty}{\partial V}&0&0&-\displaystyle\frac{1}{\tau_n}
\end{array}}\right)\right|_{X_\infty },\label{jac}
\end{align}
where we used the fact that

$$
\left.\frac{\partial}{\partial V}\left(\frac{y_\infty-y}{\tau_y}\right)\right|_{y\equiv y_\infty}=\left.\frac{\frac{\partial y_\infty}{\partial V}\tau_y-(y_\infty-y)\frac{\partial \tau_y}{\partial V}}{\tau_y^2}\right|_{y\equiv y_\infty}=\frac{1}{\tau_y}\frac{\partial y_\infty}{\partial V}.
$$
Note that the location and stability of $X_\infty $ is depending on the different system parameters, e.g. varying $G_\mathrm{L}$ changes the location and the stability of $X_\infty $, while varying $C_\mathrm{m}$ changes only the stability. To determine the stability of the equilibrium one has to calculate the solution(s) of the characteristic polynomial

\begin{align}
\mathrm{det}\left(\mathcal{J}-\lambda\mathds{1}_4\right)
=& \lambda^4+a_1\lambda^3+a_2\lambda^2+a_3\lambda+a_4=0,\label{charactersitic_equ}
\end{align}
i.e. the eigenvalues of $\mathcal{J}$, where 

$$a_1:=\displaystyle\frac{1}{\tau_h}+\displaystyle\frac{1}{\tau_m}+\displaystyle\frac{1}{\tau_n}-\displaystyle\frac{\partial \mathcal{F}}{\partial V},$$
$$a_2:=\frac{1}{\tau_h}\left(\frac{1}{\tau_m}-\displaystyle\frac{\partial h_\infty}{\partial V}\displaystyle\frac{\partial \mathcal{F}}{\partial h}-\frac{\partial \mathcal{F}}{\partial V}\right)+\frac{1}{\tau_n}\left(\frac{1}{\tau_h}-\displaystyle\frac{\partial n_\infty}{\partial V}\displaystyle\frac{\partial \mathcal{F}}{\partial n}-\frac{\partial \mathcal{F}}{\partial V}\right)+\frac{1}{\tau_m}\left(\frac{1}{\tau_n}-\displaystyle\frac{\partial m_\infty}{\partial V}\displaystyle\frac{\partial \mathcal{F}}{\partial m}-\frac{\partial \mathcal{F}}{\partial V}\right)$$
\begin{align*}a_3=&\frac{1}{\tau_h}\left(\frac{1}{\tau_m\tau_n}-\left(\frac{1}{\tau_m}+\frac{1}{\tau_n}\right)\left(\frac{\mathcal{F}}{\partial V}+\frac{\partial h_\infty}{\partial V}\displaystyle\frac{\partial \mathcal{F}}{\partial h}\right)-\frac{1}{\tau_n}\frac{\partial n_\infty}{\partial V}\frac{\partial \mathcal{F}}{\partial n}-\frac{1}{\tau_m}\frac{\partial m_\infty}{\partial V}\frac{\partial \mathcal{F}}{\partial m}\right)
\\
&-\frac{1}{\tau_m\tau_n}\left(\frac{\partial \mathcal{F}}{\partial V}+\frac{\partial n_\infty}{\partial V}\displaystyle\frac{\partial \mathcal{F}}{\partial n}+\frac{\partial m_\infty}{\partial V}\displaystyle\frac{\partial \mathcal{F}}{\partial m}\right)
\end{align*}
and

$$a_4:=-\frac{1}{\tau_h\tau_m\tau_n}\left(\displaystyle\frac{\partial h_\infty}{\partial V}\displaystyle\frac{\partial \mathcal{F}}{\partial h}+\displaystyle\frac{\partial m_\infty}{\partial V}\displaystyle\frac{\partial \mathcal{F}}{\partial m}+\displaystyle\frac{\partial n_\infty}{\partial V}\displaystyle\frac{\partial \mathcal{F}}{\partial n}-\frac{\partial \mathcal{F}}{\partial V}\right).$$
In addition, the Routh--Hurwitz criterion implies that all characteristic exponents $\lambda_i$, $i=1,\dots,4$ have negative real parts if and only if the conditions

$$\Delta_1=a_1>0,~\Delta_2=a_1a_2-a_3>0,~\Delta_3=a_3\cdot \Delta_2>0~\text{and}~\Delta_4=\Delta_3-a_1^2a_4>0$$ hold true, see~\cite{ShilnikovPartI,ShilnikovPartII,Yu2005}. Furthermore, if all Hurwitz minors satisfy $\Delta_i>0$ for $i=1,\cdots,\mathfrak{n}-1$ and $\Delta_\mathfrak{n}=0$, where $\mathfrak{n}$ denotes the dimension of the system, we know that the system exhibits an Andronov--Hopf bifurcation. Using $\Delta_4=0$, we get

\begin{align*}0
=&\lambda^4+a_1\lambda^3+a_2\lambda^2+a_3\lambda+a_4=\lambda^4+a_1\lambda^3+a_2\lambda^2+a_3\lambda+\frac{a_1a_2-a_3}{a_1}\frac{a_3}{a_1}
\\
=&\left(\lambda^2+\frac{a_3}{a_1}\right)\left(\lambda^2+a_1\lambda+\frac{a_1a_2-a_3}{a_1}\right),
\end{align*} 
i.e. the equilibrium has a pair of purely imaginary eigenvalues $\lambda_{1,2}=\ii\omega_0$ with $\omega_0=\displaystyle\frac{a_3}{a_1}>0$ and two further eigenvalues 

$$\lambda_{3,4}=\displaystyle\frac{-a_1\pm\sqrt{a_1^2-4\frac{a_1a_2-a_3}{a_1}}}{2}=\displaystyle\frac{-\Delta_1\pm\sqrt{\Delta_1^2-4\frac{\Delta_2}{\Delta_1}}}{2}.$$ 
In general, an Andronov--Hopf bifurcation corresponds to the birth of a limit cycle, when the equilibrium changes stability via a pair of purely imaginary eigenvalues. Usually, an Andronov--Hopf bifurcation is considered as a trigger to oscillatory behaviour in dynamical systems and may cause normal AP and cardiac arrhythmia in a cardiac cell model.

In case that $a_\mathfrak{n}=0$, then the system exhibits a fold or saddle--node or limit point bifurcation, i.e. the equilibrium has at least one eigenvalue which is equal to zero. A limit point bifurcation is a collision and disappearance of two equilibria in dynamical systems, which is a turning point.

\subsection{A monodomain model of the Purkinje cardiac cell model}\label{sec:monodomain}
\noindent Besides the study of cardiac cell models an important focus is the behaviour and dynamics of several cells, i.e. the dynamics on a macro-scale ($cm$), where many cells are connected together and interacting with each other. To this end, we consider the following monodomain model, i.e.  extension of the ODE model~\eqref{model} to a PDE--ODE model including an additional diffusion term: 

\begin{align}
    \begin{split}
    \begin{split}
       C_m\frac{\partial V}{\partial t}  &= \frac{\lambda}{1+\lambda}\frac{1}{\chi}\nabla\cdot\left(M_i\nabla V\right)- \left(I_\mathrm{Na}+I_\mathrm{K}+I_\mathrm{L}\right)\\
       \frac{\partial y}{\partial t}&=\frac{y_\infty(V)-y}{\tau_y(V)}, \quad y=h,m,n
       \end{split}\quad \text{in}\ \Omega,\\
        0&= \vec{\nu}\cdot\left(M_i\nabla V\right) \quad \text{on}\ \partial\Omega,
    \end{split}
    \label{monodomain}
\end{align}
where $\Omega$ is a bounded domain, $M_i$ denotes the intracellular conductivity tensor, $\lambda$ the extra- to intracellular conductivity ratio, $\chi$ is the membrane surface area per unit volume and $\vec{\nu}$ is the unit normal, cf.~\cite{computing_heart,Nielsen,MR3187193}. System \eqref{monodomain} is a monodomain model, meaning that equal anisotropy rates, i.e. $M_e=\lambda M_i$ in $\frac{mS}{cm}$, are assumed in the more complex bidomain model. Here, $\lambda\in\R$ is constant and $M_e$ denotes the intracellular conductivity tensor. Furthermore, we use $\frac{\lambda}{1+\lambda}\frac{M_i}{\chi}=\frac{1}{360}\ mS$, which is the diffusion constant originally used in \cite{Bernus}, unless otherwise stated.



For a better understanding of the behaviour of cells on a macro-scale level, we follow the approach from~\cite{PDE-ODEI,PDE-ODEII,PDE-ODEIII} to derive a linearised system of model~\eqref{monodomain}. First of all, we know that on rectangle-like domains $\Omega:=[0,\ell_1]\times\dots\times [0,\ell_\mathfrak{n}]\subset\R^\mathfrak{n}$ the eigenvalues and eigenfunctions of the Neumann problem

\begin{align}
\begin{cases}
-\Delta u_k(x)=\mu_k u_k(x)\quad  &x\in\Omega,
\\
 \vec{\nu}\cdot(\nabla u_k(x))=0&  x\in\partial\Omega
\end{cases}\label{spectral_problem}
\end{align}
are 

$$
\mu_k^{(i)}=\left(\frac{\pi k}{\ell_i}\right)^2\quad\text{and}\quad u_k^{(i)}=\cos\left(\frac{\pi kx}{\ell_i}\right)\qquad k=0,1,2,\dots, \quad i=1,\dots,\mathfrak{n},
$$
see~\cite{eigenvalues}. In general, the eigenvalue of the spectral problem~\eqref{spectral_problem} can be derived by multiplying the first equation of system~\eqref{spectral_problem}  by $u_k$, integrating over $\Omega$, using Green's formula and applying the boundary condition. Then, one gets for the Neumann problem~\eqref{spectral_problem} the following eigenvalues:

$$
\mu_k=\frac{\|\nabla u_k\|_{L^2(\Omega)}^2}{\| u_k\|_{L^2(\Omega)}^2}.
$$ 
%

Now, we consider the linearised system of the monodomain equation~\eqref{monodomain} around an equilibrium $X_\infty$ of the ODE system \eqref{eq}. It has the form

\begin{align}
\frac{\partial}{\partial t}
\begin{pmatrix}
V\\ h\\ m\\ n
\end{pmatrix}
=\mathcal{D}\Delta
\begin{pmatrix}
V\\ h\\ m\\ n
\end{pmatrix}
+\mathcal{J}
\begin{pmatrix}
V\\ h\\ m\\ n
\end{pmatrix},
\label{eq:pde-linear}
\end{align}
where $\mathcal{D}$ denotes the $4\times4$ diffusion matrix with almost everywhere zero entries except the first one, which is

$$
\mathcal{D}_{11}=\frac{\lambda}{1+\lambda}\frac{M_i}{\chi}\frac{1}{C_m},
$$
while $\mathcal{J}$ is the Jacobian as stated in~\eqref{jac}. Define the linear operator $\mathcal{L}$ as

\begin{align*}
\mathcal{L}
\begin{pmatrix}
V\\ h\\ m\\ n\end{pmatrix}
:=\mathcal{D}\Delta
\begin{pmatrix}
V\\ h\\ m\\ n
\end{pmatrix}
+\mathcal{J}
\begin{pmatrix}
V\\ h\\ m\\ n
\end{pmatrix}.
\end{align*}
Then, consider the following characteristic equation

\begin{align*}
\mathcal{L}
\begin{pmatrix}
\psi_1\\ \vdots\\ \psi_4
\end{pmatrix}=\mu
\begin{pmatrix}
\psi_1\\ \vdots\\ \psi_4
\end{pmatrix},
\end{align*}
where $(\psi_1,\cdots,\psi_4)^T$ is the eigenfunction of $\mathcal{L}$ corresponding to the eigenvalue $\mu$. Thus, let

\begin{align*}
\begin{pmatrix}
\psi_1\\ \vdots\\ \psi_4
\end{pmatrix}=\sum_{k=0}^\infty
\begin{pmatrix}
{V_k}\\ h_k\\m_k\\ {n_k}
\end{pmatrix}\cos\left(\frac{\pi kx}{\ell}\right),
\end{align*}
where $V_k, h_k, m_k$ and $n_k$ are time-dependent coefficients. We can then conclude that

\begin{align*}
\sum_{k=0}^\infty \mathcal{J}_k
\begin{pmatrix}
{V_k}\\ h_k\\m_k\\ {n_k}
\end{pmatrix}=\mu\sum_{k=0}^\infty
\begin{pmatrix}
{V_k}\\ h_k\\m_k\\ {n_k}
\end{pmatrix},
\end{align*}
with 

\begin{align}
\mathcal{J}_k=\left(\scalemath{1}{\begin{array}{cccc} 
\mathcal{J}_{k_{11}}&\left.\displaystyle\frac{\partial \mathcal{F}}{\partial h}\right|_{X_\infty }&\left.\displaystyle\frac{\partial \mathcal{F}}{\partial m}\right|_{X_\infty }&\left.\displaystyle\frac{\partial \mathcal{F}}{\partial n}\right|_{X_\infty }\\[1.5ex]
\displaystyle\frac{1}{\tau_h}\displaystyle\frac{\partial h_\infty}{\partial V}&-\displaystyle\frac{1}{\tau_h}&0&0\\[1.5ex]
\displaystyle\frac{1}{\tau_m}\displaystyle\frac{\partial m_\infty}{\partial V}&0&-\displaystyle\frac{1}{\tau_m}&0\\[1.5ex]
\displaystyle\frac{1}{\tau_n}\displaystyle\frac{\partial n_\infty}{\partial V}&0&0&-\displaystyle\frac{1}{\tau_n}
\end{array}}\right),\label{J_k}
\end{align}
where 

$$
\mathcal{J}_{k_{11}}=\left.\displaystyle\frac{\partial \mathcal{F}}{\partial V}\right|_{X_\infty }-\left(\frac{\pi k}{\ell}\right)^2\frac{\lambda}{1+\lambda}\frac{M_i}{\chi}\frac{1}{C_m}\quad k=0,1,2,3,\dots.
$$
Keep in mind that system \eqref{spectral_problem} has eigenvalues

$$
0=\mu_0<\mu_1=\left(\frac{\pi }{\ell}\right)^2<\mu_2=4\left(\frac{\pi }{\ell}\right)^2<\mu_3=9\left(\frac{\pi }{\ell}\right)^2<\cdots\longrightarrow \infty.
$$
Hence, the linearised system \eqref{eq:pde-linear} has infinitely many Jacobians $\mathcal{J}_{k}$. 



We continue by deriving an ODE model to represent the behaviour of the linearised system \eqref{eq:pde-linear} for each mode $k=0,1,2,3,\dots$ in close proximity to the equilibrium $V_\infty$. This is done by ensuring that the resulting system has the same equilibrium as the Noble model~\eqref{model} and the Jacobian $\mathcal{J}_k$.


Then, we are in a situation where we can analyse the behaviour of a single cell on a macro-scale including the influence of the diffusion term (dependent on $k$) of the monodomain model~\eqref{monodomain}. This allows us to gain intuition on how the cells interact. The resulting ODE system reads as follows:


\begin{align}
\frac{\dd}{\dt}
\begin{pmatrix}
{V_k}\\ h_k\\m_k\\ {n_k}
\end{pmatrix}=
\begin{pmatrix}
-\displaystyle\frac{I_\text{Na}+I_\text{K}+I_\text{L}}{C_m}-\left(\frac{\pi k}{\ell}\right)^2\frac{\lambda}{1+\lambda}\frac{M_i}{\chi}\frac{(V_k-V_\infty)}{C_m}\\ (h_\infty(V_k)-h_k)/\tau_h(V_k)\\ (m_\infty(V_k)-m_k)/\tau_m(V_k)\\ (n_\infty(V_k)-n_k)/\tau_n(V_k)\end{pmatrix},\label{PDE-ODE-system}
\end{align}
where $V_\infty$ is an equilibrium for equation~\eqref{eq}. We have designed the location of the equilibrium of system~\eqref{PDE-ODE-system} to be the same as for the Noble model~\eqref{model}. However, considering the stability analysis from Section~\ref{sec:model}, the stability of the equilibrium may be different. The first entry of the Jacobian $\mathcal{J}$ in~\eqref{jac} is replaced by $\mathcal{J}_{k_{11}}$. Thus, also the coefficients $a_j$, $j=1,\dots,4$ of the characteristic polynomial~\eqref{charactersitic_equ} are changed. Obviously, this will affect the stability of the system~\eqref{PDE-ODE-system} dependent on the parameters $\lambda,~M_i,~\chi,~\ell$ and $k$. This indicates that the cellular behaviour of model~\eqref{model} is not (necessarily) one-to-one transferred to the behaviour and dynamics of the monodomain equation~\eqref{monodomain}. Nevertheless, it is a good starting point to study the dynamics of a single cell before extending the analysis to the macro-scale. Do however note that the mode $k=0$ does give us the same dynamics as that of the ODE system \eqref{model}.

In the discrete setting, if we choose $k=\frac{1}{\mathfrak{h}^2}$, where $\mathfrak{h}$ denotes the grid size, and $k=\frac{1}{\mathfrak{h}^2}$ is the biggest eigenvalue for the discrete laplacian, the term

$$
-(\pi k)^2\frac{\lambda}{1+\lambda}\frac{M_i}{\ell^2\chi}\frac{1}{C_m}=-\frac{\pi^2}{\mathfrak{h}^4}\frac{\lambda}{1+\lambda}\frac{M_i}{\ell^2\chi}\frac{1}{C_m}=-\frac{\pi^2}{\mathfrak{h}^4}\frac{1}{\ell^2}\frac{1}{4320}\frac{mScm^2}{\mu F}
$$
in \eqref{PDE-ODE-system} tends to $-\infty$ and blows up as $\mathfrak{h} \to 0$, which should stabilise \eqref{PDE-ODE-system} for large modes $k$ (or refined grid size $\mathfrak{h}$). Looking at the coefficients $a_j$, $j=1,\dots,4$ of the characteristic polynomial~\eqref{charactersitic_equ}, we can see that $a_1,$ $a_2$ and $a_3$ will tend to $\infty$, while $a_4$ will tend to $-\infty$ as $\mathfrak{h}$ tends to zero and $k$ tends to $\infty$. This implies that 

$$\Delta_1=a_1>0,~\Delta_2=a_1a_2-a_3>0,~\Delta_3=a_3\cdot \Delta_2>0~\text{and}~\Delta_4=\Delta_3-a_1^2a_4>0$$
and thus, the numerics of the PDE will stabilise for decreasing grid size, which is to be expected. 

The same analysis can be used on 2D domains $\Omega=[0,\ell]\times[0,\ell_*]$, $\ell,\ell_*>0$, by considering the slightly different term

$$
-(\pi k)^2\frac{\lambda}{1+\lambda}\frac{M_i}{\chi}\frac{1}{C_m}\left(\frac{1}{\ell^2}+\frac{1}{\ell_*^2}\right).
$$
A similar modification applies to a 3D cube or cuboid. For a more general geometry one can expect that the additional term deriving from the linearisation of the monodomain model~\eqref{monodomain} is more complicated. However, the general discussion above is expected to remain valid.


\subsection{\sus{A ventricular cardiac cell model}}\label{sec:bernus_descr} 
\noindent As mentioned, we first investigate the Noble model~(\ref{model}). Indeed, we will see that this model has some limitations. Therefore, we will also study \sus{a human ventricular cell model, which is more advanced due to the number of included ion currents.} The description of the system in~\cite{Bernus} for epicardial cells is similar to system~(\ref{model}), but it contains more ion currents and reads as follows:

\begin{align}
	C_m\frac{\dd V}{\dt}=-I_\mathrm{ion}+I_\mathrm{stimulus}, \label{model_bernus}
\end{align}
where $I_\mathrm{stimulus}$ denotes an external stimulus and 

$$I_\mathrm{ion}=I_\mathrm{Cab}+I_\mathrm{NaCa}+I_\mathrm{Nab}+I_\mathrm{Ca}+I_\mathrm{K}+I_\mathrm{NaK}+I_\mathrm{Na}+I_\mathrm{K1}+I_\mathrm{to}$$ is depending on the fast sodium current $I_\mathrm{Na}=G_\mathrm{Na}m^3v^2(V-E_\mathrm{Na}),$ the slow calcium current $I_\mathrm{Ca}=\frac{3}{5}G_\mathrm{Ca}d_\infty(V)f(V-E_\mathrm{Ca}),$ the transient outward current $I_\mathrm{to}=G_\mathrm{to}r_\infty(V)to(V-E_\mathrm{to}),$ the delayed rectifier $K$ current $I_\mathrm{K}$ and the inward rectifier ${K1}$ current $I_\mathrm{K1}=G_\mathrm{K1}{K1}_\infty(V)(V-E_\mathrm{K}),$ respectively. 

In~\cite{Bernus} the authors studied system \eqref{model_bernus} with a delayed rectifier current $I_\mathrm{K}=G_\mathrm{K}x^2(V-E_\mathrm{K})$, while we are considering the delayed rectifier current $I_\mathrm{K}=I_{\mathrm{K}_\mathrm{r}}+I_{\mathrm{K}_\mathrm{s}}$ from~\cite{PB}. Here, the current $I_{\mathrm{K}_\mathrm{r}}=G_{\mathrm{K}_\mathrm{r}}x_r\mathrm{rik}(V)(V-E_\mathrm{K})$ denotes the rapidly activating current, while $I_{\mathrm{K}_\mathrm{s}}=G_{\mathrm{K}_\mathrm{s}}x_s^2(V-E_\mathrm{K})$ is the slowly activating current. Including the fast and slow potassium current will makes the dynamics more realistic. 

Furthermore, system~(\ref{model_bernus}) contains the background currents $I_\mathrm{Cab}$ and $I_\mathrm{Nab}$, the $\mathrm{Na}^+/\mathrm{Ca}^{2+}$ exchanger current $I_\text{NaCa}$  and the $\mathrm{Na}^+/\mathrm{K}^+$ pump current $I_\text{NaK}$, cf. Figure~\ref{fig:cell}(a). Notice that the system is depending on 9 gating variables, i.e. $m$, $v$, $d$, $f$, $r$, $to$, $x_r$, $x_s$ and ${K1}$, satisfying the differential equation~(\ref{diff_gating_variable}), where $d$, $r$ and ${K1}$ are assumed to be equal to their steady states. We will consider all gating variables as state variables, therefore the dimension of the system is 10. 

Moreover, we use $C_\mathrm{m}=1.534\ \frac{\mu F}{cm^2}$, cf.~\cite{PB}, while the individual ionic conductances are given by $G_\mathrm{Na} = 16\ \frac{mS}{cm^2}$, $G_\mathrm{Ca} = 0.064\ \frac{mS}{cm^2}$, $G_\mathrm{to} = 0.3\ \frac{mS}{cm^2}$, $G_{\mathrm{K}_\mathrm{r}}= 0.015\ \frac{mS}{cm^2}$, $G_{\mathrm{K}_\mathrm{s}}= 0.02\ \frac{mS}{cm^2}$ and $G_\mathrm{K1} = 2.5\ \frac{mS}{cm^2}$. The equilibria and the Jacobian are similarly determined as for the Noble model~(\ref{model}). The only difference is that we have 6 ODEs more to consider. This increases the $4\times4$ matrix $\mathcal{J}$ to a $10\times 10$ matrix. 
Following the same approach as in Section~\ref{sec:monodomain}, we can derive from the monodomain model related to system~\eqref{model_bernus} with $I_\mathrm{stimulus}=40\frac{\mu A}{cm^2}$ and a duration of 2 seconds, i.e.

\begin{align}
\begin{split}
       C_m\frac{\partial V}{\partial t}  &= -\frac{\lambda}{1+\lambda}\frac{1}{\chi}\nabla\cdot\left(M_i\nabla V\right)-I_\mathrm{ion}+I_\mathrm{stimulus},\\
       \frac{\partial y}{\partial t}&=\frac{y_\infty(V)-y}{\tau_y(V)}, \quad y=m,v,d,f,r,to,x_r,x_s,K1
\end{split}\quad\text{in}~\Omega
\label{monodomain_bernus}
\end{align}
with Neumann boundary condition $\vec{\nu}\cdot\left(M_i\nabla V\right)=0$ on $\partial\Omega$, the following ODE system
\begin{align}
\frac{\dd}{\dt}
\begin{pmatrix}
V_k\\ m_k\\ v_k\\ d_k\\ f_k\\r_k\\to_k\\{x_r}_k\\{x_s}_k\\K1_k
\end{pmatrix}=
\begin{pmatrix}
-\displaystyle\frac{I_\mathrm{ion}-I_\mathrm{stimulus}}{C_m}-\left(\frac{\pi k}{\ell}\right)^2\frac{\lambda}{1+\lambda}\frac{M_i}{\chi}\frac{(V_k-V_\infty)}{C_m}\\ (m_\infty(V_k)-m_k)/\tau_m(V_k)\\ (v_\infty(V_k)-v_k)/\tau_v(V_k)\\ (d_\infty(V_k)-d_k)/\tau_d(V_k)\\ (f_\infty(V_k)-f_k)/\tau_f(V_k)\\ (r_\infty(V_k)-r_k)/\tau_r(V_k)\\ (to_\infty(V_k)-to_k)/\tau_{to}(V_k)\\ (x_{r_\infty}(V_k)-{x_r}_k)/\tau_{x_r}(V_k)\\ (x_{s_\infty}(V_k)-{x_s}_k)/\tau_{x_s}(V_k)\\ (K1_\infty(V_k)-{K1}_k)/\tau_{K1}(V_k)\end{pmatrix}\label{PDE-ODE-system_big}
\end{align}
where $V_\infty$ the equilibrium of the voltage $V$ of system~\eqref{model_bernus}. Note that the stability analysis for system~\eqref{PDE-ODE-system} and the previous discussion also holds true for system~\eqref{PDE-ODE-system_big}.

\subsection{Numerical methods}\label{sec:numericalmethods}
\noindent For our simulations we will use MATLAB R2019b and the ode solver \textit{ode15s} with a relative tolerance of $10^{-13}$ and an absolute tolerance of $10^{-18}$. For the monodomain models, the \textit{pdepe} solver is used. Moreover, as initial values for system~(\ref{model}) we will use $V_0=-79.04\ mV$, $h_0=0.81$, $m_0=0.045$ and $n_0=0.52$, while for the second model~(\ref{model_bernus}) we utilise $V_0=-93.3701\ mV$, $m_0=0.0004$, $v_0=0.9990$, $f_0=0.8797$, $x_{r_0}=0.0042$, $to_0=0.9999$, $d_0=0.0000$, $r_0=0.0000$, $K1_0=0.0419$ and $x_{s_0}=0.0912$, as long we do not specify anything else. The desired bifurcation diagrams will be derived utilising the MATLAB toolboxes MATCONT and CL\_MATCONT~\cite{Dhooge,Dhooge1,Govaerts}, which are numerical continuation packages for interactive bifurcation analysis of dynamical systems. 

\section{Dynamics of the Noble model}\label{sec:bif}
\begin{figure}[h]
\centering
\subfigure[Trajectory of system~(\ref{model}) with $G_\mathrm{L}=0.075\ \frac{mS}{cm^2}$.]{{\includegraphics[width=0.45\textwidth]{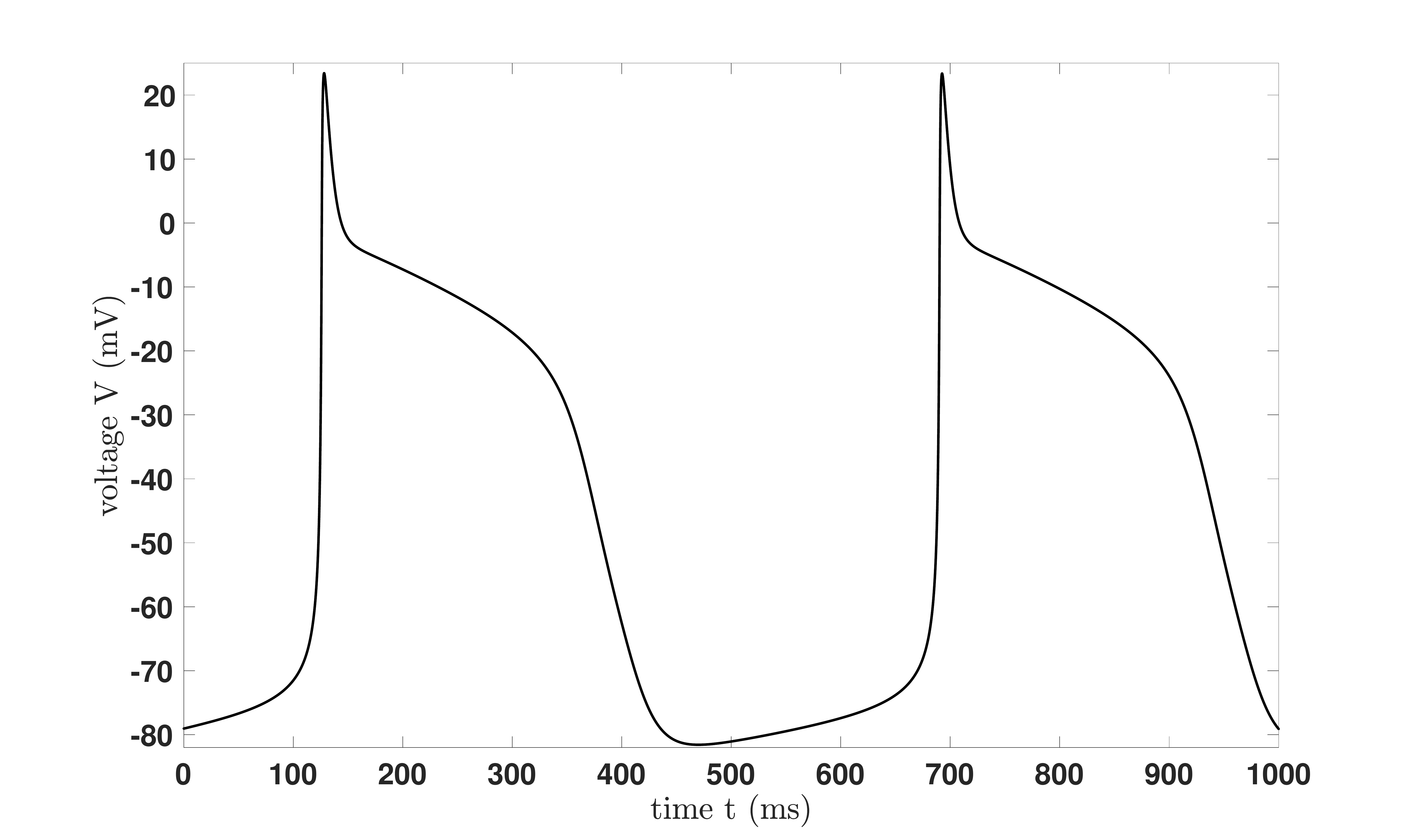}}}\subfigure[Trajectory of system~(\ref{model}) with $G_\mathrm{L}=0.18\ \frac{mS}{cm^2}$.]{{\includegraphics[width=0.45\textwidth]{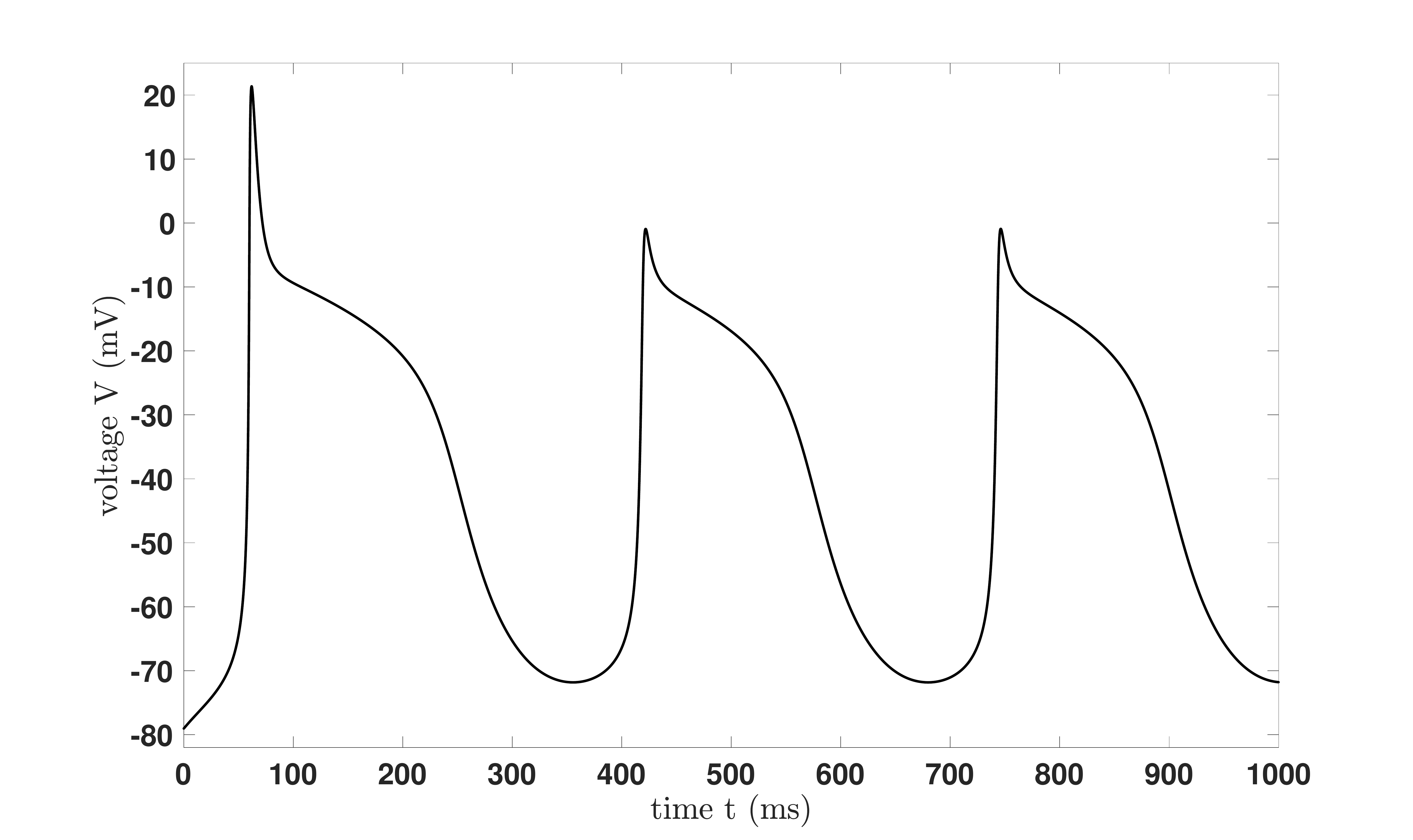}}}
\caption{Different action potentials: Simulation of system~(\ref{model}) (default setting) with two different values of the leak conductance $G_\mathrm{L}$, where these values are chosen according to the observation in~\cite{Noble}.} \label{fig:trajectory}
\end{figure}

\noindent In~\cite{Noble}, the author mentioned a change in the dynamics of system~(\ref{model}) by varying the leak conductance $G_\mathrm{L}$. We will show that changing $G_\mathrm{L}$ has influence on the period of the AP as well as if the system converges into a stable equilibrium or not. In Figure~\ref{fig:trajectory}, the trajectory of system~(\ref{model}) is given for two different values of the leak conductance $G_\mathrm{L}$, i.e. $G_\mathrm{L}=0.075\ \frac{mS}{cm^2}$ and $G_\mathrm{L}=0.18\ \frac{mS}{cm^2}$. In Figure~\ref{fig:trajectory}(a), system~(\ref{model}) reveals a periodic trajectory representing two APs of a cardiac single cell, while in Figure~\ref{fig:trajectory}(b) the trajectory needs certain amount of time to reach a stable periodic pattern also representing APs. This shows that system~(\ref{model}) is sensitive with respect to $G_\mathrm{L}$, which may have an influence on the amplitude and the period $T$ of the trajectory given by $V(t+T)-V(t)=0$. In addition, the initial value also has an influence on the appearing pattern (at least locally).

\subsection{Bifurcation analysis of system~\eqref{model} with respect to $G_\mathrm{L}$} \label{sec:bif_gl}
\noindent We now systematically investigate the behaviour described above utilising numerical bifurcation analysis as done in~\cite{AE_MMOs,Ae_control}. Our first step is to determine a bifurcation diagram with respect to the leak current, i.e. we choose the leak conductance $G_\mathrm{L}$ as the bifurcation parameter. Notice that we will only consider positive values of $G_\mathrm{L}$, since they are the physiologically relevant ones. Mathematically and numerically we can easily extend the bifurcation diagram also to negative values of $G_\mathrm{L}$. 

We start by determining the equilibrium curve, i.e. we calculate for different values of $G_\mathrm{L}$ the corresponding  equilibria, which gives us the desired equilibrium curve. Moreover, we determine the stability of each equilibria. The equilibrium curve we can easily calculate with a self-written routine. However, to derive a detailed enough bifurcation diagram efficiently, it is advisable to use a robust continuation algorithm as the one mentioned in Section~\ref{sec:numericalmethods}. 

\begin{figure}[h!]
\centering
{\includegraphics[width=0.9\textwidth]{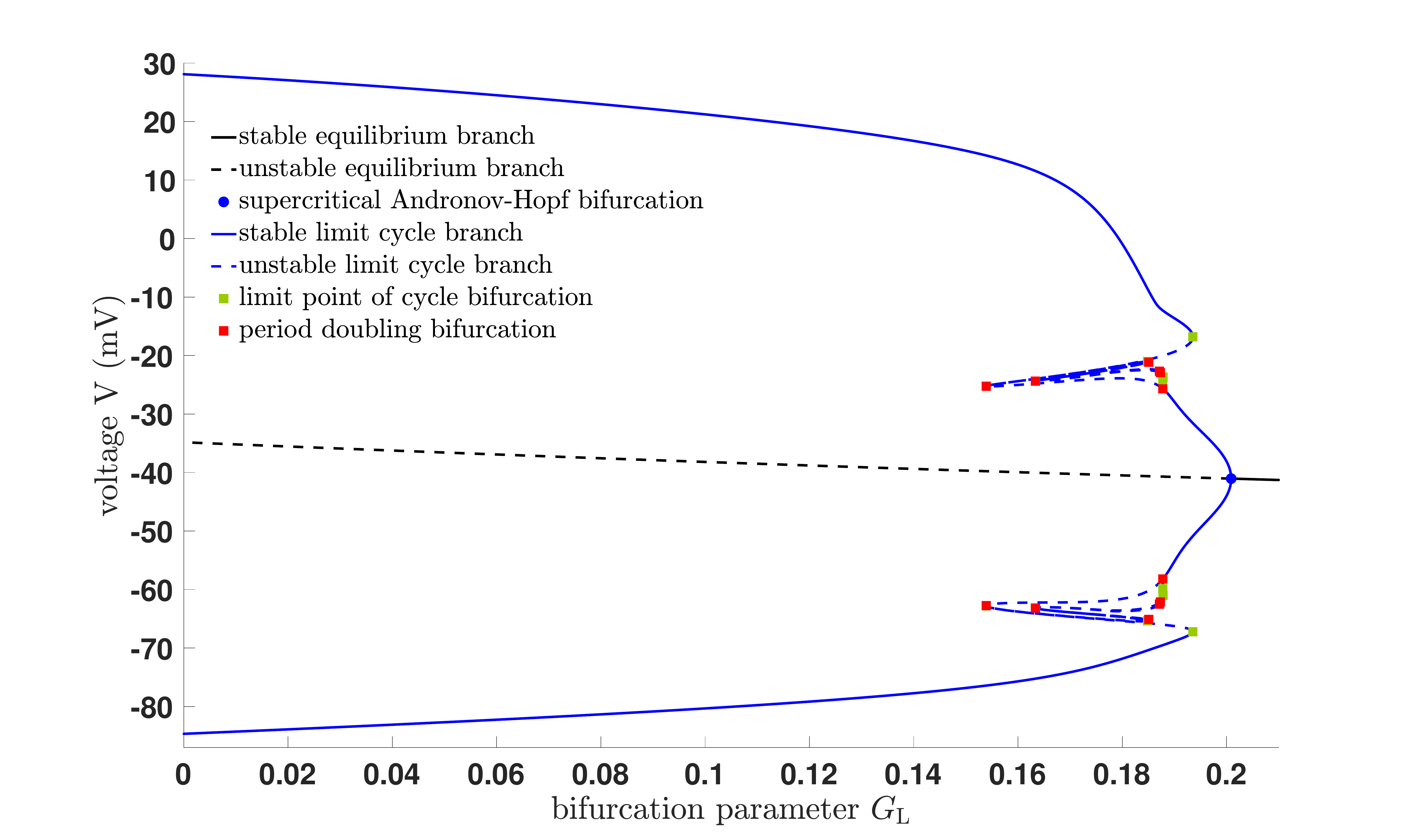}}
\caption{Bifurcation diagram in 2D: Projection onto the $(G_\mathrm{L},V)$-plane showing the equilibrium curve and the first two limit cycle branches.} \label{fig:bif_2D_GL}
\end{figure}
\begin{figure}[h]
\centering
\subfigure{{\includegraphics[width=0.5\textwidth]{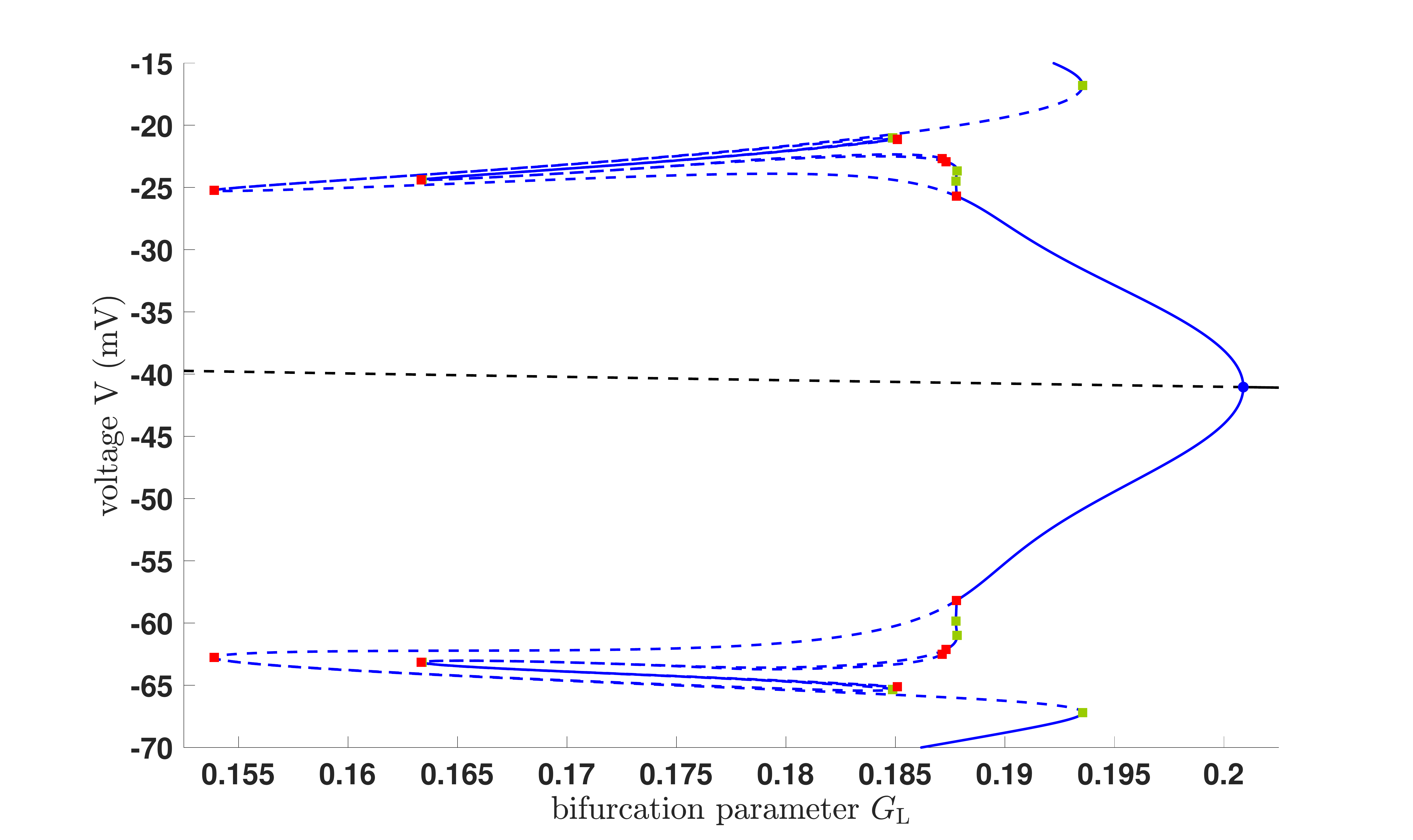}}}\subfigure{{\includegraphics[width=0.5\textwidth]{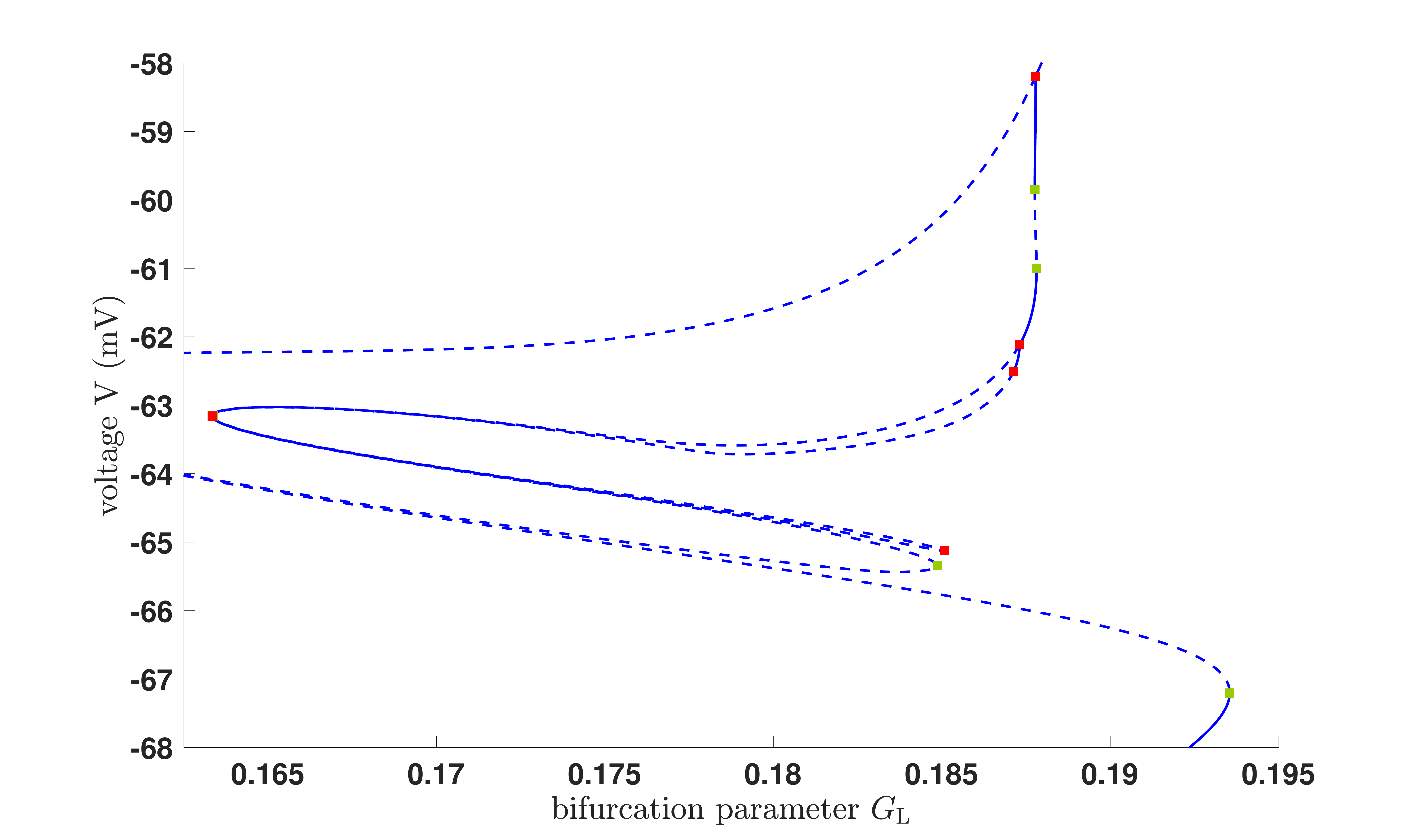}}}
\caption{Bifurcation diagram in 2D: Zoom of Figure~\ref{fig:bif_2D_GL} also including the third limit cycle branch.} \label{fig:bif_2D_GL_zoom}
\end{figure}

The bifurcation diagram of system~(\ref{model}) related to $G_\mathrm{L}$ exhibits an unstable (black dashed line) and a stable (black solid line) equilibrium branch, see Figure~\ref{fig:bif_2D_GL}. The equilibrium curve changes stability via a supercritical Andronov--Hopf bifurcation (blue dot, $G_\mathrm{L}\approx 0.200883 \frac{mS}{cm^2}$) with a negative first Lyapunov coefficient. From the supercritical Andronov--Hopf bifurcation a stable limit cycle branch (solid blue line) bifurcates, which becomes unstable (dashed blue line) via a period doubling bifurcation (solid red square, $G_\mathrm{L}\approx 0.187785 \frac{mS}{cm^2}$). Then, the limit cycle branch gains stability again via a limit point of cycle bifurcation (solid green square, $G_\mathrm{L}\approx0.193546 \frac{mS}{cm^2}$). Furthermore, this limit cycle branch contains a second period doubling bifurcation, which is also connected to the first one via a second limit cycle branch, cf. Figure~\ref{fig:bif_2D_GL} and Figure~\ref{fig:bif_2D_GL_zoom}.

The bifurcation diagram in Figure~\ref{fig:bif_2D_GL} only includes the first two limit cycle branches, as including further details would not be particularly visible. Indeed, further branches exist as we will see below. 

\begin{figure}[h]
\centering
{\includegraphics[width=0.9\textwidth]{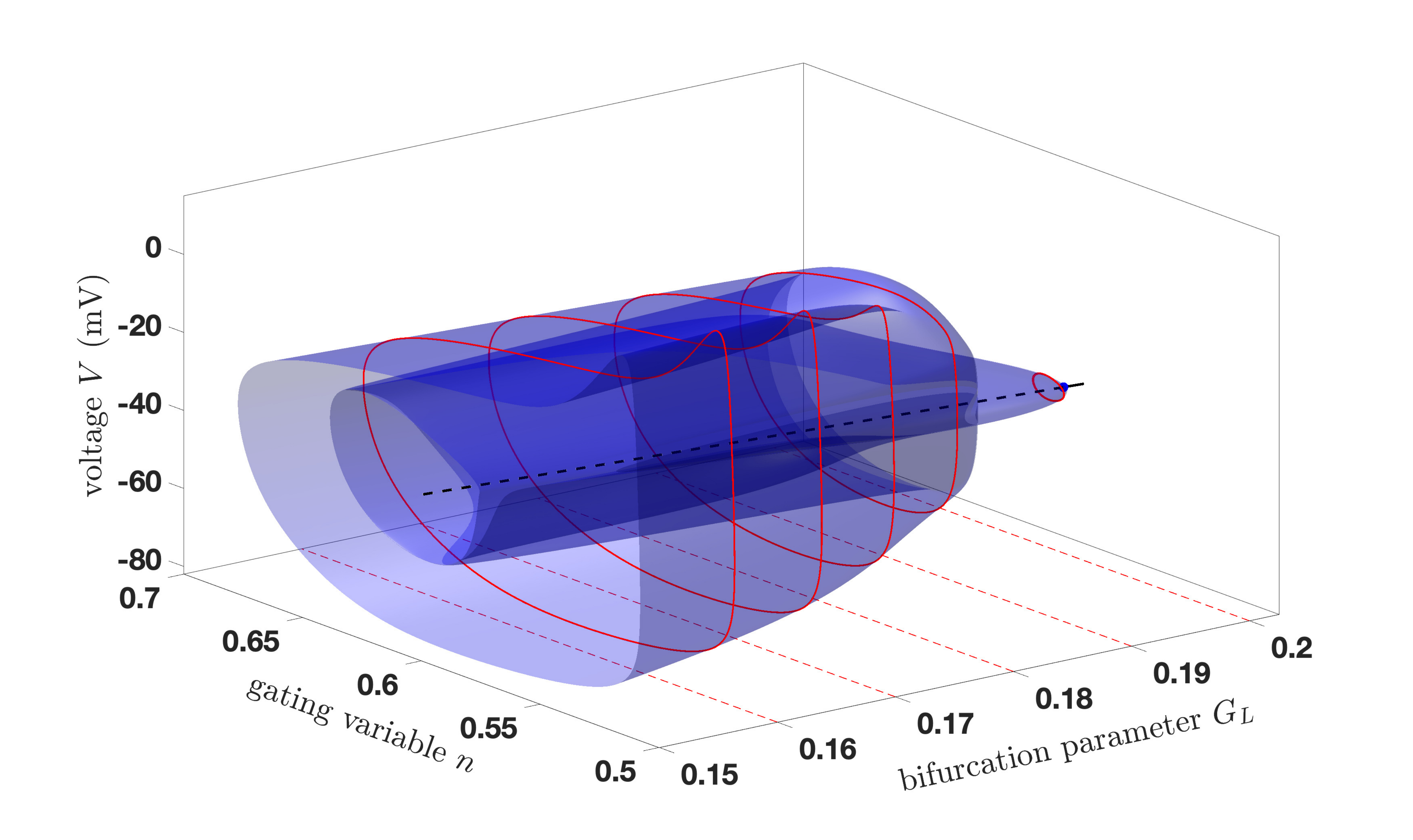}} 
\caption{Bifurcation diagram in 3D: Projection onto the $(G_\mathrm{L},n,V)$-plane including several trajectories of system~(\ref{model}) for different values of $G_\mathrm{L}$.} \label{fig:bif_3D_GL}
\end{figure}

Together with Figure~\ref{fig:bif_2D_GL_zoom} and Figure~\ref{fig:bif_3D_GL}, a nice graphical explanation of the observations in~\cite{Noble} appear. We can identify values of $G_\mathrm{L}$ for which system~(\ref{model}) oscillates or has a stable equilibrium. Even more, we can determine the maximal amplitude of an oscillation corresponding to $G_\mathrm{L}$. In addition, we also get the corresponding period of each limit cycle. We have that the period $T$ for $G_\mathrm{L}=0.075\ \frac{mS}{cm^2}$ is approximately $564.1345ms$, while for $G_\mathrm{L}=0\ \frac{mS}{cm^2}$ we have $T\approx 839.5015ms$ and for $G_\mathrm{L}=0.18\ \frac{mS}{cm^2}$ we have $T\approx 324.2749ms$. Thus, we know that the pattern of the trajectory repeats faster if $G_\mathrm{L}$ increases before it converges into a stable equilibrium for $G_\mathrm{L}$ too large, which explains the observations from~\cite{Noble}. Furthermore, from Figure~\ref{fig:bif_2D_GL} we observe that the maximal amplitude is decreasing while $G_\mathrm{L}$ increases. 

In Figure~\ref{fig:bif_2D_GL_zoom} a third limit cycle branch is also included and we focus on a smaller range of $G_\mathrm{L}$ values, where interesting dynamics may appear. The next additional limit cycle branches behave very similarly to the third one, i.e. they are bifurcating from a period doubling bifurcation, they lose stability via a further period doubling bifurcation, and stay unstable until they converge into a period doubling bifurcation of the previous limit cycle branch. Only the second limit cycle branch behaves slightly different --- it bifurcates from the first period doubling bifurcation ($G_\mathrm{L}\approx0.187785 \frac{mS}{cm^2}$), becomes unstable via a limit point of cycle bifurcation, gains stability again via a second limit point of cycle bifurcation and finally, loses stability via a further period doubling bifurcation ($G_\mathrm{L}\approx 0.187308 \frac{mS}{cm^2}$). Then, it stays unstable until it converges into a period doubling bifurcation of the first limit cycle branch, cf. Figure~\ref{fig:bif_2D_GL_zoom}. 

Starting from the second limit cycle branch, system~(\ref{model}) exhibits a stable period doubling cascade, which is usually a route to chaos, cf. e.g.~\cite{KBE}. However, in the bifurcation diagram in Figure~\ref{fig:bif_3D_GL}, where several trajectories (red lines) are highlighted for different values of $G_\mathrm{L}$, no irregular or chaotic behaviour nor additional oscillations can be seen. 

The reason is quite simple: The system exhibits a stable limit cycle branch from the first limit cycle branch, which ``surrounds'' the interesting dynamics of the Noble model~(\ref{model}), cf.~Figure~\ref{fig:bif_2D_GL} and Figure~\ref{fig:bif_2D_GL_zoom}. Hence, for the standard initial condition we do not get any sudden change in the dynamics. The trajectory will jump on to the ``outer'' stable limit cycle branch and stays there.
\begin{figure}[h]
\centering
\subfigure{{\includegraphics[width=0.5\textwidth]{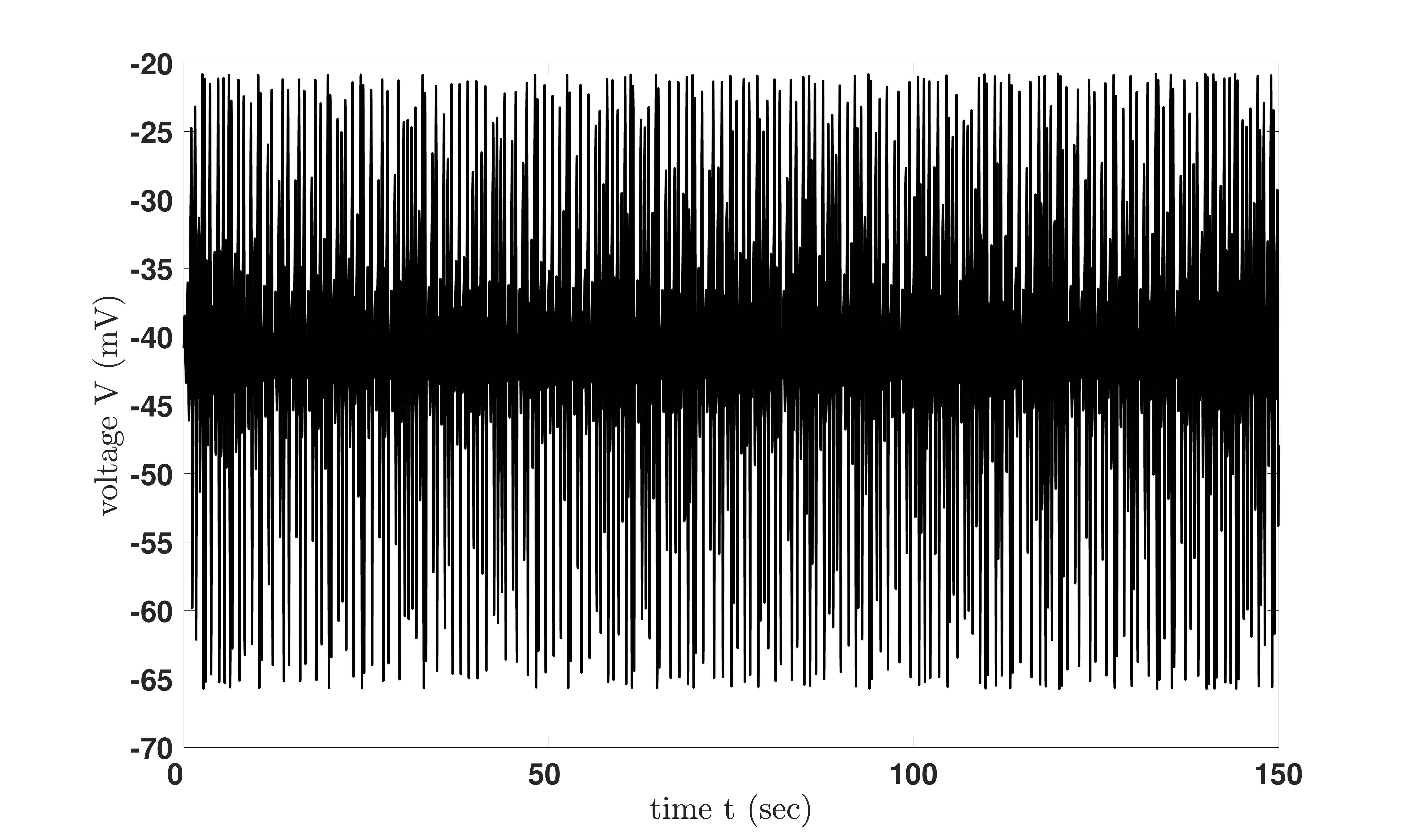}}}\subfigure{{\includegraphics[width=0.5\textwidth]{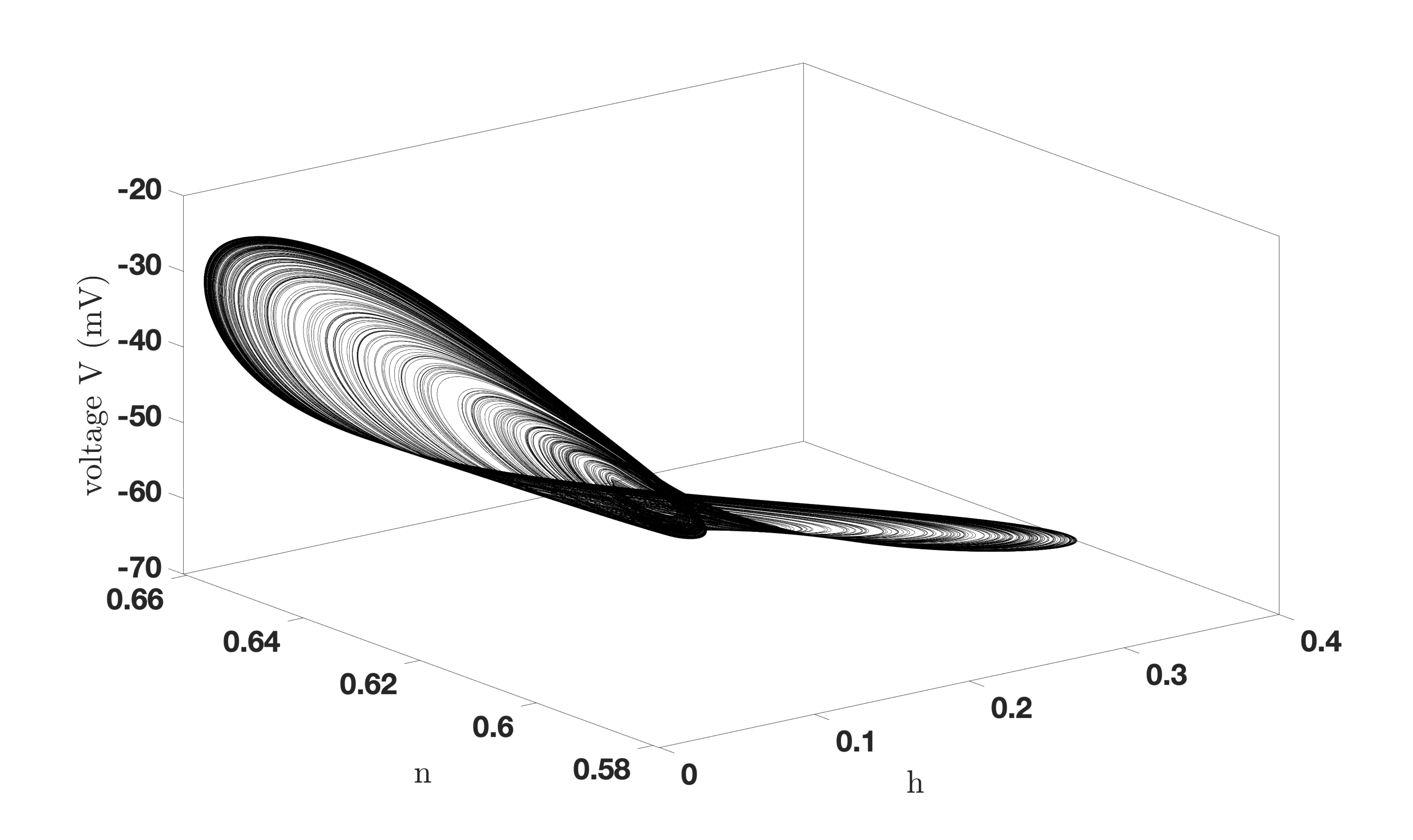}}}
\caption{Simulation of system~(\ref{model}) with $G_\mathrm{L}=0.1845\  \frac{mS}{cm^2}$ and initial values $V_0=-40.8454\ mV$, $h_0= 0.0268$, $m_0=0.3233$ and $n_0=0.5852$.} \label{fig:chaos}
\end{figure}

\begin{figure}[h]
\centering
\subfigure[$G_\mathrm{L}= 0.18778\  \frac{mS}{cm^2}$.]{{\includegraphics[width=0.4\textwidth]{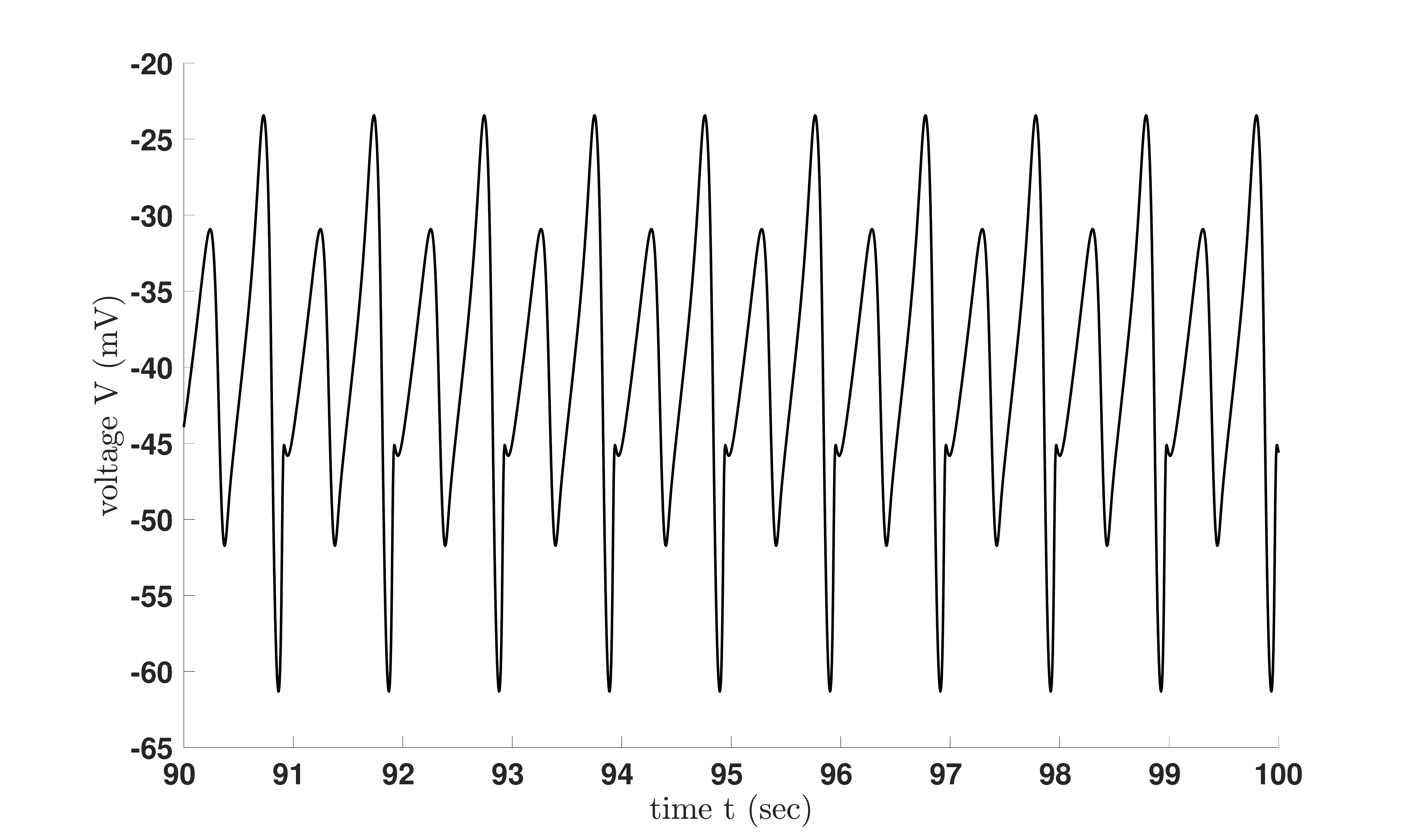}}}\subfigure[$G_\mathrm{L}= 0.18778\  \frac{mS}{cm^2}$.]{{\includegraphics[width=0.4\textwidth]{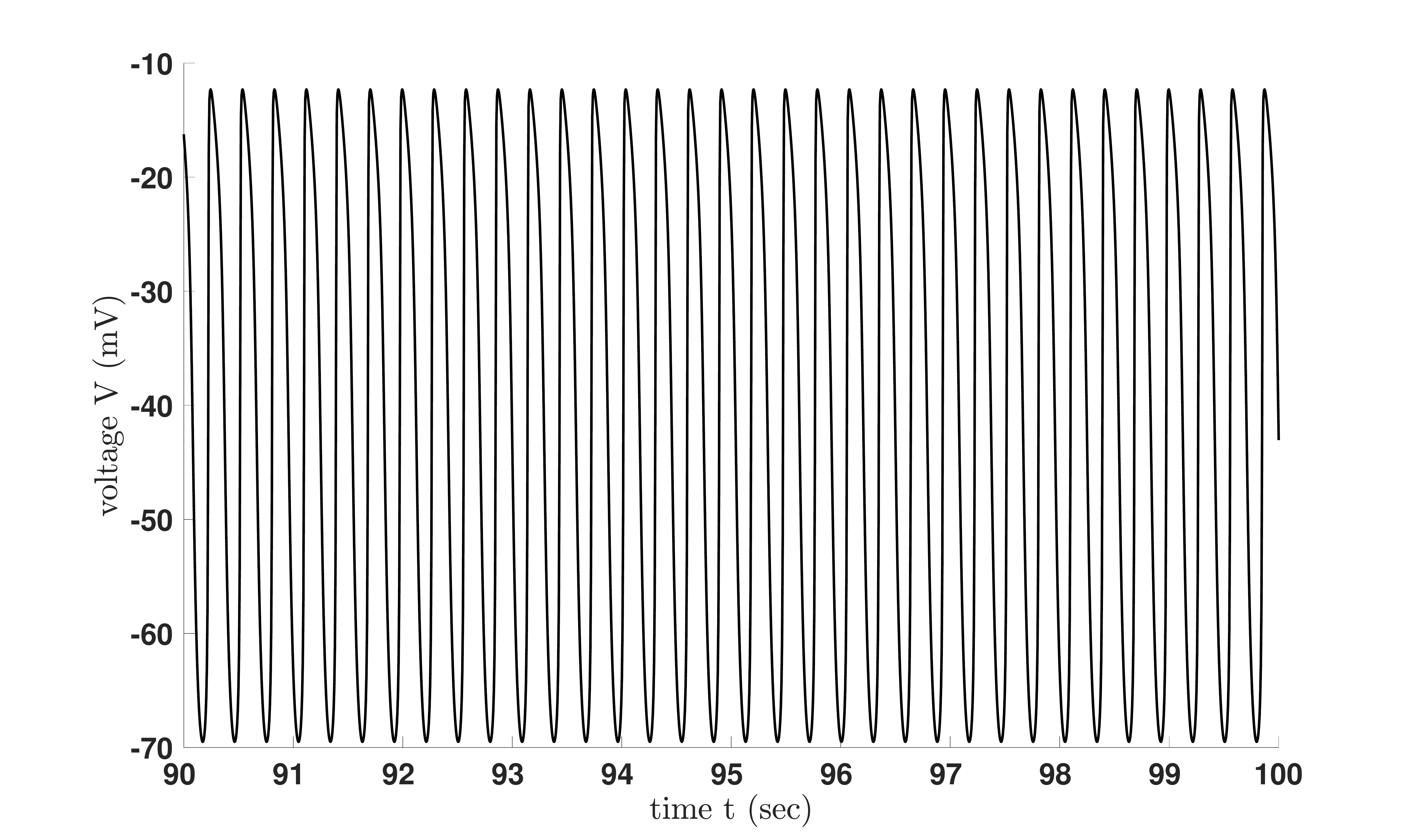}}}
\subfigure[$G_\mathrm{L}= 0.1865\  \frac{mS}{cm^2}$.]{{\includegraphics[width=0.4\textwidth]{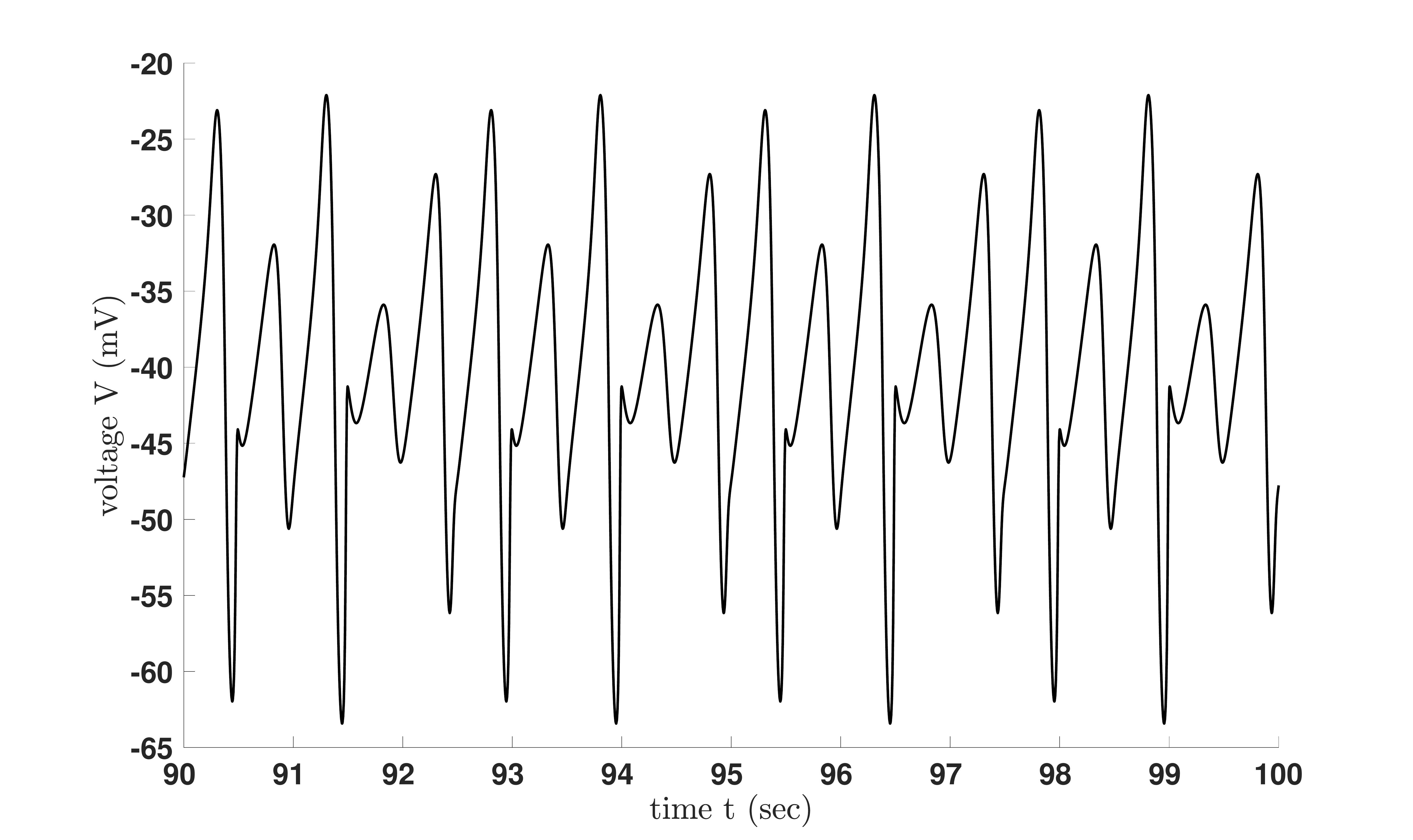}}}\subfigure[$G_\mathrm{L}= 0.1865\  \frac{mS}{cm^2}$.]{{\includegraphics[width=0.4\textwidth]{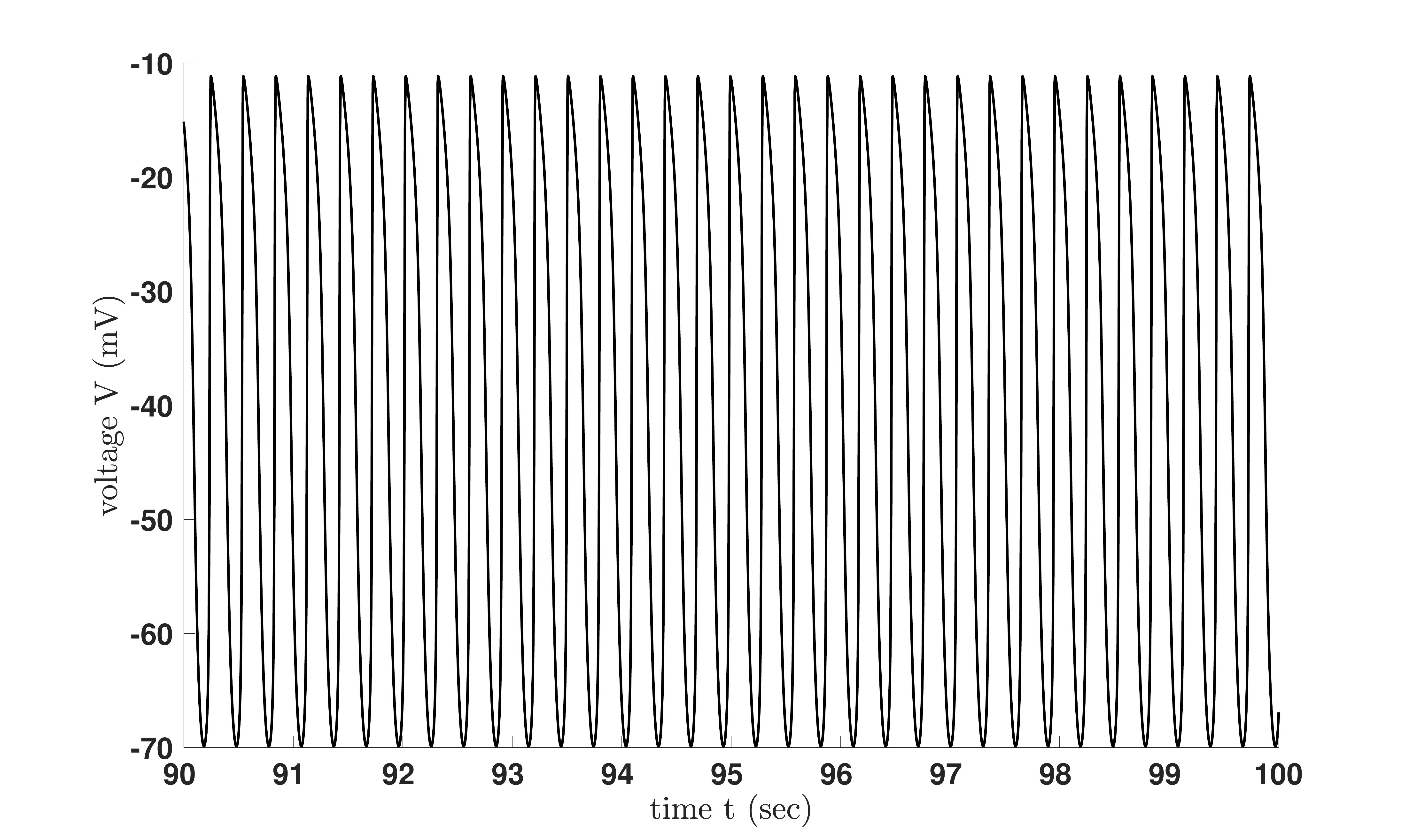}}}
\caption{Comparison of the effect of the choice of initial values on the trajectory of system~(\ref{model}): (a) shows the trajectory of model~(\ref{model}) at the first period doubling bifurcation with initial values $V_0=-40.8454\ mV$, $h_0= 0.0268$, $m_0=0.3233$ and $n_0=0.5852$, while in (b) the default initial values are used. (c) shows the trajectory of model~(\ref{model}) with $G_\mathrm{L}= 0.1865\  \frac{mS}{cm^2}$ and initial values $V_0=-40.8454\ mV$, $h_0= 0.0268$, $m_0=0.3233$ and $n_0=0.5852$, while in (d) the default initial values are used again.} \label{fig:comparison}
\end{figure}

However, if we choose an initial value from inside our limit cycle branches, e.g. one of the unstable equilibria, we will have a sudden change in the dynamics of the system close to the period doubling bifurcations and the period doubling cascade. For instance, Figure~\ref{fig:chaos} presents a simulation of system~(\ref{model}) over $150s$, which shows an irregular and chaotic behaviour. Hence, not only the choice of bifurcation parameter is essential, but also the choice of initial values, cf. Figure~\ref{fig:comparison}.
\subsection{Bifurcation analysis with respect to the potassium $I_\mathrm{K}$ current}
\noindent Next, we analyse the dynamics of system~(\ref{model}) with respect to the potassium $I_\mathrm{K}$ current. To this end, we choose one of the two potassium conductances $G_{\mathrm{K}_1}$ and $G_{\mathrm{K}_2}$. Notice that usually the $I_\mathrm{K}$ current depends on one gating variable, cf.~\cite{Xie,Bernus,LR} or one has a splitting into a fast $I_{\mathrm{K}_r}$ and a slow $I_{\mathrm{K}_s}$ current, cf.~\cite{TNNP04,TP06,PB}. Here the situation is different, since $I_\mathrm{K}$ is modelled as the sum of two currents, which is not split into a fast $I_{\mathrm{K}_r}$ and slow $I_{\mathrm{K}_s}$ part, nor into a potassium current and a background current, cf.~\cite{Bernus,TNNP04,PB}. Nevertheless, we will study the influence of a deficit in the potassium current in \eqref{model}. 

\begin{figure}[h]
\centering
\subfigure{{\includegraphics[width=0.5\textwidth]{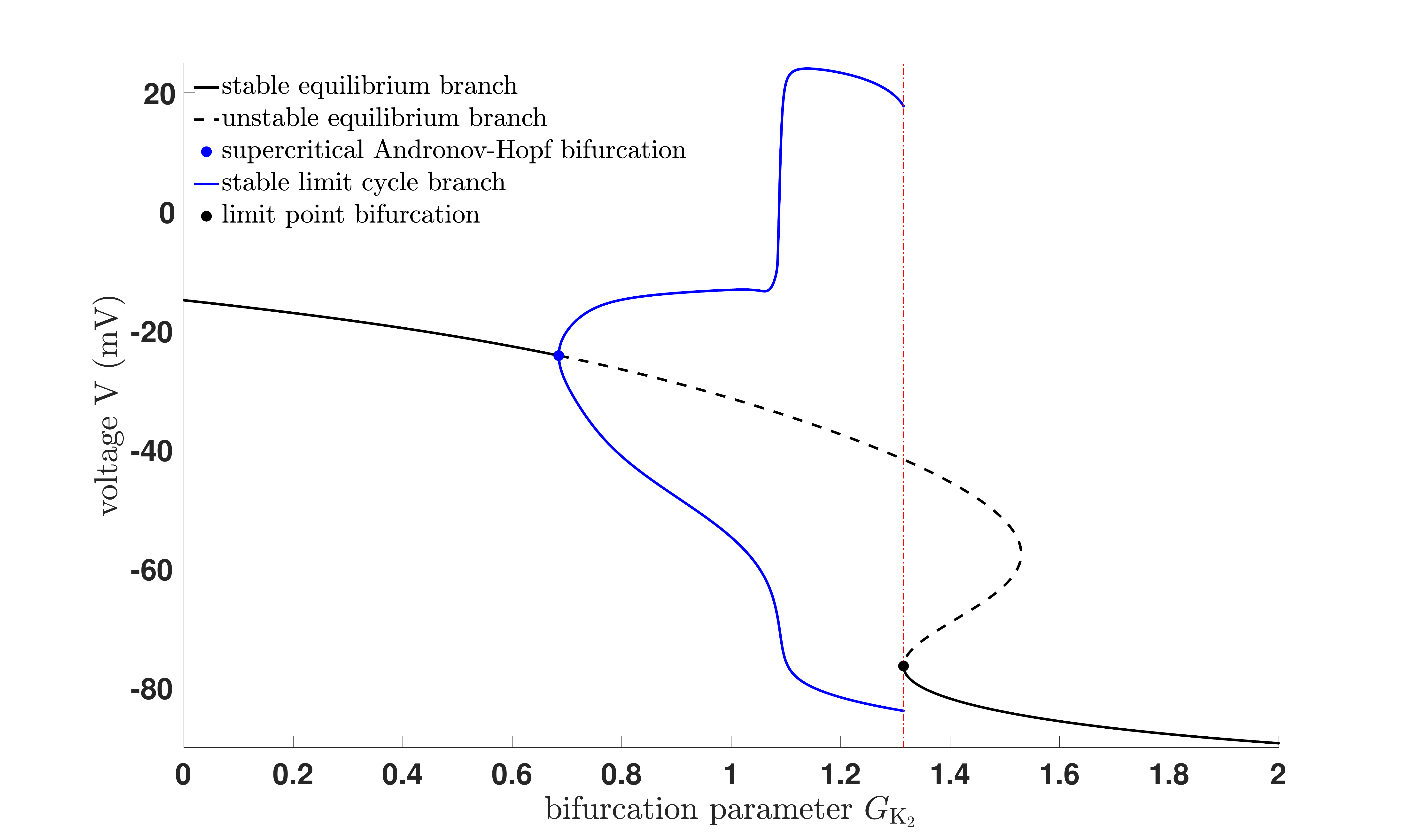}}}\subfigure{{\includegraphics[width=0.5\textwidth]{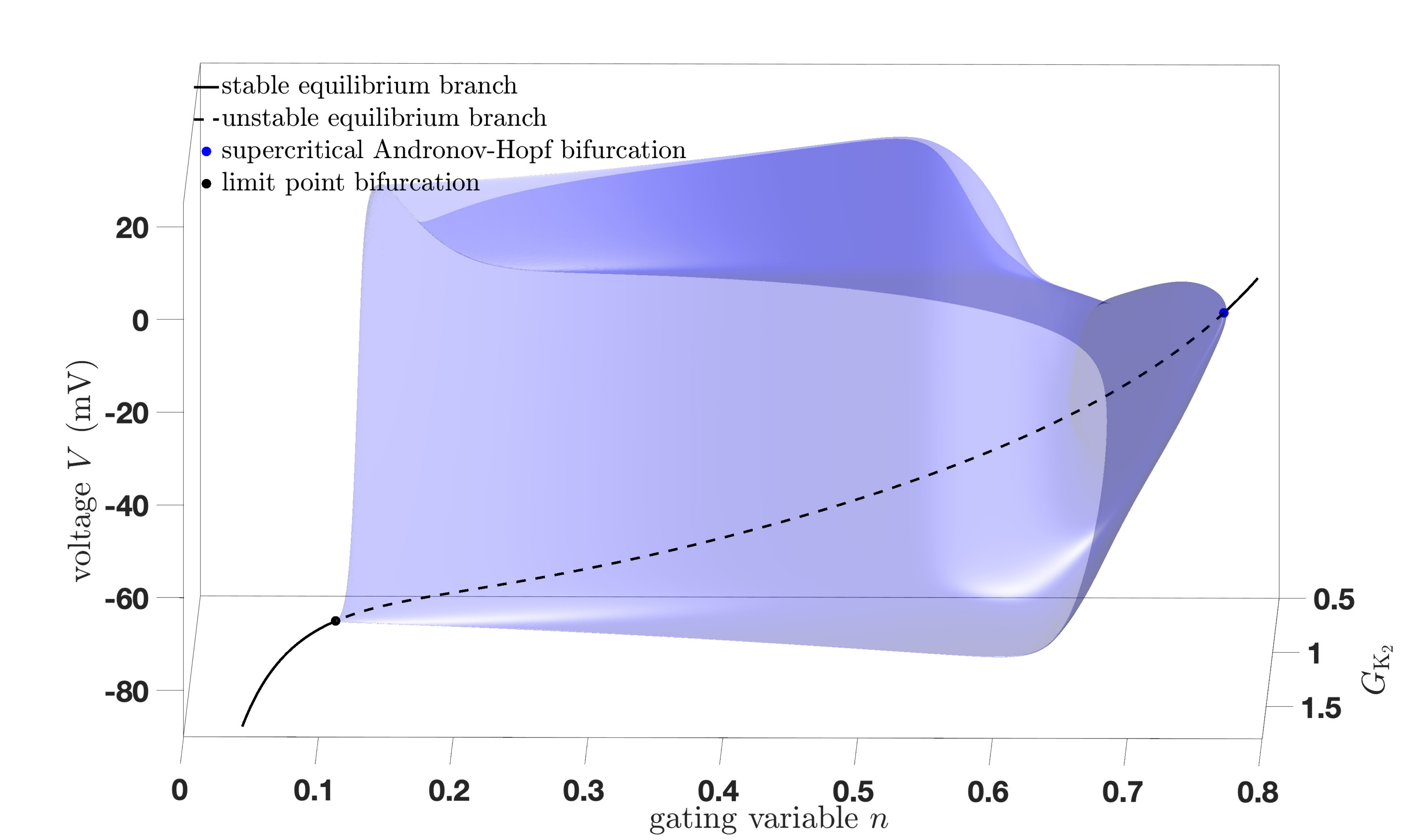}}}
\caption{Bifurcation diagram with $G_{\mathrm{K}_2}$ chosen as bifurcation parameter: From a supercritical Andronov--Hopf bifurcation $G_{\mathrm{K}_2}\approx0.6851\ \frac{mS}{cm^2}$ a stable limit cycle branch bifurcates and determine at a limit point bifurcation $G_{\mathrm{K}_2}\approx1.3147\ \frac{mS}{cm^2}$.} \label{fig:GK2}
\end{figure}

It is well known that a deficit in the potassium current may induce EADs. This was among others verified for simplistic models in~\cite{Tran,Xie}. Using bifurcation analysis we can conclude from Figure~\ref{fig:GK2} that no EADs appear via a deficit in $G_{\mathrm{K}_2}$. In the case where we combine $G_{\mathrm{K}_1}$ and $G_{\mathrm{K}_2}$, i.e. consider only one potassium conductance $G_\mathrm{K}=G_{\mathrm{K}_1}=G_{\mathrm{K}_2}$, system~\eqref{model} behaves similarly.

\begin{figure}[h]
\centering
\subfigure{{\includegraphics[width=0.5\textwidth]{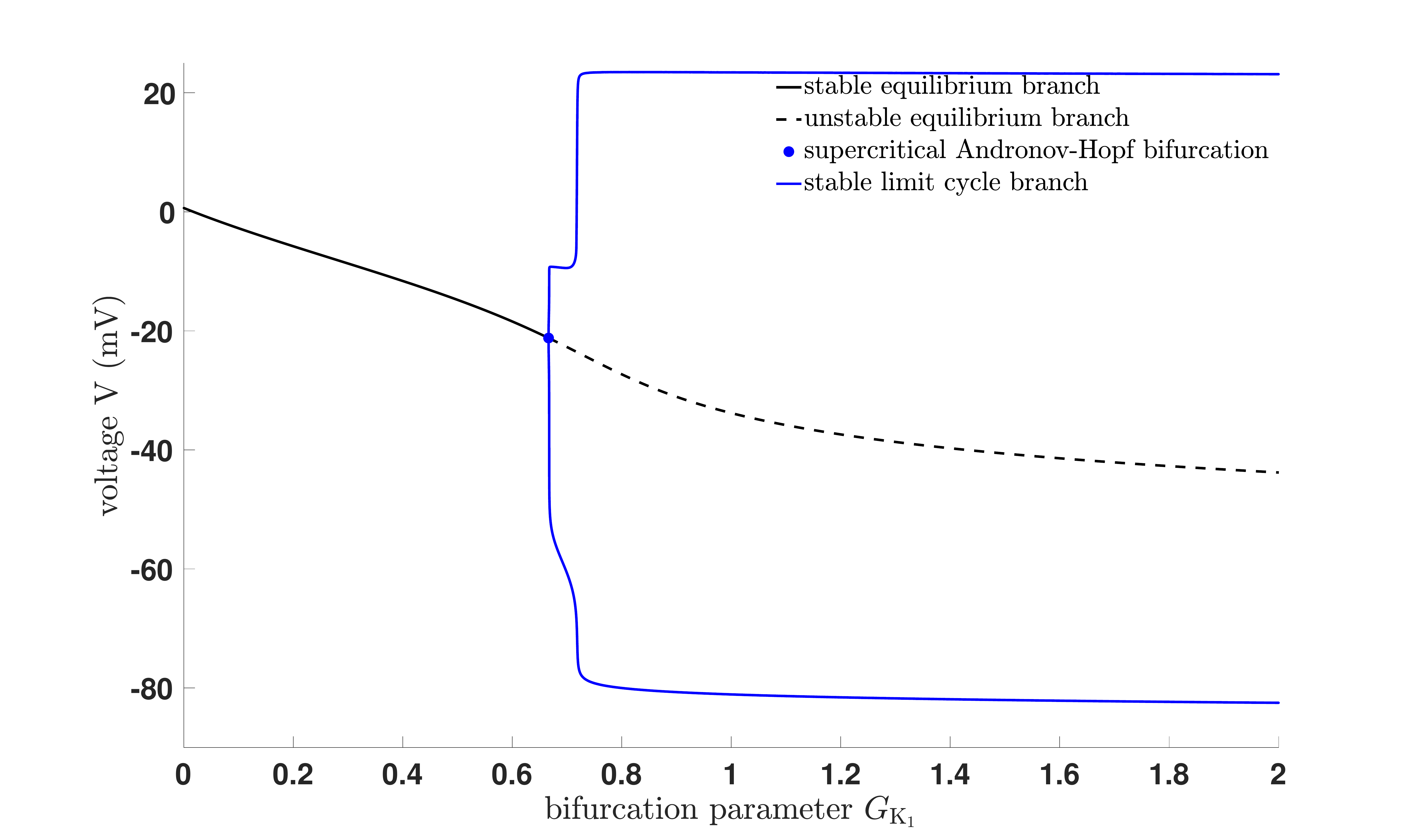}}}\subfigure{{\includegraphics[width=0.5\textwidth]{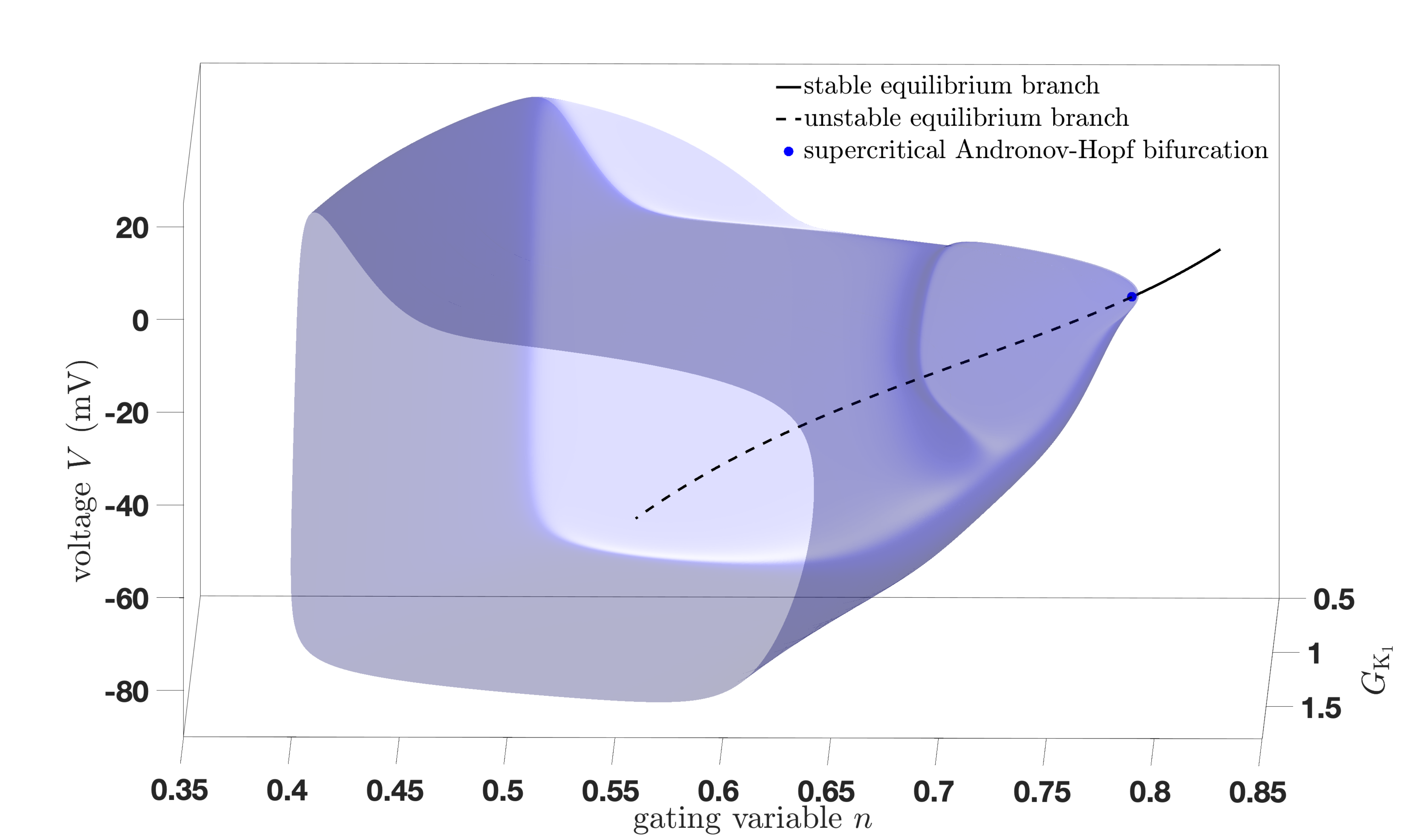}}}
\caption{Bifurcation diagram with $G_{\mathrm{K}_1}$ chosen as bifurcation parameter: From a supercritical Andronov--Hopf bifurcation $G_{\mathrm{K}_1}\approx0.6664\ \frac{mS}{cm^2}$ a stable limit cycle branch bifurcates.} \label{fig:GK1}
\end{figure}

Figure~\ref{fig:GK1} shows that the behaviour of system~\eqref{model} does not change dramatically under variations of $G_{\mathrm{K}_1}$ compared to variations of $G_{\mathrm{K}_2}$. We see that the model~(\ref{model}) does not exhibit EADs via a reduced potassium current. This was also checked for $G_{\mathrm{L}}$ between $0\ \frac{mS}{cm^2}$ and $0.2\ \frac{mS}{cm^2}$.

One explanation is that the calcium current $I_\mathrm{Ca}$ has an important influence on the behaviour of a cardiac cell, cf.~\cite{Ae_control,AE_MMOs}. This current is missing in the Noble model~\eqref{model}. For each additional current one gets additional ion current interactions, and one gains additional system parameters influencing the dynamics of the system. From~\cite{Ae_control} it is known that the ion current interaction between the potassium and the calcium current is important for the occurrence of EADs. This might be a reason for the appearance of different behaviour from what one would expect in a cardiac muscle cell. 
\\[2ex]
In conclusion, we see that we need a more detailed model to study more diverse behaviour (including EADs), see e.g.~\cite{Rodriguez_Purkinje}. The above investigation shows that Hodgkin--Huxley (type) models are sensitive with respect to their system parameters, and also with respect to their initial values. Thus, to capture all dynamics it is essential to consider both the system parameters and the initial values. The challenge is to derive a model which exhibits all dynamics of interest, and which at the same time is simple enough to allow its behaviour to be studied efficiently. This study indicates that model~\eqref{model} exhibits physiologically relevant APs or very fast oscillations with small amplitudes with respect to the potassium current $I_\mathrm{K}$. The system also exhibits chaotic behaviour with respect to the leak current $I_\mathrm{L}$. However, the voltage is always negative, see Figure \ref{fig:chaos}, and this seems to be non-physiological \cite{QU_review}. 

In Section \ref{sec:Bernus_bif} we will see that by slightly increasing the complexity of the model considered, more of the expected dynamics appear. 

\subsection{Effects on the macro-scale}
\noindent {In Section \ref{sec:monodomain} the analysis of the linearised system~\eqref{PDE-ODE-system} showed that one cannot expect that the cellular behaviour of a single cell model is one-to-one transferred to the behaviour and dynamics of the corresponding monodomain equation. Therefore, we briefly visualise how the interaction of an ensemble of cells belonging to two different regimes might play out in this section. Based on the discussion in the last two paragraphs of Section \ref{sec:monodomain}, we focus on a $1 \, cm$ one-dimensional cable, i.e. $x \in [0,1]$, for simplicity. }


To see the additional effects of the cell interactions at the macro-scale, we set the initial condition to partly belong to the chaotic regime of the Noble ODE model \eqref{model}, and partly to the stable regime:

\begin{align}
\left[V_0(x), h_0(x), m_0(x), n_0(x)\right] = 
\begin{cases}
[-40.8454\ mV, 0.0268, 0.3233, 0.5852], \quad & x \in \mathcal{D}\\
[-79.04\ mV, 0.81, 0.045, 0.52], \quad & x \in [0,1]\setminus \mathcal{D},
\end{cases}
\label{eq:macro_init}
\end{align}
where $\mathcal{D} \in [0,1]$. Note that this is not an equilibrium of the monodomain model \eqref{monodomain}. We fix $G_\mathrm{L}=0.1845\,\frac{mS}{cm^2}$ to allow for both chaotic and stable behaviour.  

\begin{figure}[h]
\centering
\subfigure{{\includegraphics[width=0.37\textwidth]{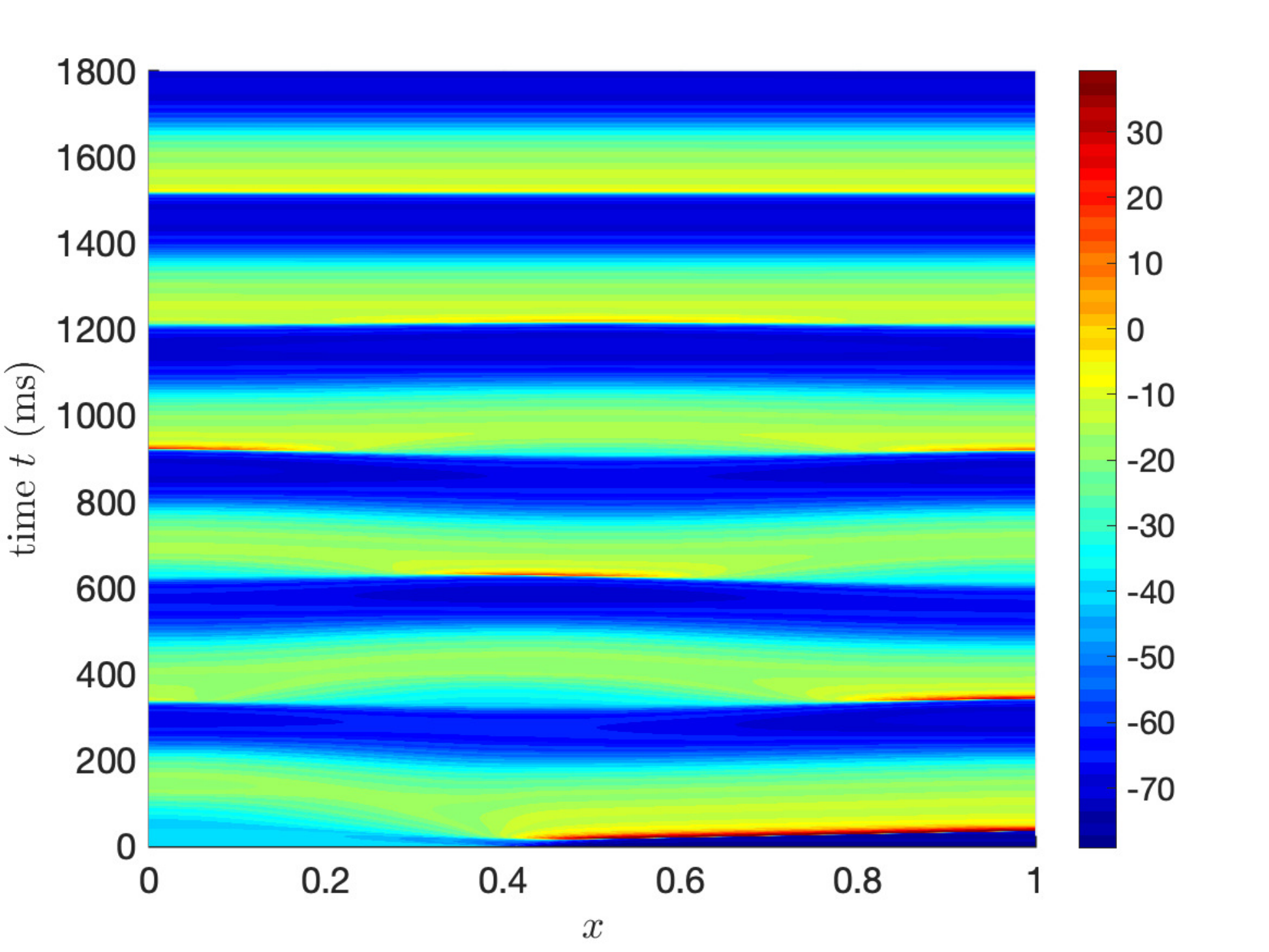}}}
\subfigure{{\includegraphics[width=0.325\textwidth]{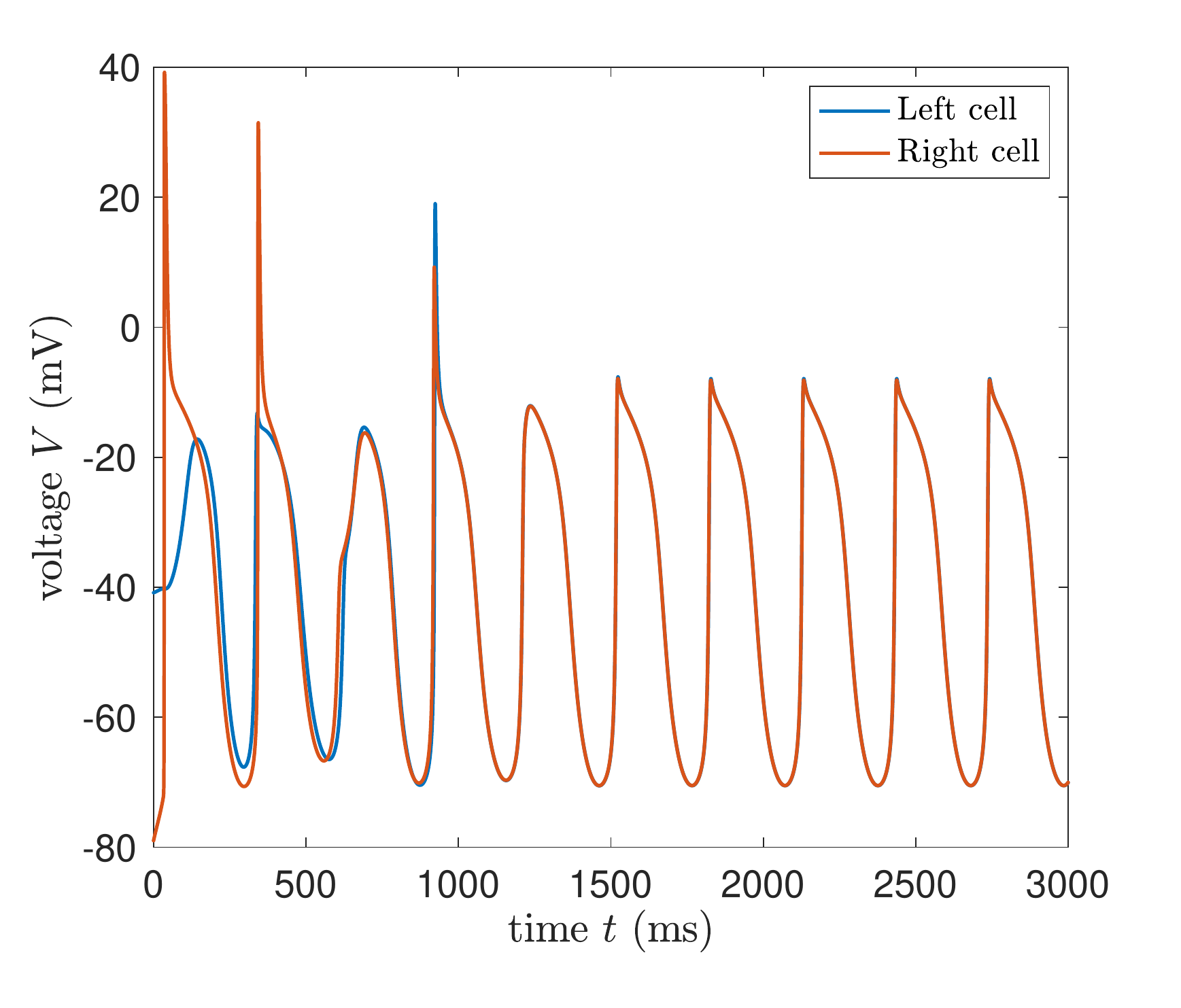}}}
\caption{Simulation of system~(\ref{monodomain}) with $G_\mathrm{L}=0.1845\,\frac{mS}{cm^2}$ and initial condition \eqref{eq:macro_init} with $\mathcal{D} = [0,\sqrt{2}-1)$. Computed with $256$ cells.} \label{fig:mono-stable}
\end{figure}

\begin{figure}[h]
\centering
\subfigure{{\includegraphics[width=0.325\textwidth]{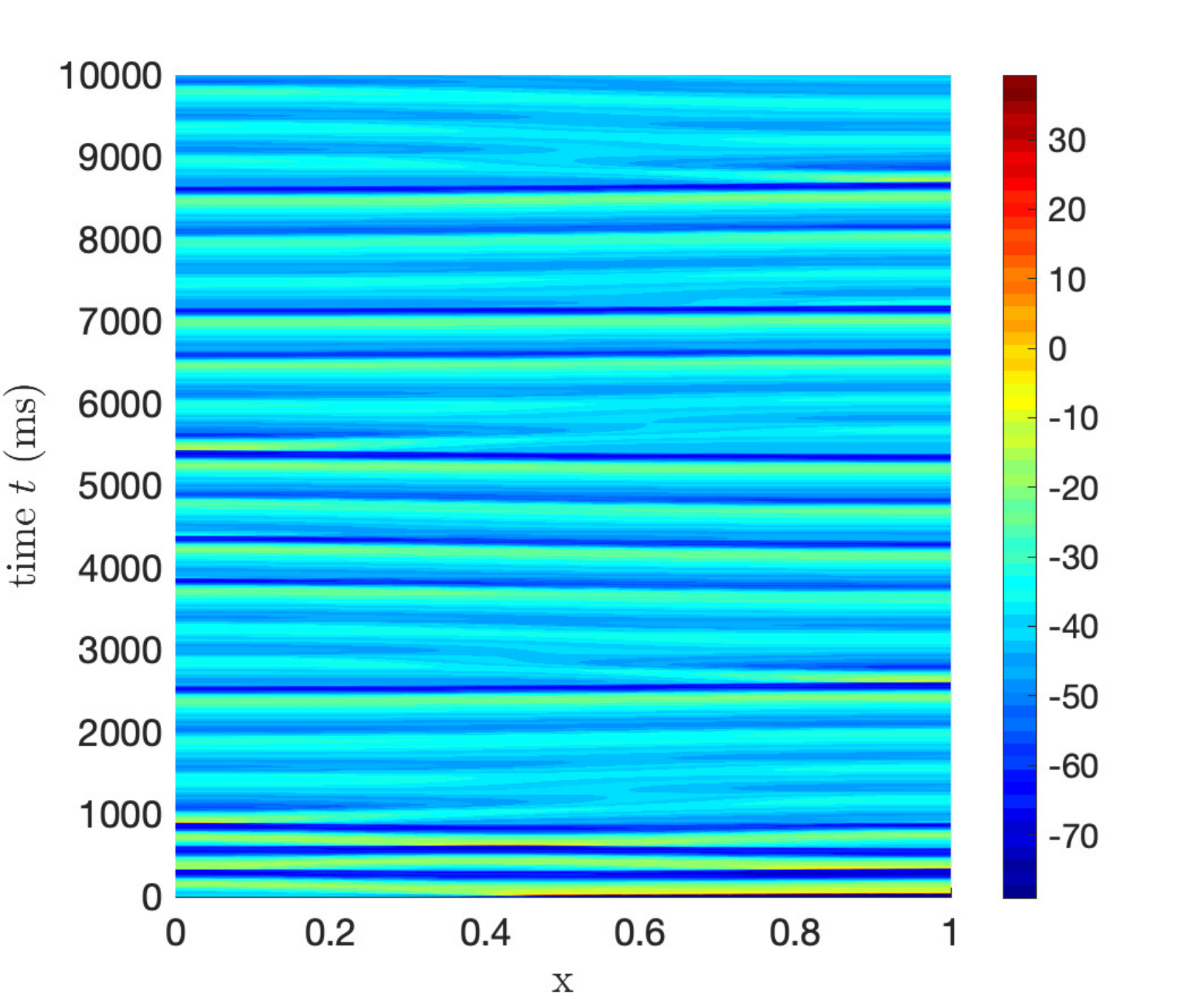}}}
\subfigure{{\includegraphics[width=0.325\textwidth]{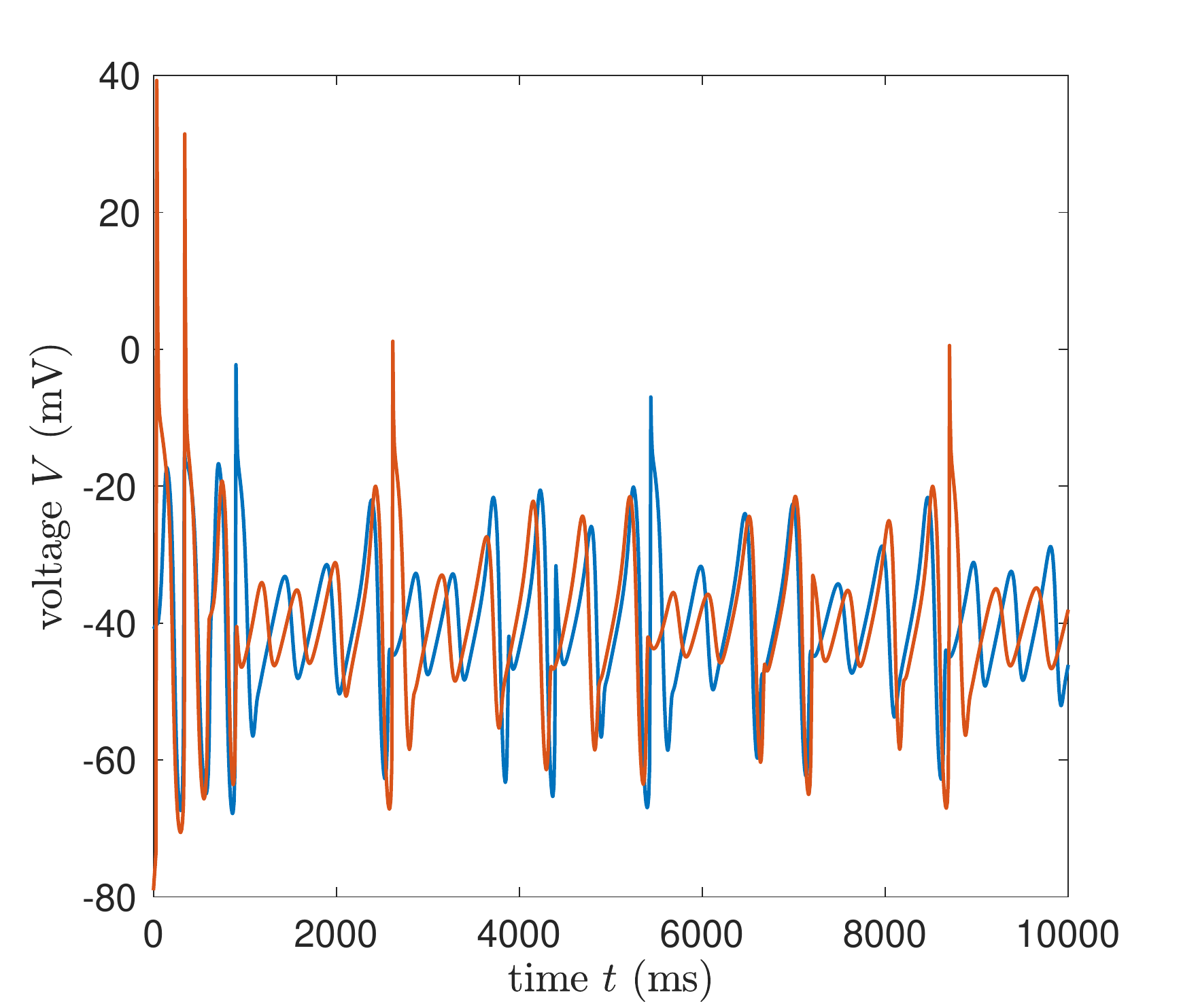}}}
\subfigure{{\includegraphics[width=0.325\textwidth]{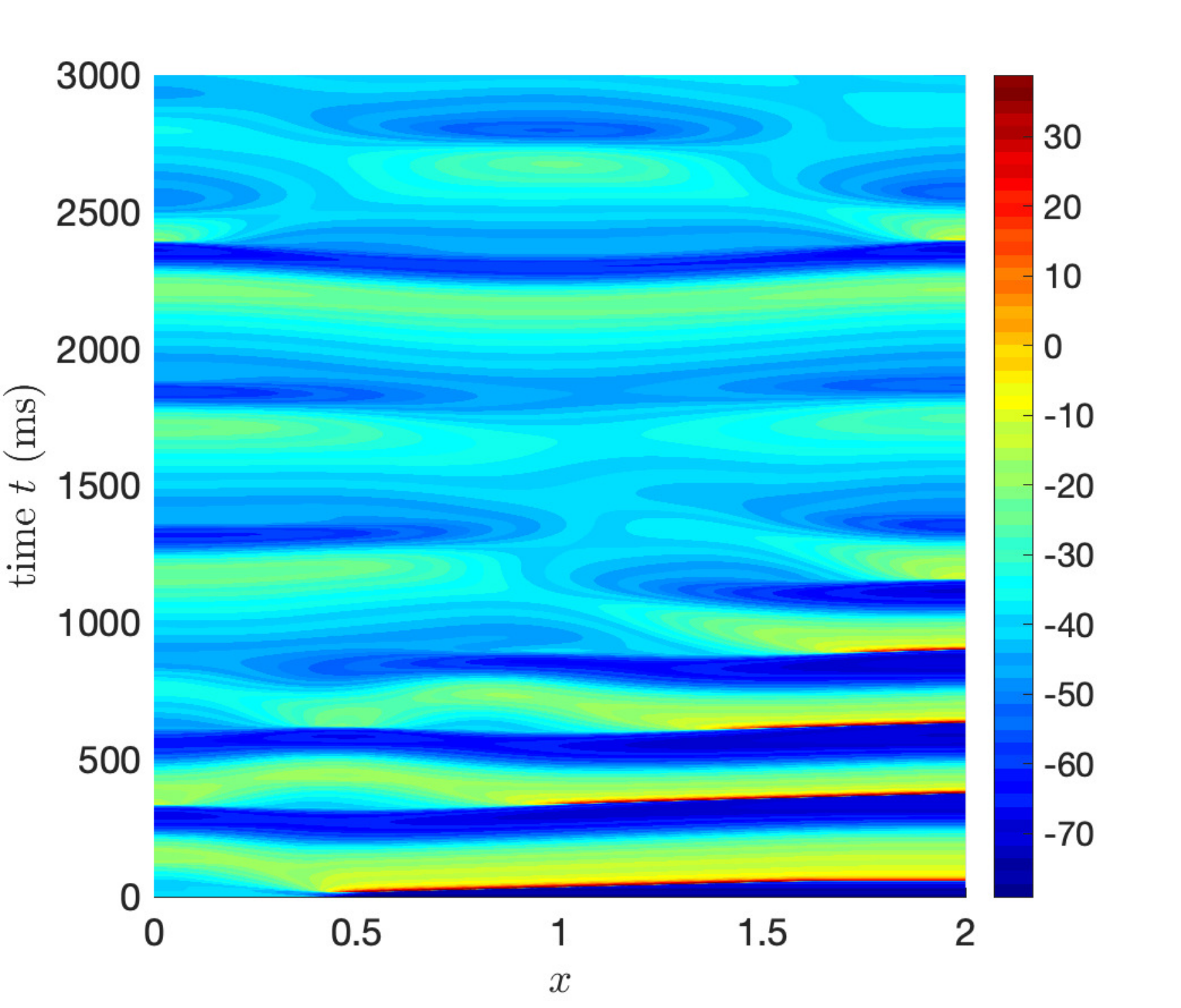}}}
\caption{Simulation of system~(\ref{monodomain}) with $G_\mathrm{L}=0.1845\,\frac{mS}{cm^2}$ and initial condition \eqref{eq:macro_init} with $\mathcal{D} = [0,\sqrt{2}-1)$. Computed with $256$ cells. The rightmost plot is computed on $[0,2]$ to better illustrate the dynamics as chaos takes over (all other parameters are the same).} \label{fig:mono-chaos}
\end{figure}

Letting $\mathcal{D} = [0,\sqrt{2}-1)$, Figure \ref{fig:mono-stable} shows that the stable behaviour suppresses the chaotic one and a travelling wave dynamic appear as time progresses. Increasing the initial chaotic domain slightly to $\mathcal{D} = [0,\sqrt{2}-0.99)$, one can observe from Figure \ref{fig:mono-chaos} that the chaotic behaviour prevails.

\begin{figure}[h]
\centering
\subfigure[Diffusion parameter $\frac{1}{360}$\,\it{mS}.]{{\includegraphics[width=0.33\textwidth]{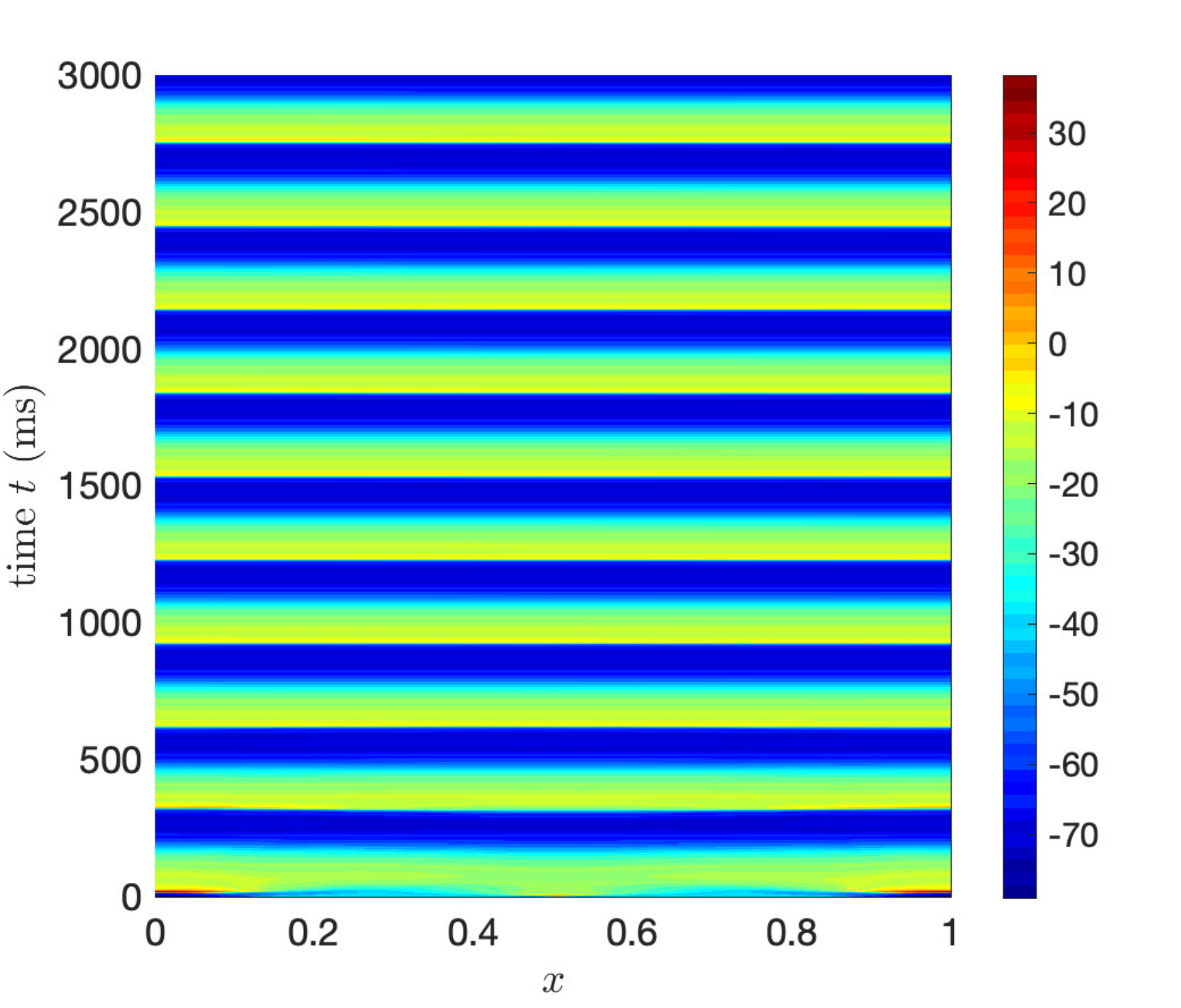}}}
\subfigure[Diffusion parameter $0.00005$\,\it{mS}.]{{\includegraphics[width=0.33\textwidth]{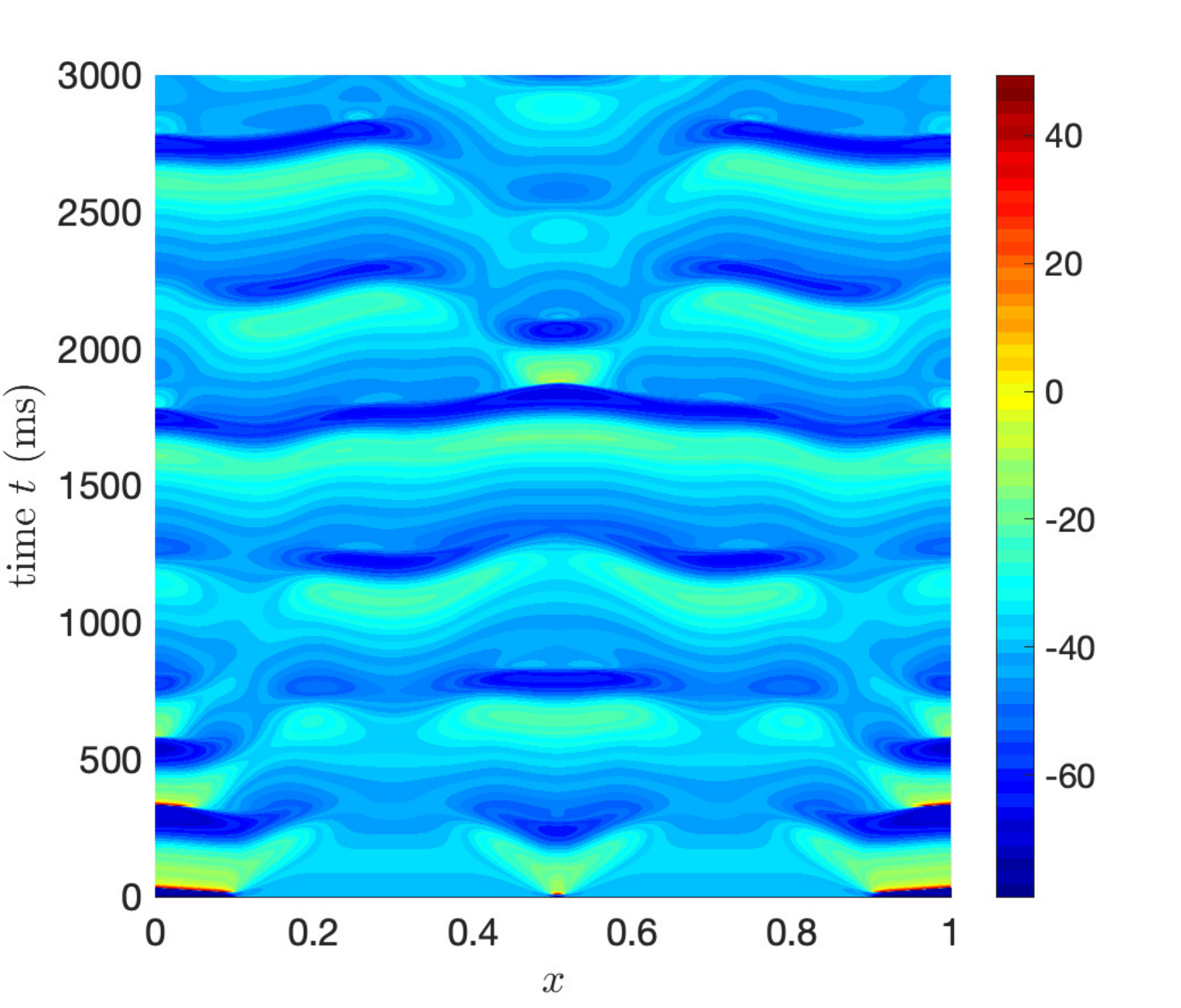}}}
\caption{Simulation of system~(\ref{monodomain}) with $G_\mathrm{L}=0.1845\,\frac{mS}{cm^2}$ and initial condition \eqref{eq:macro_init} with $\mathcal{D}=(0.1,0.5)\cup(0.51,0.9)$ for two different diffusion parameters. Computed with $256$ cells.} \label{fig:mono-diffusion}
\end{figure}

Moving the chaotic regime to the middle of the interval, the system can tolerate a much larger area of initially chaotic cells, see Figure \ref{fig:mono-diffusion}(a). Here the initial data is set to \eqref{eq:macro_init} with $\mathcal{D}=(0.1,0.5)\cup(0.51,0.9)$.  However, reducing the diffusion parameter, the behaviour of the system turns chaotic, see Figure \ref{fig:mono-diffusion}(b).

One can see from above that there is a critical mass of initially chaotic cells, depending on the diffusion parameter and the placement of the initially chaotic region, for inducing chaotic behaviour along the whole cable. These results are in concert with earlier observations noting that enough cells have to be triggered for chaotic behaviour to prevail \cite{WEISS, DELANGE2012365}. Furthermore, although induced differently, similar behaviour to what we see in Figure \ref{fig:mono-diffusion} has been observed in experiments \cite{QU_review}. However, as noted in the previous section, the voltages seem to be within a non-physiological range. 

\section{Dynamics of the modified Bernus model}\label{sec:Bernus_bif}
\noindent 
\sus{In this section we investigate the slightly modified version of the human ventricular cardiac cell model from~\cite{Bernus}, i.e. system~\eqref{model_bernus}.} As described in Section \ref{sec:bernus_descr}, this model contains more ion currents, pumps and exchangers compared to model~(\ref{model}), including the missing calcium current $I_\text{Ca}$ and the fast and slow potassium current, $I_{\mathrm{K}_\mathrm{r}}$ and $I_{\mathrm{K}_\mathrm{s}}$. Therefore, one may expect that this model is more realistic and exhibits different and diverse dynamics. System~\eqref{model_bernus} contains the important $I_\mathrm{Ca}$ current as well as the fast potassium current $I_{\mathrm{K}_\mathrm{r}}$, which are important for EADs to establish. Hence, system~\eqref{model_bernus} may exhibits EADs dependent on the choice of system parameters. 

\begin{figure}[h]
\centering
\subfigure[regular action potential]{{\includegraphics[width=0.5\textwidth]{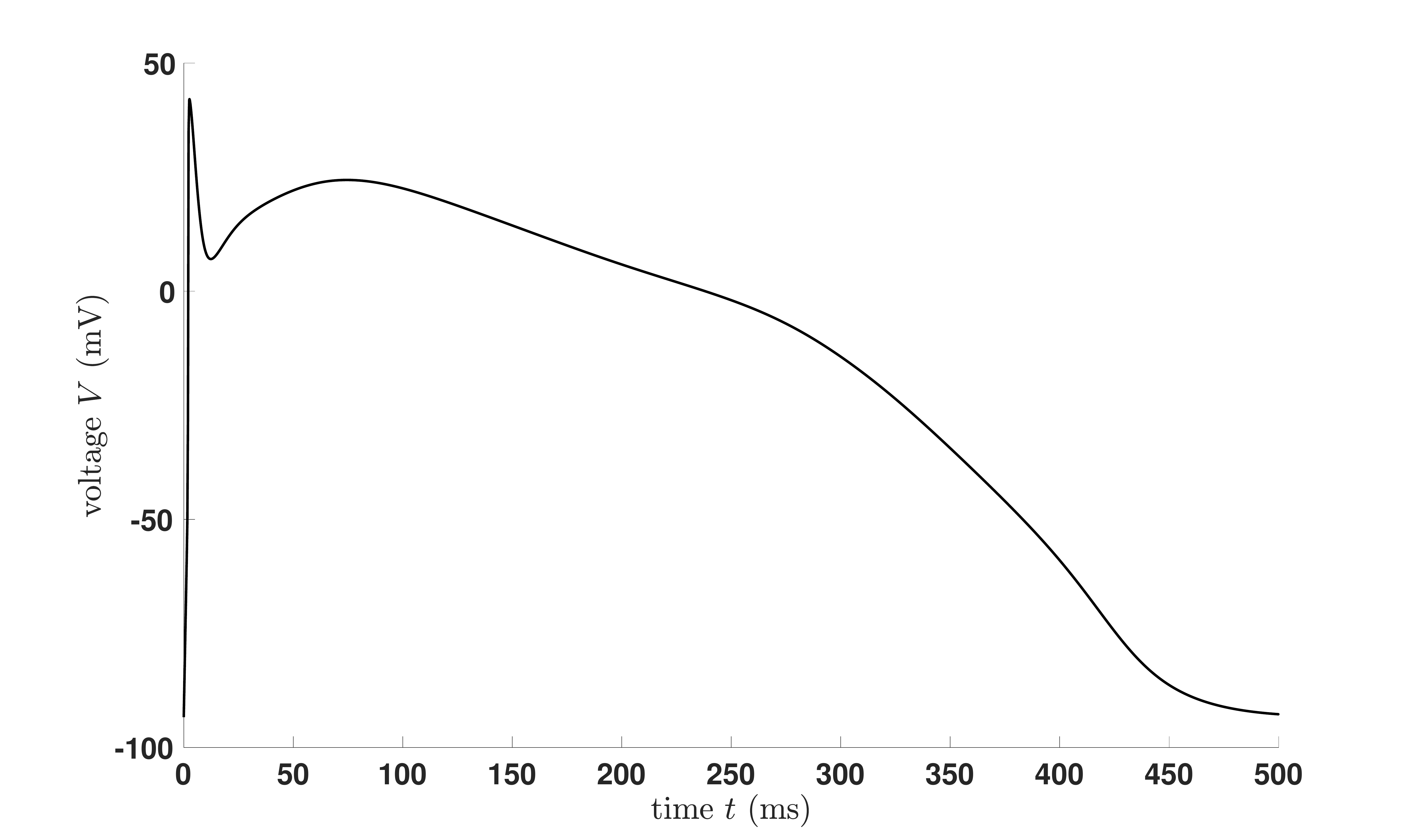}}}\subfigure[different EADs]{{\includegraphics[width=0.5\textwidth]{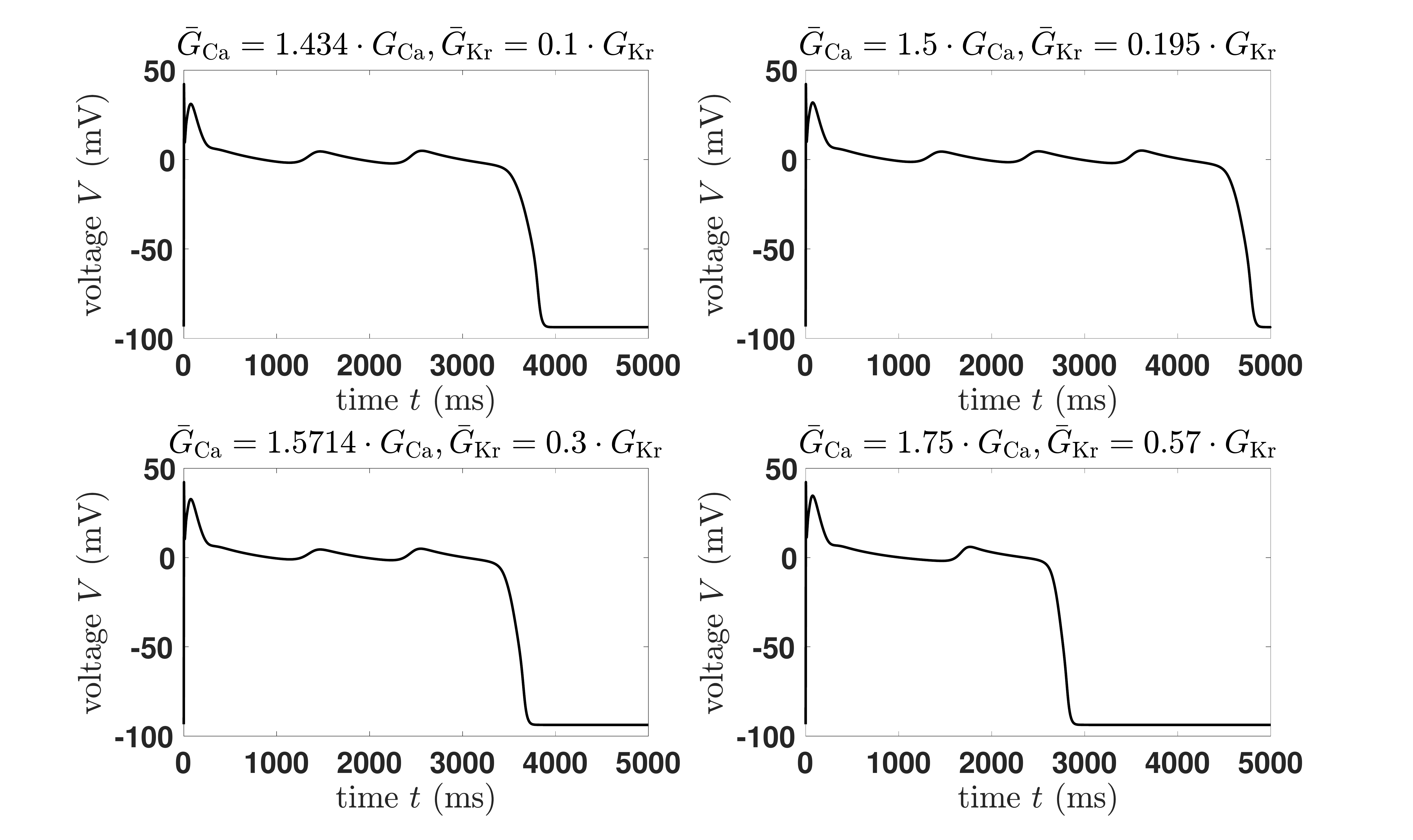}}}\caption{Comparison of trajectories of system~\eqref{model_bernus}. a) normal action potential (standard setting). b) different EADs induced by different $I_{\mathrm{K}_\mathrm{r}}$ reduction and $I_\mathrm{Ca}$ enhancements.} \label{fig:EADs}
\end{figure}

Furthermore, we know that a reduced (fast) potassium current and/or an enhanced calcium current may leads to EADs, cf.~\cite{VN}. In addition, in~\cite{Ae_control} it is shown how combinations of reduced and enhanced potassium and calcium currents increases the risk of the appearance of EADs, or conversely, may control the pattern of the APs. To illustrate the behaviour of system~\eqref{model_bernus} Figure \ref{fig:EADs} shows a comparison of a normal AP and the occurrence of EAD patterns for different combinations of a reduced fast potassium and an enhanced calcium current. 

It is well known that EADs occur either in the plateau or in the repolarisation phase of the AP and they are benefited by an elongation of the AP. This may happen by an increase of one or more inward currents and/or a decrease in one or more outward currents~\cite{WEISS}. In addition, it is well established that the calcium current $I_\text{Ca}$ plays an important rule during the plateau phase, while the potassium current $I_\text{K}$ during the repolarisation phase, cf.~\cite{Ae_control}. 

Figure \ref{fig:EADs}(b) shows that EADs appear as a combination of a reduced fast potassium current $I_{\mathrm{K}_r}$ and an enhanced calcium current $I_{\mathrm{Ca}}$. 
Therefore, we restrict our analysis to the case where we have a 80\% block of the fast potassium current by introducing a new conductance $\bar{G}_{\mathrm{K}_r}=0.2\cdot G_{\mathrm{K}_r}$ and choosing $G_\mathrm{Ca}$ as bifurcation parameter.

\begin{figure}[h]
\centering
\subfigure[Different EADs appearing in system \eqref{model_bernus} with a fixed $I_{\mathrm{K}_\mathrm{r}}$ blockade of $80\%$ and different $I_\mathrm{Ca}$ enhancements.]{{\includegraphics[width=0.5\textwidth]{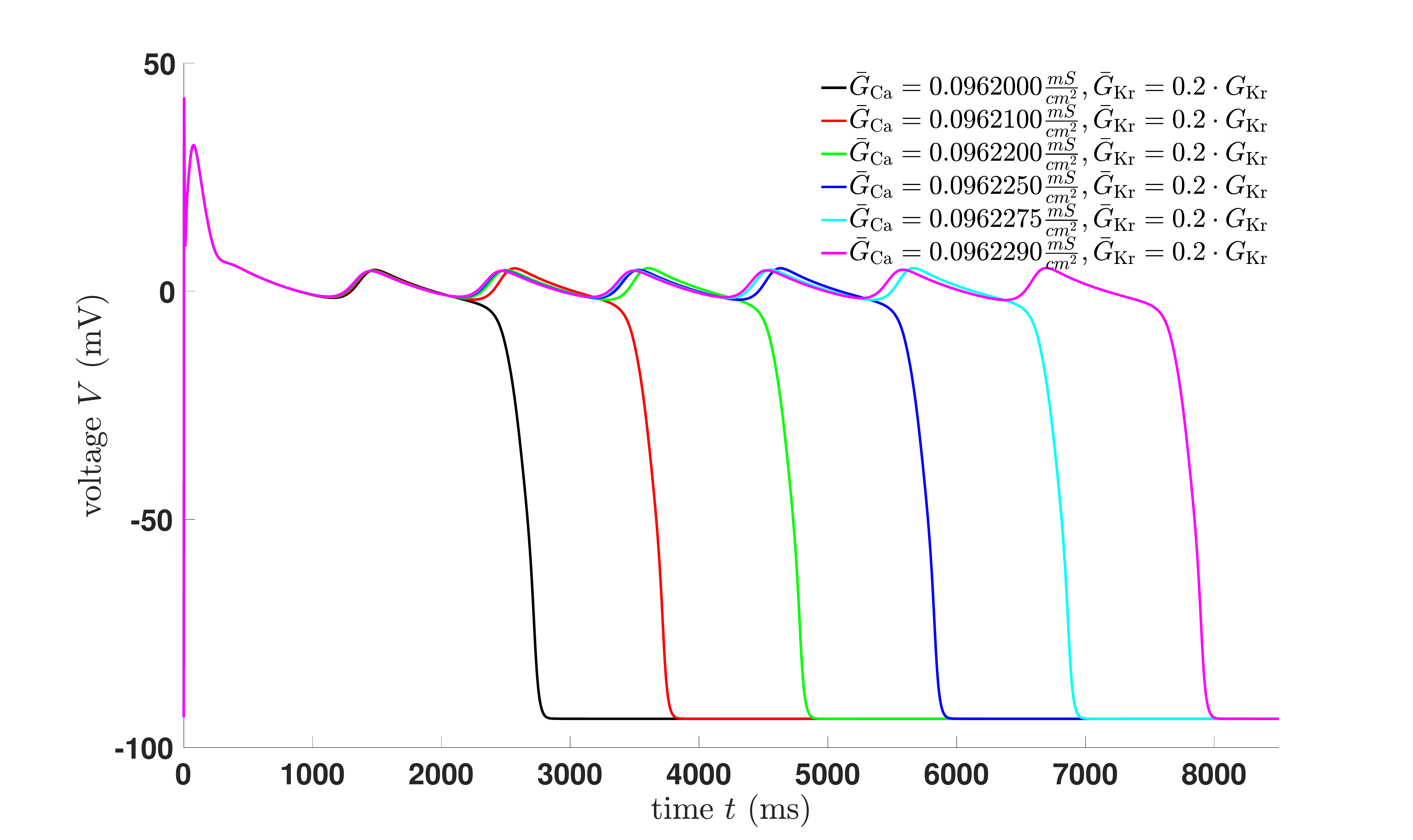}}}\subfigure[2D bifurcation diagram.]{\includegraphics[width=0.5\textwidth]{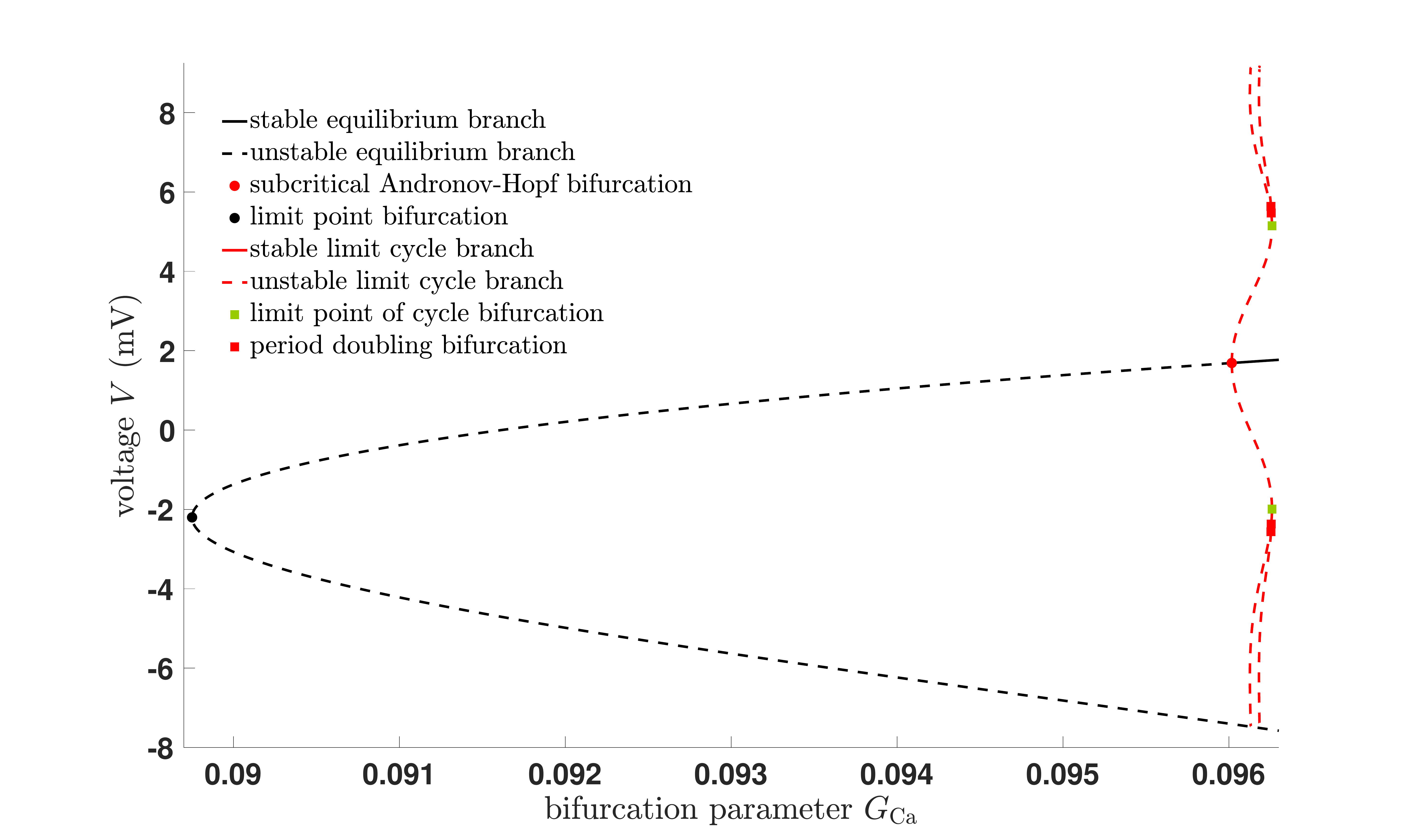}}
\subfigure[Zoom of b) around the limit cycles.]{\includegraphics[width=0.5\textwidth]{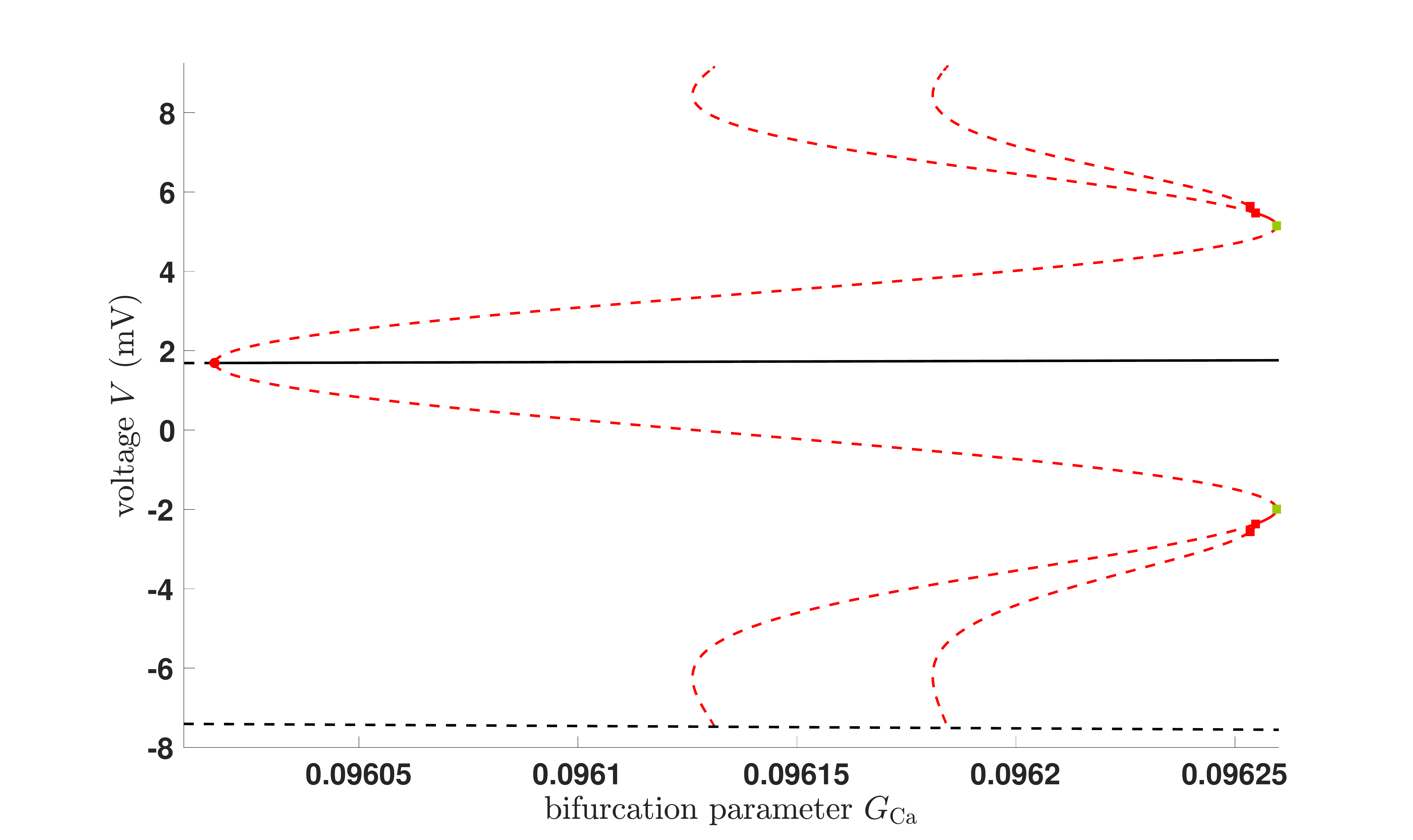}}\subfigure[3D bifurcation diagram.]{\includegraphics[width=0.5\textwidth]{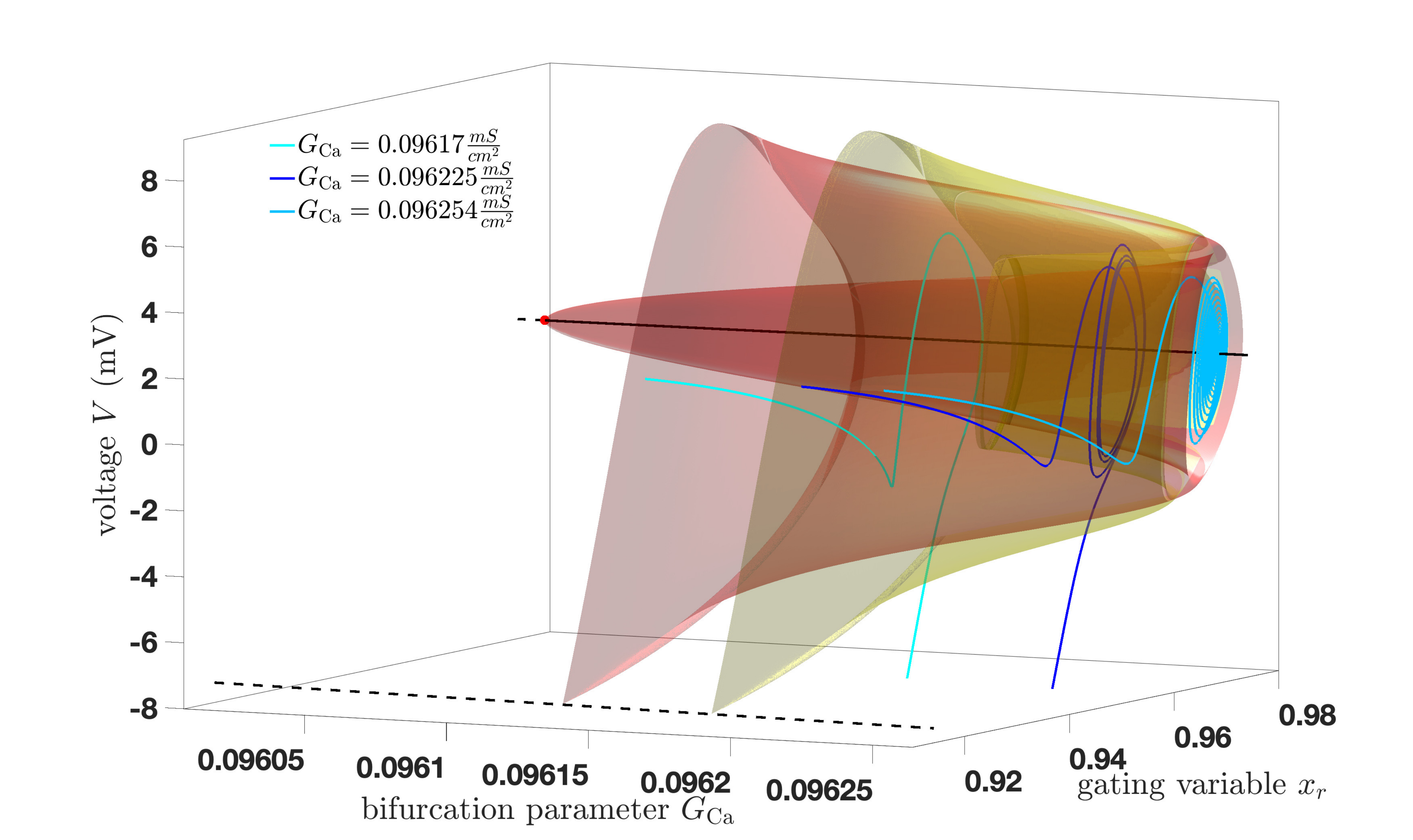}}
\caption{Bifurcation diagram of system~\eqref{model_bernus} including a 80\% block of $I_{\mathrm{K}_\mathrm{r}}$, i.e. $\bar{G}_{\mathrm{K}_\mathrm{r}}=0.2\cdot{G}_{\mathrm{K}_\mathrm{r}}$, using $G_\mathrm{Ca}$ as bifurcation parameter.} \label{fig:bif_Bernus}
\end{figure}

Similar to the analysis of the Noble model~\eqref{model} we start by determining the equilibrium curve, cf. Figure \ref{fig:bif_Bernus}(b) and (c). Here, we have again an unstable (black dashed line) and a stable (black solid line) equilibrium branch. The equilibrium curve changes stability via a subcritical Andronov--Hopf bifurcation (red dot, $G_\mathrm{Ca}\approx 0.096017\frac{mS}{cm^2}$) with a positive first Lyapunov coefficient. From the subcritical Andronov--Hopf bifurcation an unstable limit cycle branch (dashed red line) bifurcates, cf. Figure \ref{fig:bif_Bernus}(b) and (c). 

{The first limit cycle branch contains a limit point of cycle bifurcation (solid green square, $G_\mathrm{Ca}\approx 0.096259\frac{mS}{cm^2}$) and a period doubling bifurcation (solid red square, $G_\mathrm{Ca}\approx 0.096255\frac{mS}{cm^2}$) from which a stable period doubling cascade bifurcates. The first limit cycle branch changes stability via the limit point of cycle bifurcation and again via the (first) period doubling bifurcation.} The bifurcation diagram in Figure \ref{fig:bif_Bernus} contains the first two limit cycle branches, where both limit cycles terminate at the unstable equilibrium branch. The limit cycle branches are mostly unstable and therefore, not attracting. Nevertheless, they influence the dynamics of the system, i.e they at least prolong the plateau phase and may cause EADs, provided the initial stimulus is strong enough to establish an AP (as it is in the standard setting). 

In addition to Figure \ref{fig:bif_Bernus}(b) and (c), we provide in Figure \ref{fig:bif_Bernus}(d) the corresponding 3D bifurcation diagram including three different trajectories. From this it is obvious that the trajectories curl around the limit cycles resulting in a prolongation of the AP and/or EADs. Indeed, as soon as the trajectories enters the inside of the limit cycle branches, they will converge into the stable equilibrium branch. Notice that the stable equilibriums close to the Andronov--Hopf bifurcation are less attracting than others, since some of the negative eigenvalues are very small, but still negative.

\begin{figure}[h]
\centering
\includegraphics[width=0.77\textwidth]{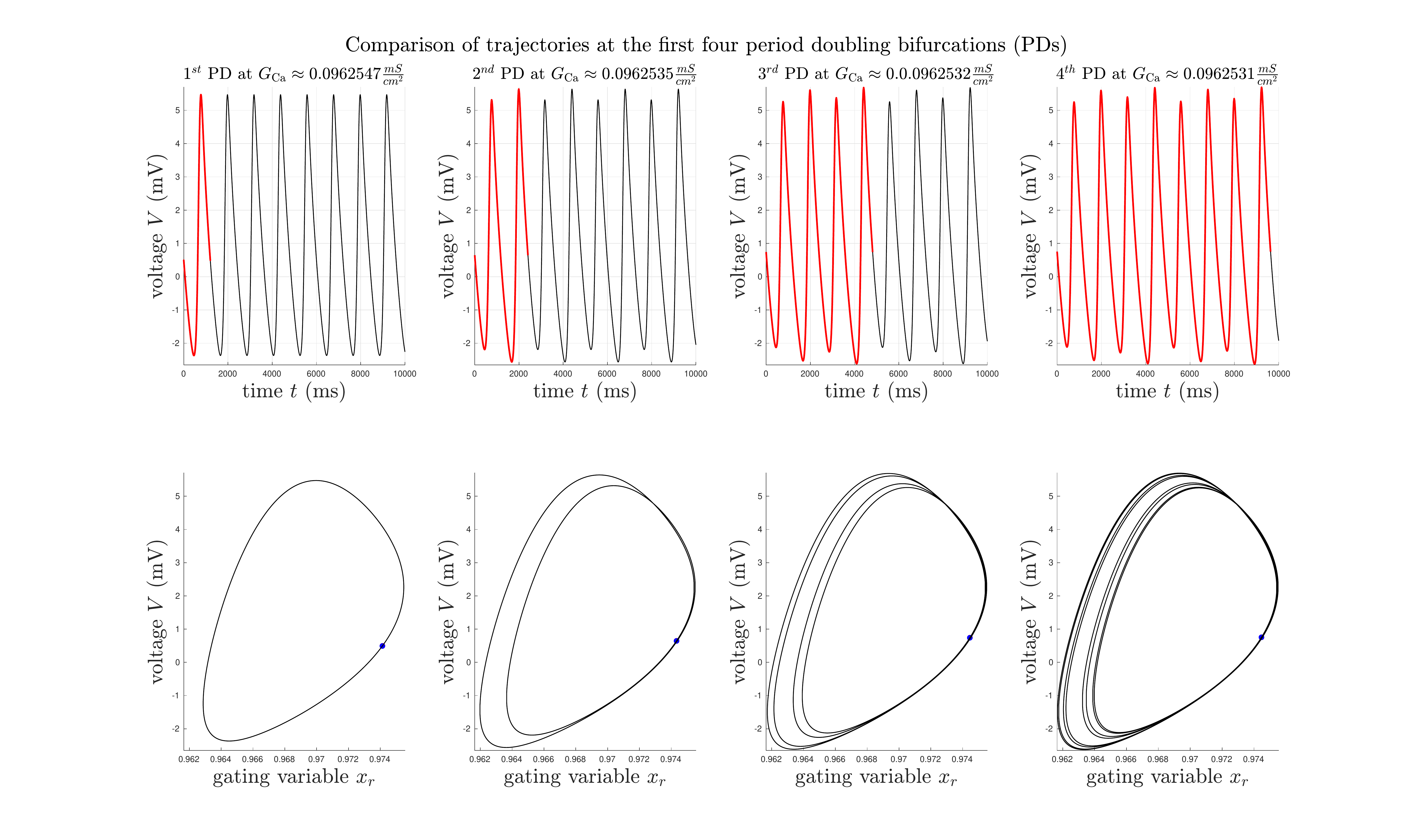}
\caption{Simulation of system~\eqref{model_bernus} over $10000\ ms$ for different $G_\mathrm{Ca}$ values with $I_\mathrm{stimulus}=0$ and initial values on the corresponding limit cycle branch.}\label{fig:pd_sim}
\end{figure}

\begin{figure}
\centering
\subfigure[Chaos scenario]{\includegraphics[width=0.35\textwidth]{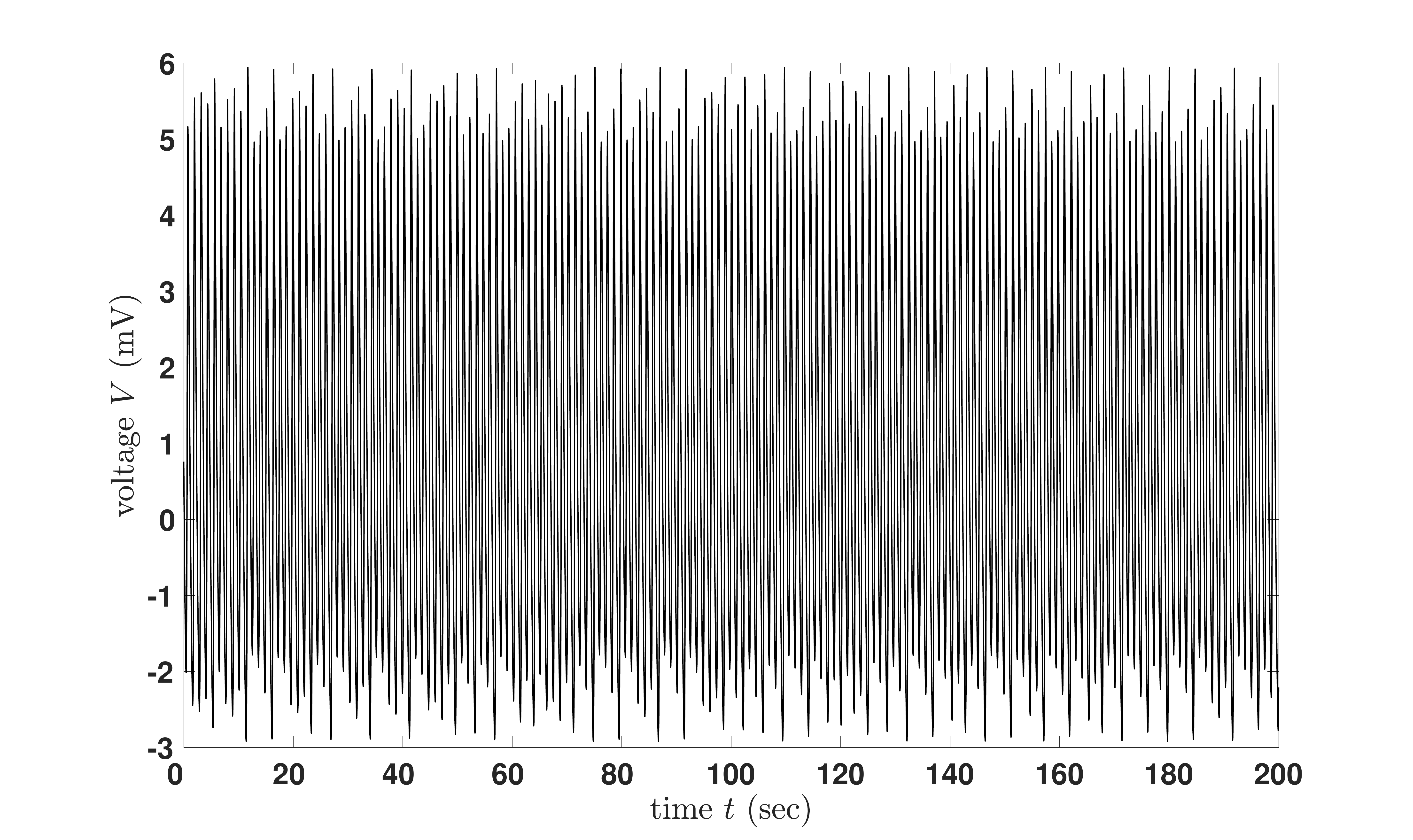}}\subfigure[Non-chaos scenario]{\includegraphics[width=0.35\textwidth]{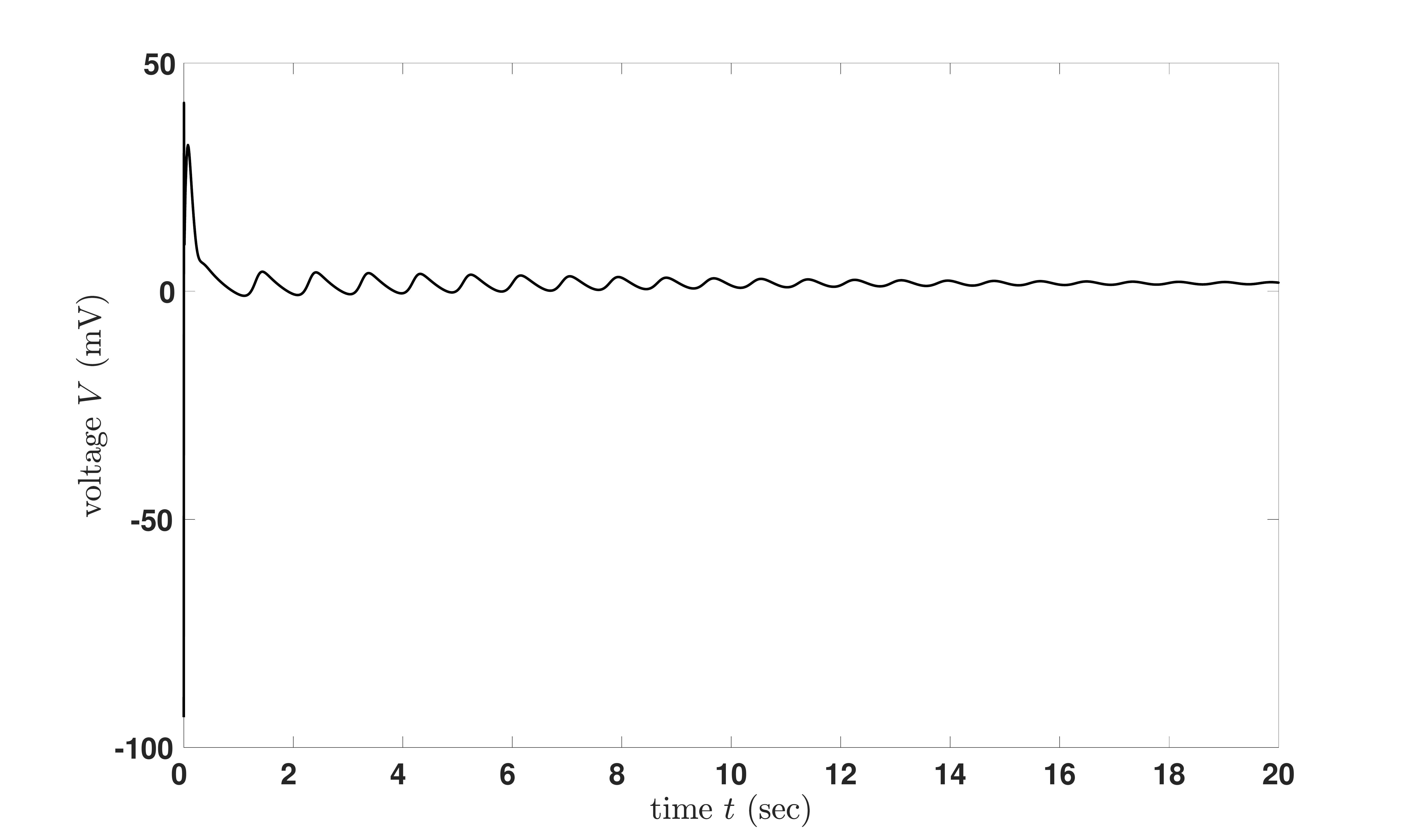}}
\caption{Simulation of system~\eqref{model_bernus} for $G_\mathrm{Ca}=0.0962518\frac{mS}{cm^2}$. (a) Chaos scenario: initial values 
$V_0=0.7589\ mV$, $m_0=0.9952$, $v_0=0.0$, $f_0=0.134$, $x_{r_0}=0.9745$, $to_0=0.0028$, $d_0=0.9141$, $r_0=0.0411$, $K1_0=0.0$, $x_{s_0}=0.5485$, and $I_\mathrm{stimulus}=0$. (b) Non-chaos scenario: standard initial values and external stimulus.}\label{fig:bernus_chaos}
\end{figure}

However, if the trajectory of system~\eqref{model_bernus} is in the basin of attraction of the stable period doubling cascade (the stable attracting parts of the limit cycle branches), the system develops self-oscillating behaviour (Figure~\ref{fig:pd_sim}), or chaos (Figure~\ref{fig:bernus_chaos}). A setting for this to happen is the combination of a $G_\mathrm{Ca}$ value of the period doubling cascade with initial values in the basin of attraction and $I_\mathrm{stimulus}=0$.

{Figure~\ref{fig:pd_sim} contains four simulations of system~\eqref{model_bernus} at the first four period doubling bifurcations of the period doubling cascade. The black fine line denotes the trajectory over $10000\ ms$, while the red line indicates the length of the trajectory with period $T$, i.e. $V(T)-V_0=0$. Additionally, Figure~\ref{fig:pd_sim} provides the corresponding phase space $(x_r,V)$ to these simulations showing a closed curve starting from a red dot and terminating at a blue one. Note that the red dot is overlaid by the blue one due to $V(T)=V_0$ and therefore, barely or not visible.}

{Finally, depending on the initial values the period doubling cascade is again a route to chaos similar to situation for the Noble model~\eqref{model}, see Figure~\ref{fig:bernus_chaos}. Again, it is clear that besides the system parameters also the initial values, and additionally the external stimulus, play a crucial role for the occurrence of certain dynamics and patterns such as (normal) AP, EADs or chaos. This indicates that a disorder in the external stimulus may also initiate a sudden death (at least on the cellular level).}

Our analysis shows on the one hand that the dynamics of a single cell model are sensitive to its system parameters, and on the other hand they are sensitive to the choice of initial values. This is in accord with our previous analysis. However, the behaviours of the Noble model~\eqref{model} and system~\eqref{model_bernus} are quite different. Note that the strength of the initial stimulus influences the initial values and influences the dynamics of system~\eqref{model_bernus}.

Notably, the occurring EADs in Figure~\ref{fig:bif_Bernus}(a) appear in a physiological feasible range, cf. \cite{Vandersickel2}, acting as a validation of the model. However, the chaotic behaviour in Figure~\ref{fig:bernus_chaos}(a) is most likely non-physiological due to the small voltage range, and we would expect cardiac death as in Figure \ref{fig:bernus_chaos}(b) to happen in the real cell.

\subsection{Effects on the macro-scale}
Based on the analysis in the previous section, we study the synchronisation behaviour of an ensemble of cells along a 1D cable of 1~\emph{cm} to gain an intuition on whether EADs can spread or not in cardiac tissue. As for the monodomain model \eqref{monodomain}, we split the cable into two parts $\mathcal{D}$ and $[0,1]\setminus \mathcal{D}$. In the domain $\mathcal{D}$, we set the calcium conductance to $G_{\mathrm{Ca}}=0.096229\frac{mS}{cm^2}$ to ensure that an EAD with six additional oscillations (the pink line in Figure \ref{fig:bif_Bernus}(a)) would occur in the ODE model \eqref{model_bernus}, while in the remaining part it is set to $G_{\mathrm{Ca}}=0.09616\frac{mS}{cm^2}$ (close to the occurence of EADs). We keep the diffusion small and set the diffusion constant to $0.00005$\,{\it{mS}}.

\begin{figure}[h]
\centering
\subfigure[$\mathcal{D}=[0.49, 0.5)$.]{{\includegraphics[width=0.24\textwidth]{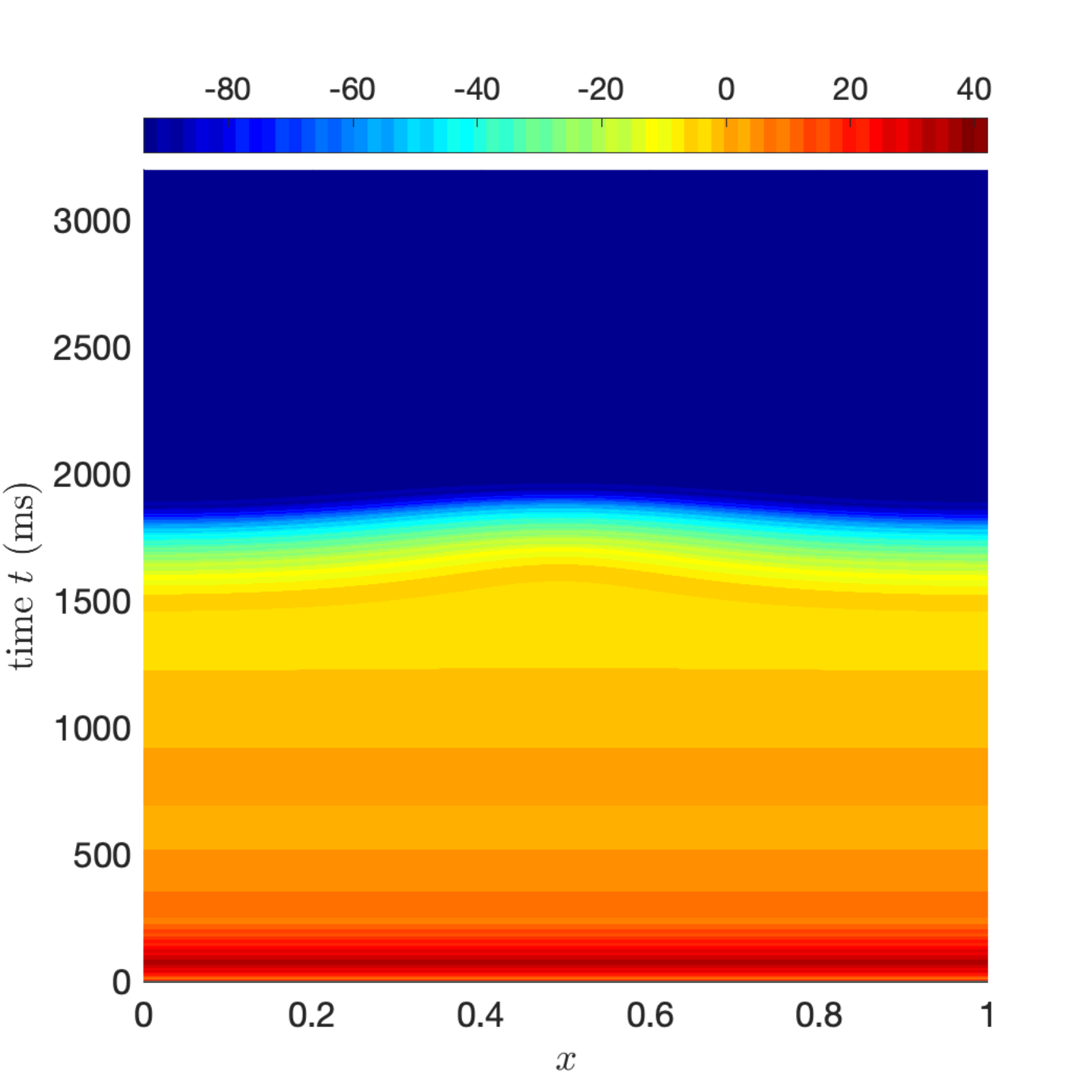}}}
\subfigure[$\mathcal{D}=[0.48, 0.5)$.]{{\includegraphics[width=0.24\textwidth]{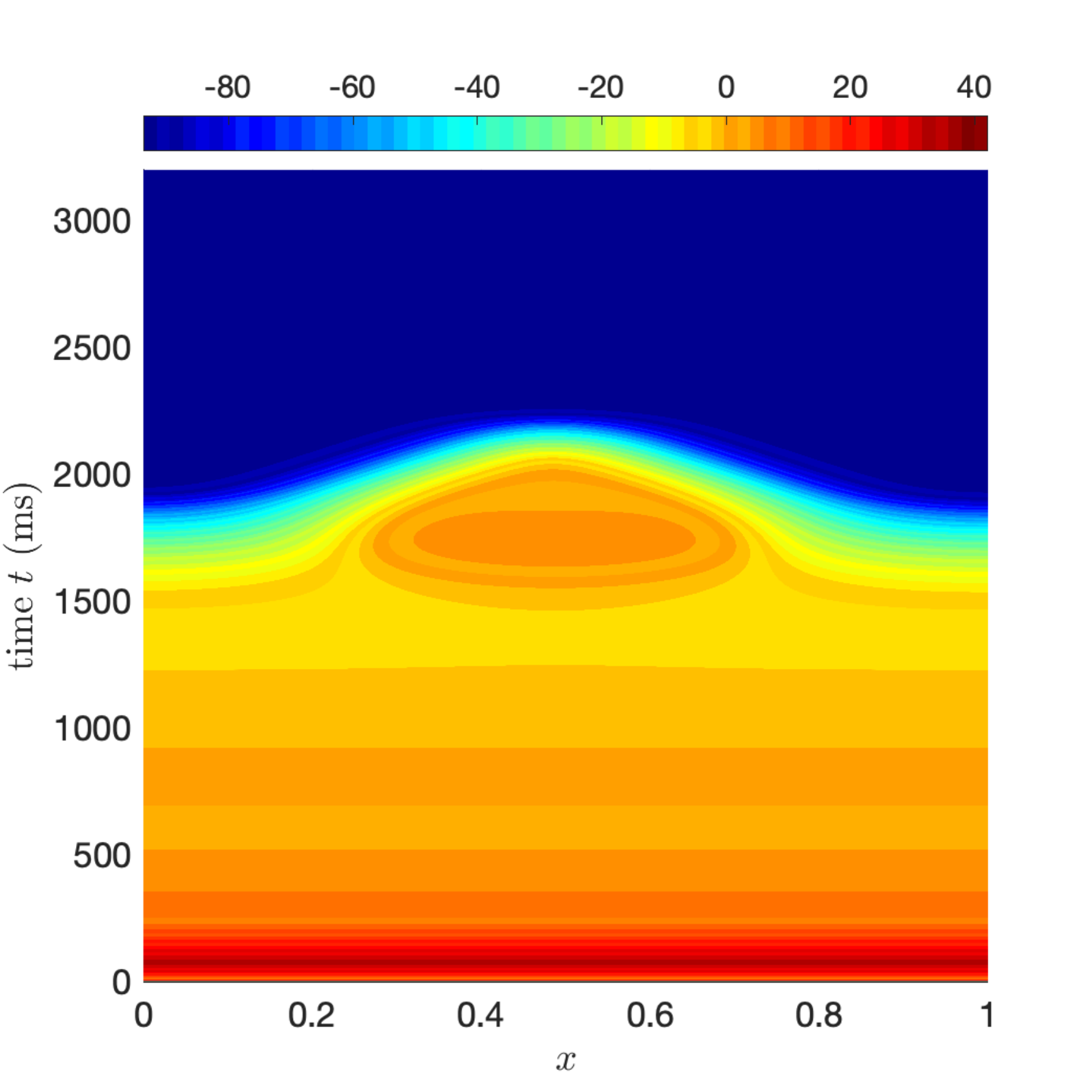}}}
\subfigure[$\mathcal{D}=[0.2,0.7)$.]{{\includegraphics[width=0.24\textwidth]{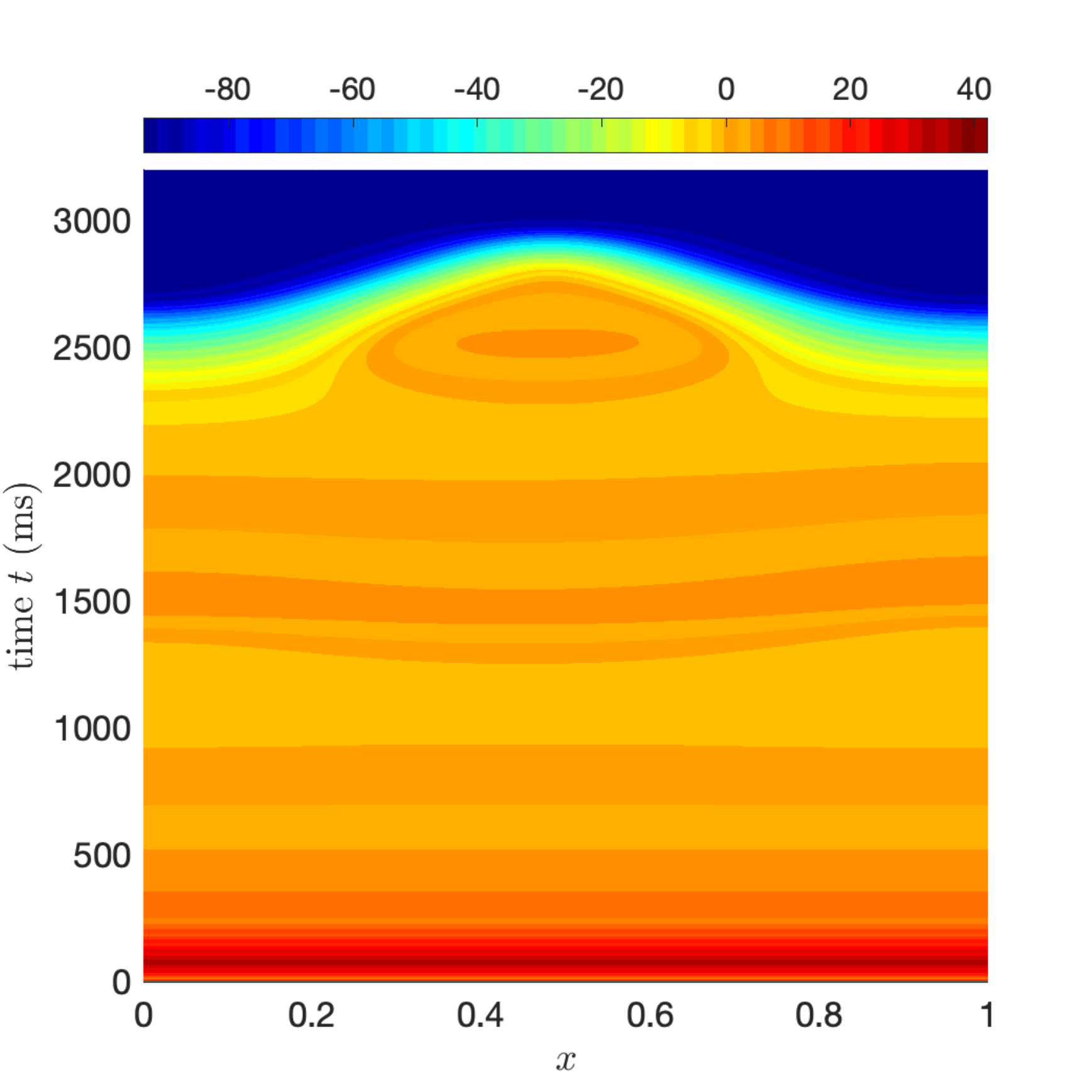}}}
\subfigure[$\mathcal{D}=[0.2,0.7)$.]{{\includegraphics[width=0.24\textwidth]{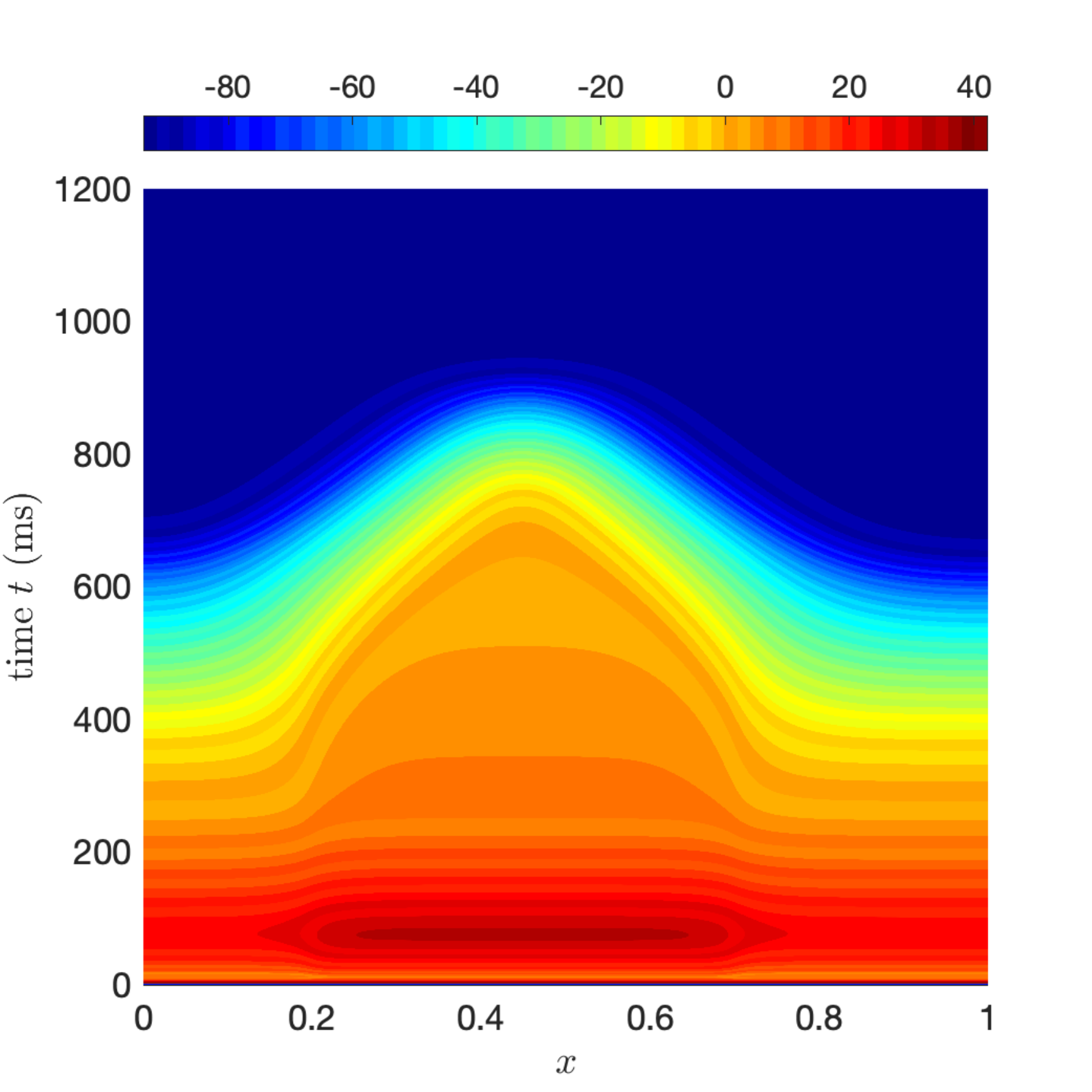}}} \\
\subfigure
{{\includegraphics[width=0.24\textwidth]{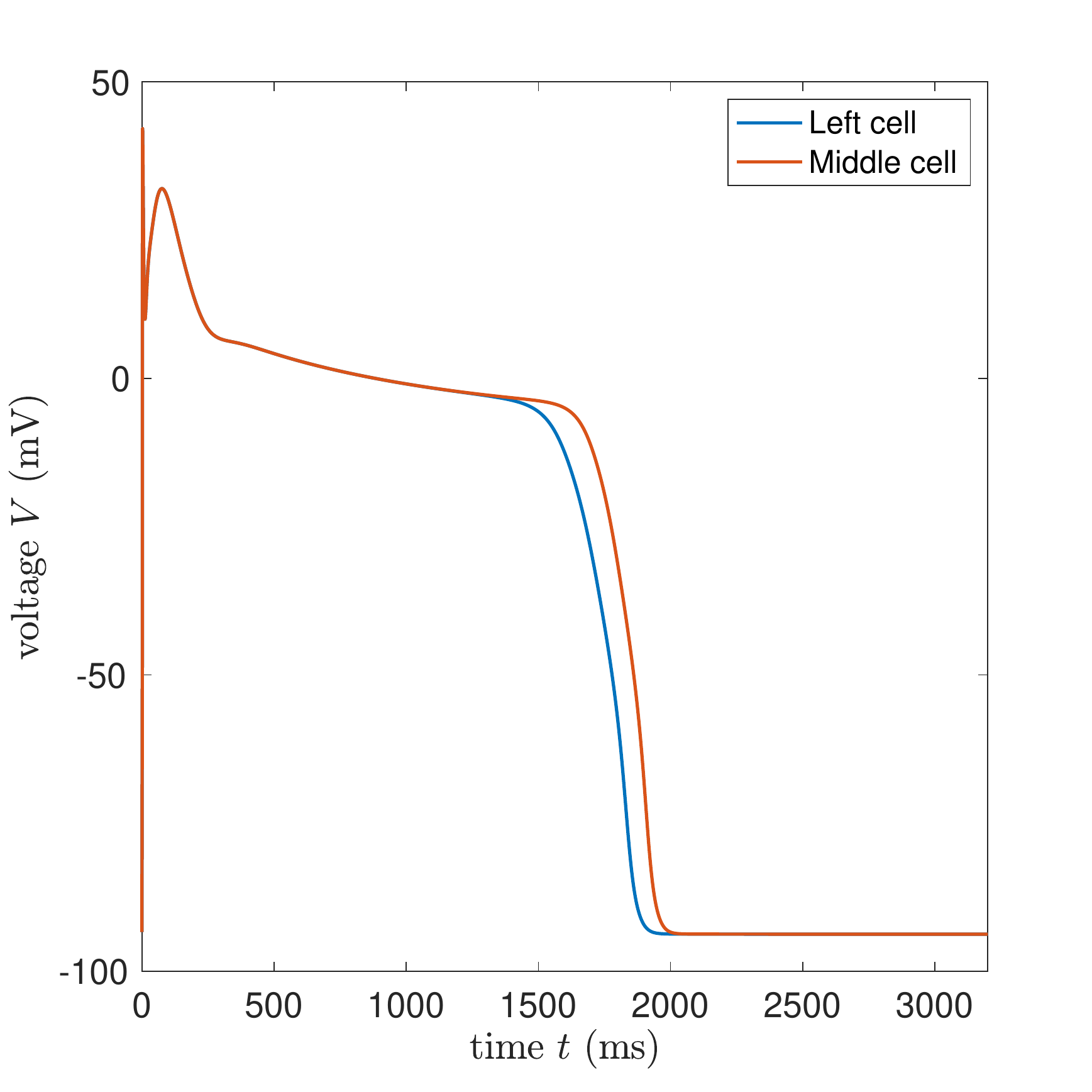}}}
\subfigure{{\includegraphics[width=0.24\textwidth]{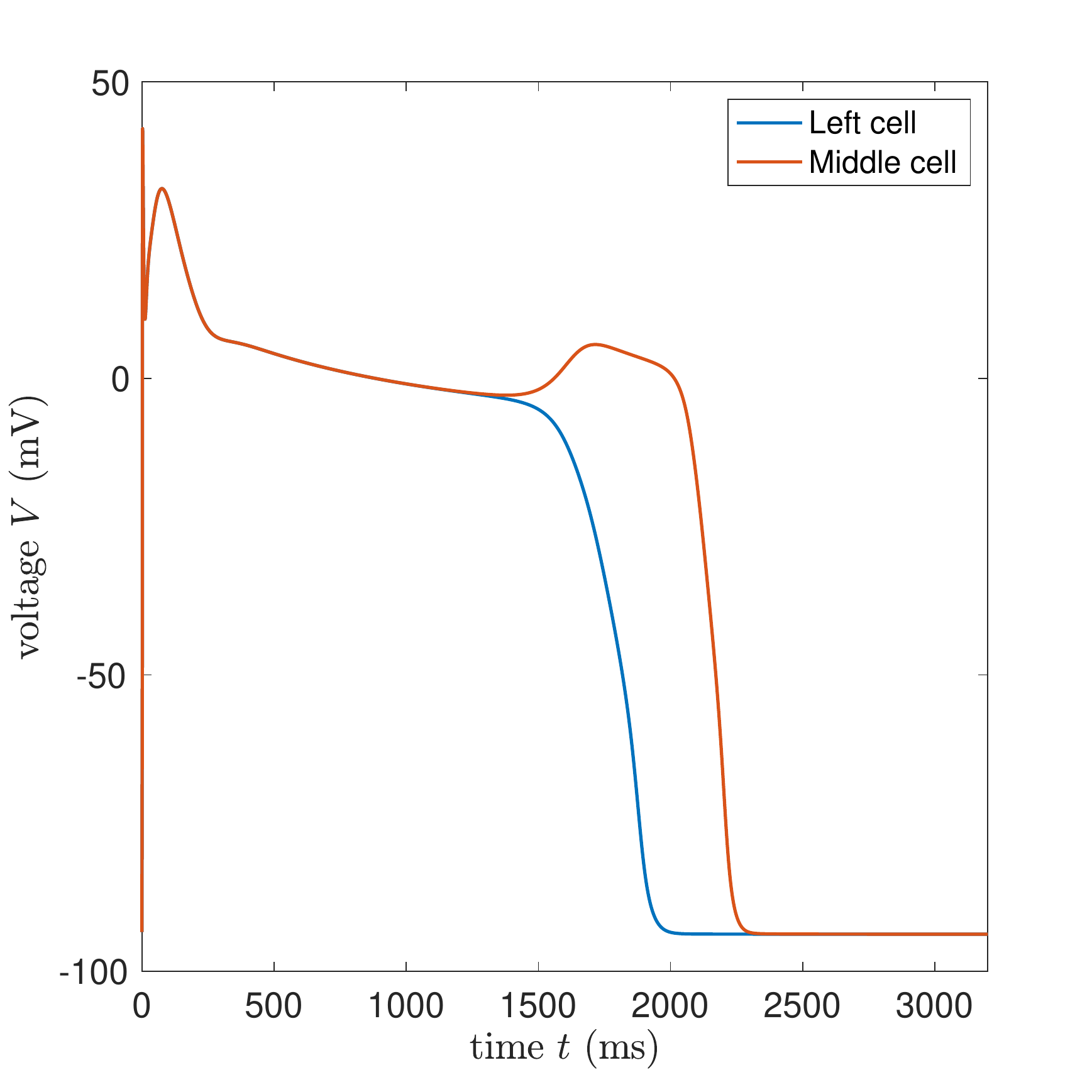}}}
\subfigure{{\includegraphics[width=0.24\textwidth]{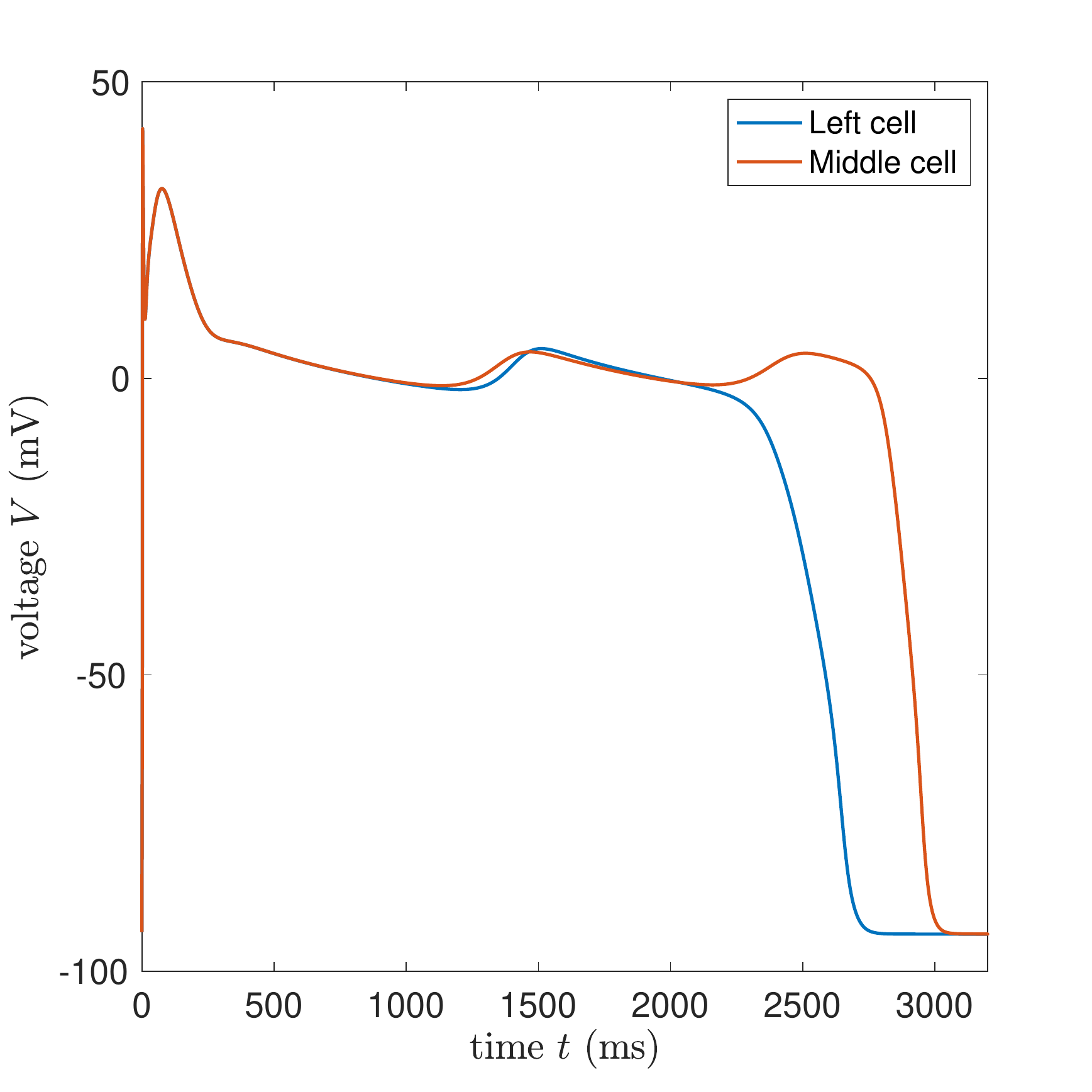}}}
\subfigure{{\includegraphics[width=0.24\textwidth]{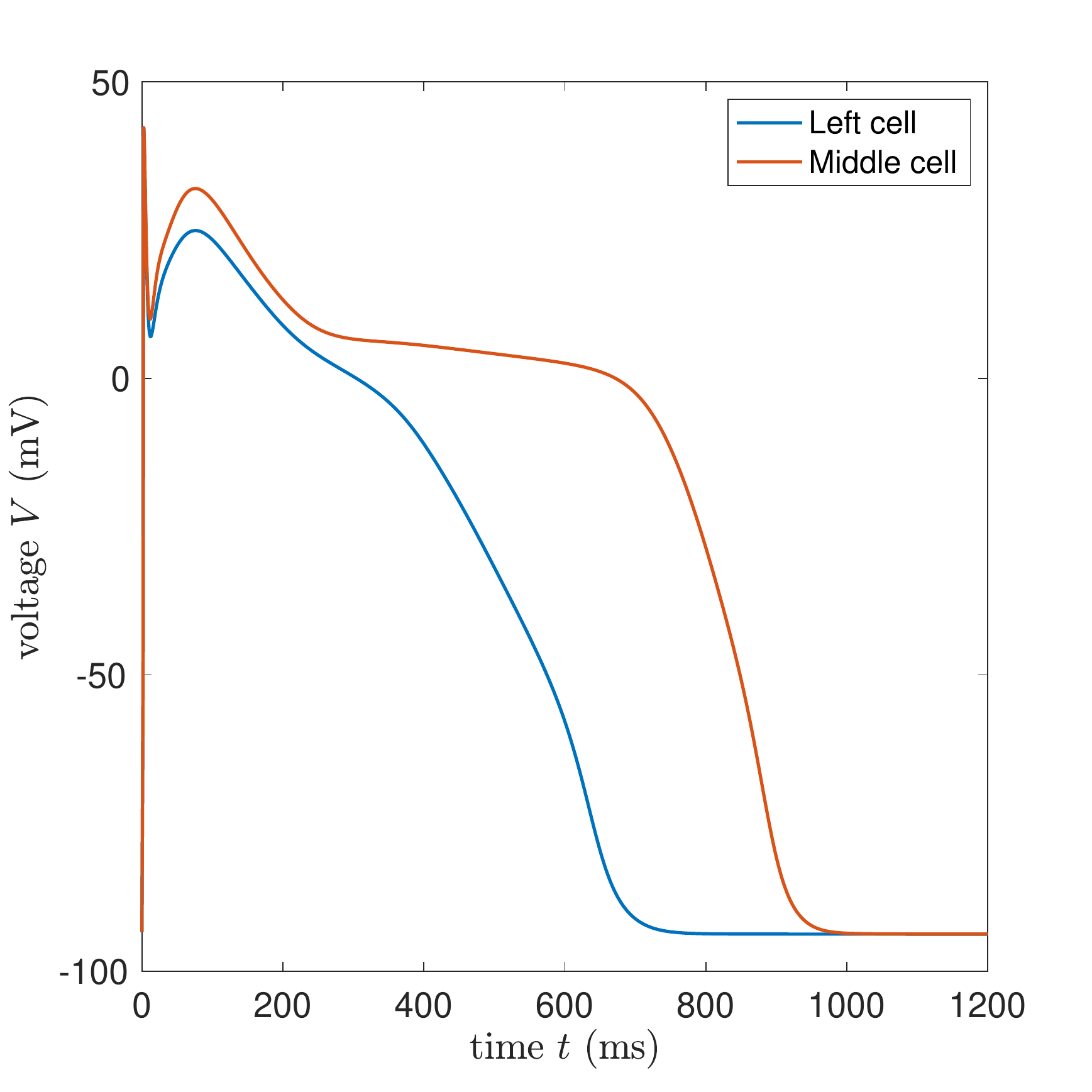}}}
\caption{Simulations of system~\eqref{monodomain_bernus} for increasing EAD domains (left to right). In $\mathcal{D}$ the calcium conductance is $G_{\mathrm{Ca}}=0.096229\frac{mS}{cm^2}$, while in the remaining parts it is set to $G_{\mathrm{Ca}}=0.09616\frac{mS}{cm^2}$ in (a)--(c), and to $G_{\mathrm{Ca}}=0.064\frac{mS}{cm^2}$ in (d). Top line: the 1D cable simulations. Bottom line: the corresponding left and middle cell (blue and red). Computed with $128$ cells.} \label{fig:mono-bernus}
\end{figure}

{
Figure~\ref{fig:mono-bernus} shows the dynamics of an ensemble of 128 cells where $1\%$, $2\%$, and $50\%$ of the cells are set to the six-oscillation setting for the ODE. We observe that no EADs occur when $1\%$ of the cells are EAD prone (Figure~\ref{fig:mono-bernus}(a)), while there is one additional small oscillation on parts of the cable for $2\%$ EAD prone cells (Figure~\ref{fig:mono-bernus}(b)). Increasing the percentage of EAD prone cells to $50\%$, we see that there are small additional oscillations along the whole cable (Figure~\ref{fig:mono-bernus}(c)). All three experiments show fewer small additional oscillations than in the ODE case.}

{
In Figure~\ref{fig:mono-bernus}(a)--(c) the cells surrounding the EAD prone cells ($x \in [0,1]\setminus \mathcal{D}$) are very close to establishing EAD behaviour. Hence, only a very small percentage of EAD prone cells are needed for EADs to occur along the cable. If the surrounding cells are further away from the EAD setting ($G_{\mathrm{Ca}}=0.064\frac{mS}{cm^2}$, i.e. the standard setting), we can observe that no EADs occur even if $50\%$ of the cells are set to EAD inducing behaviour, see Figure~\ref{fig:mono-bernus}(d).
}

\begin{figure}[h]
\centering
\subfigure[Diffusion const. $\frac{1}{360}$\,{\it{mS}}.]{{\includegraphics[width=0.28\textwidth]{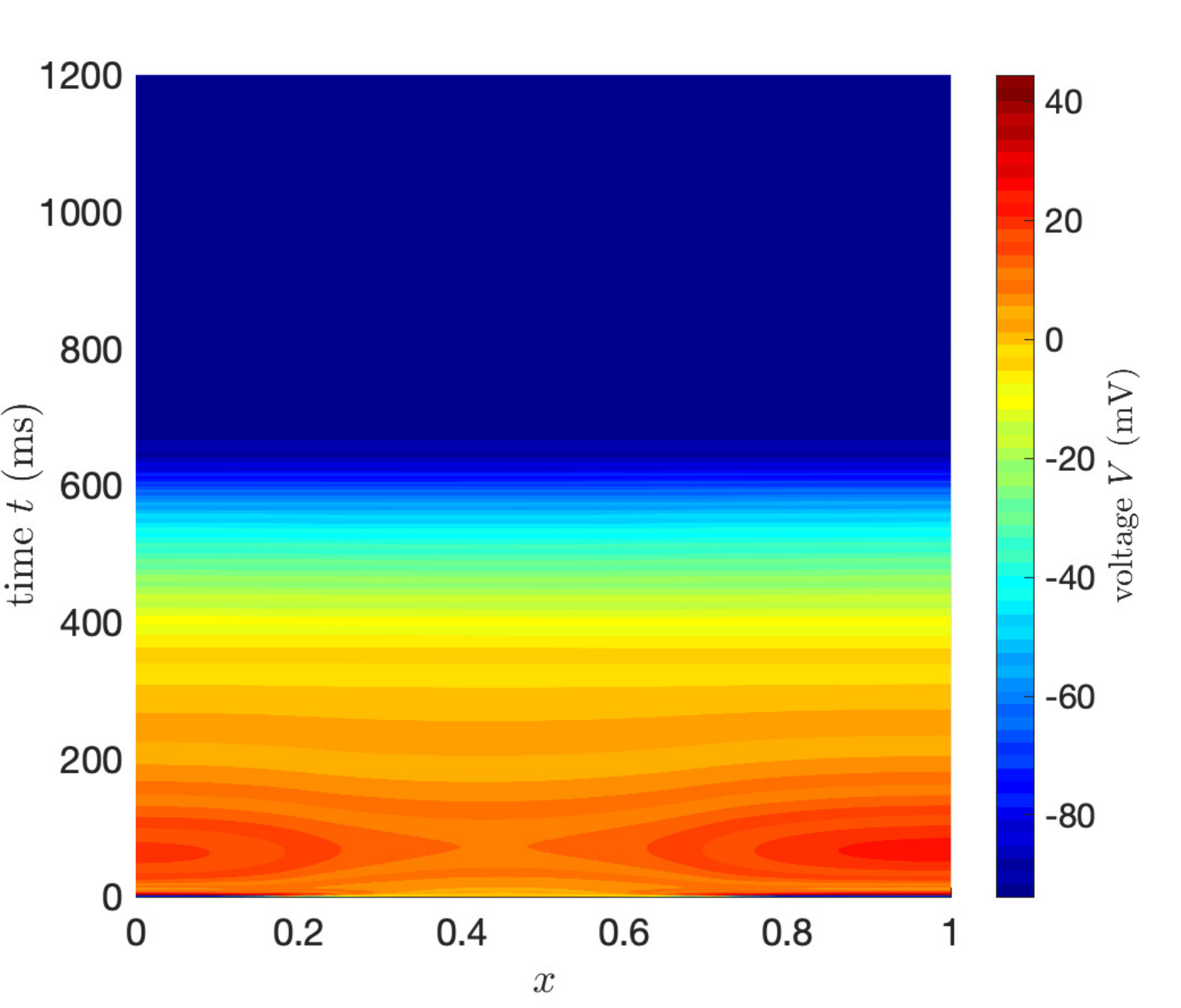}}}
\subfigure[Diffusion const. $0.00005$\,{\it{mS}}.]{{\includegraphics[width=0.28\textwidth]{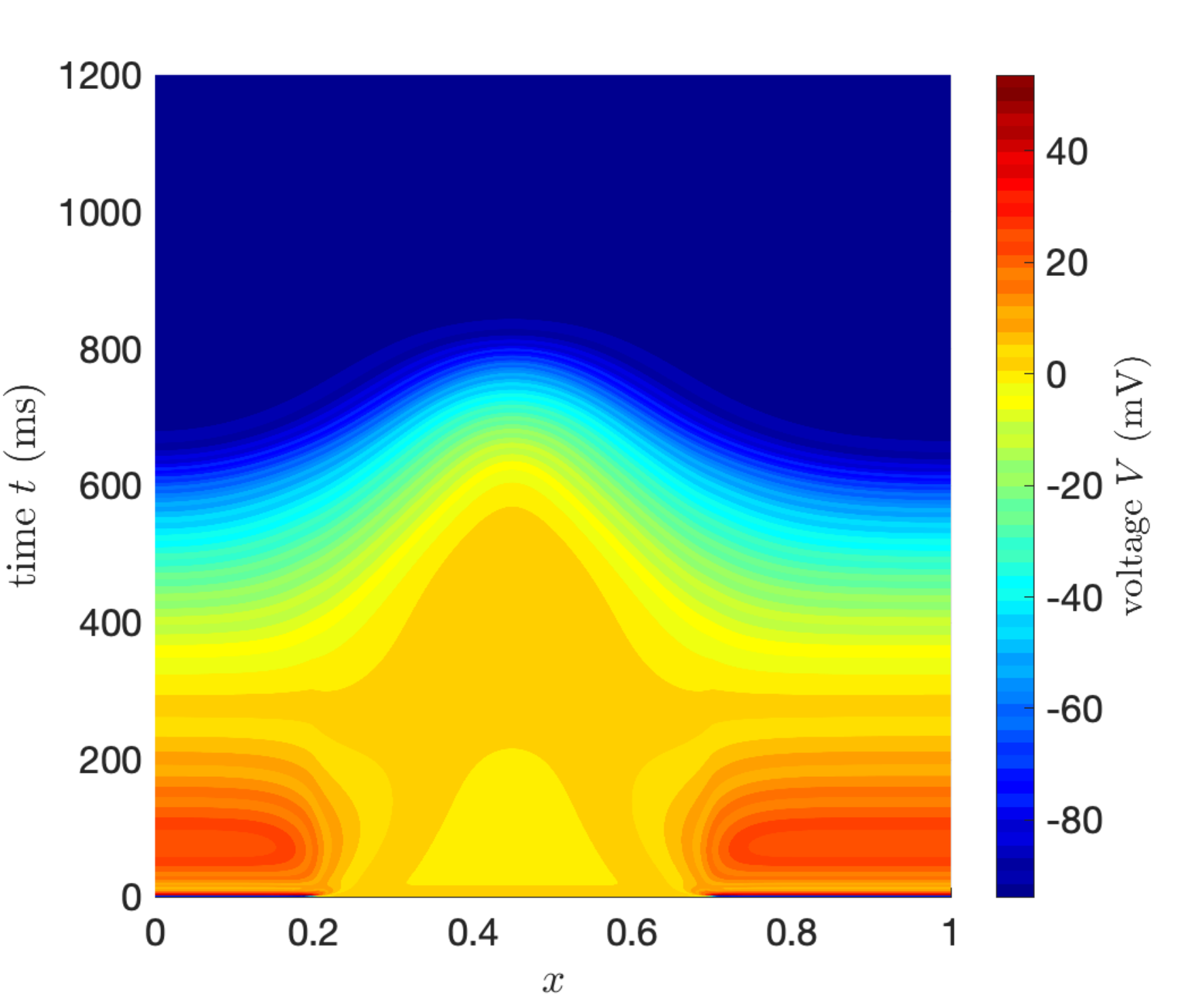}}}
\subfigure[Diffusion const. $0.00005$\,{\it{mS}}.]{{\includegraphics[width=0.28\textwidth]{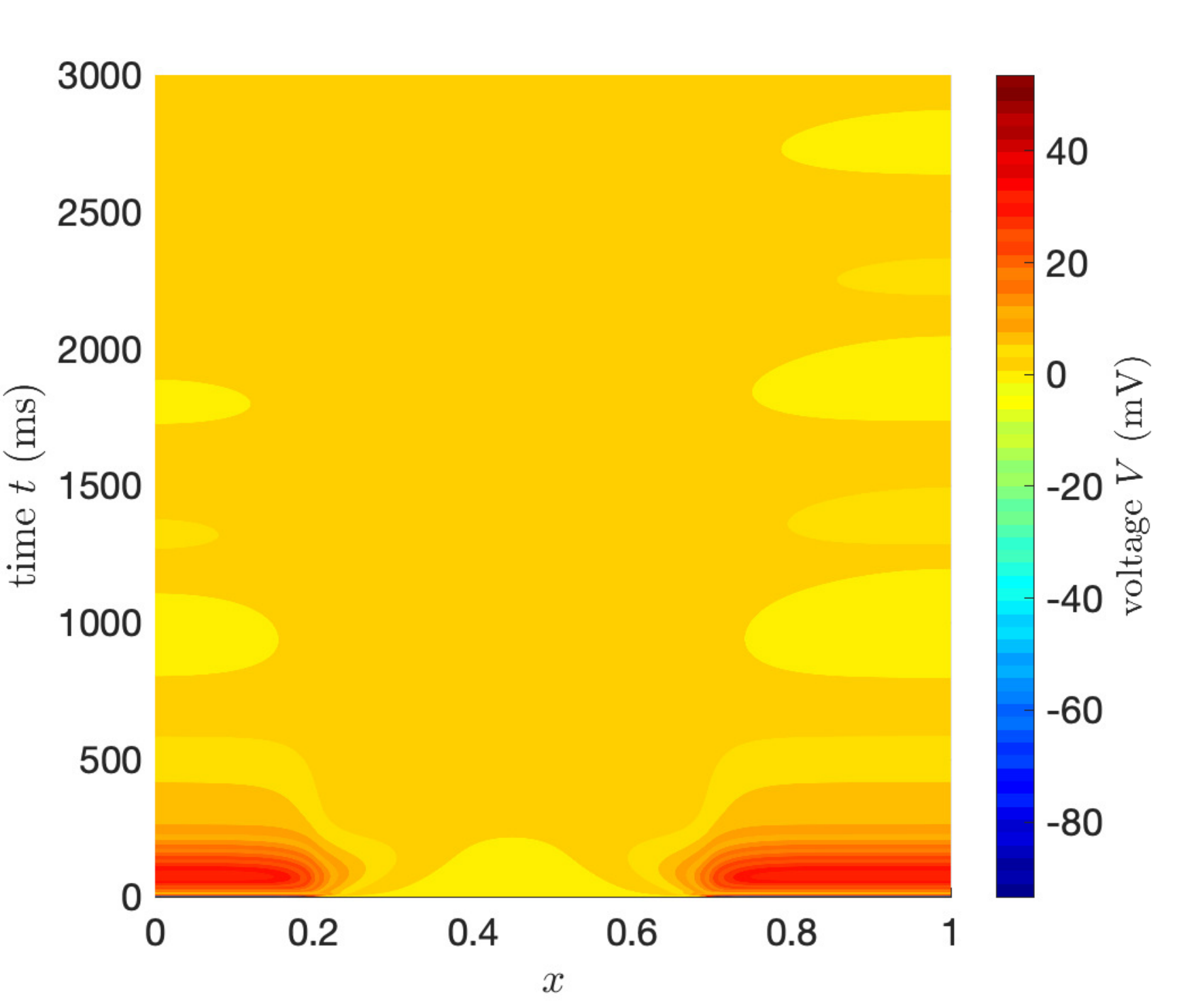}}}\\
\subfigure
{{\includegraphics[width=0.28\textwidth]{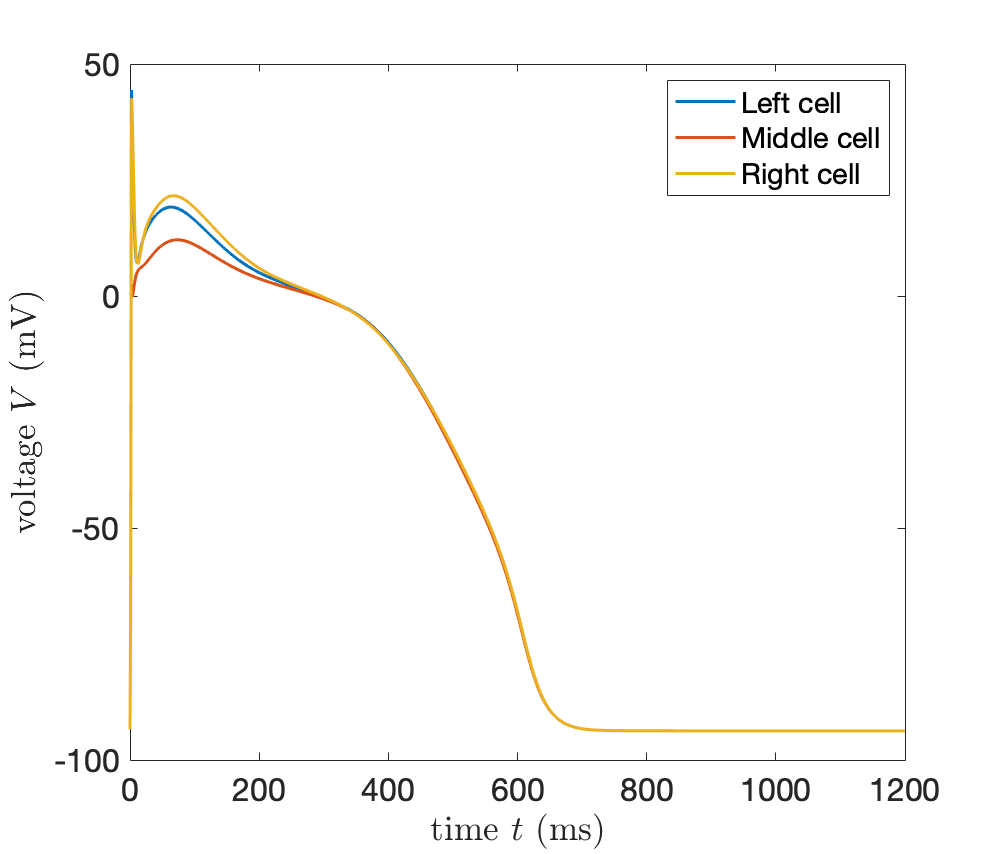}}}
\subfigure{{\includegraphics[width=0.28\textwidth]{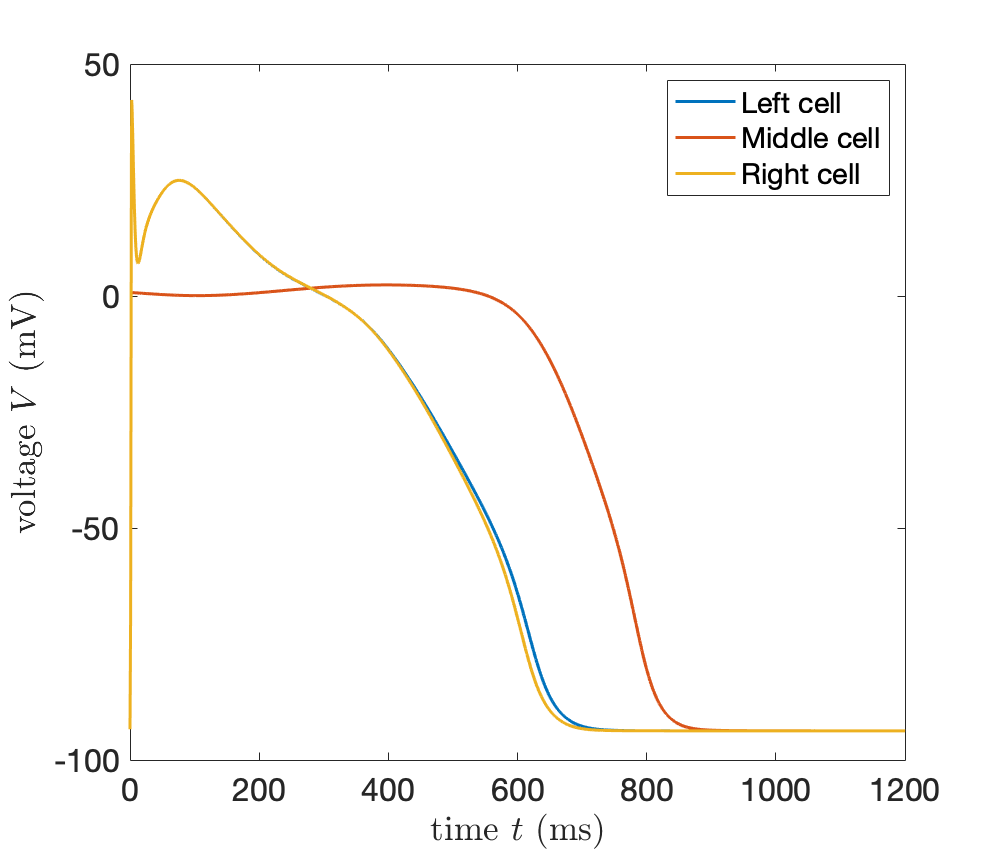}}}
\subfigure{{\includegraphics[width=0.28\textwidth]{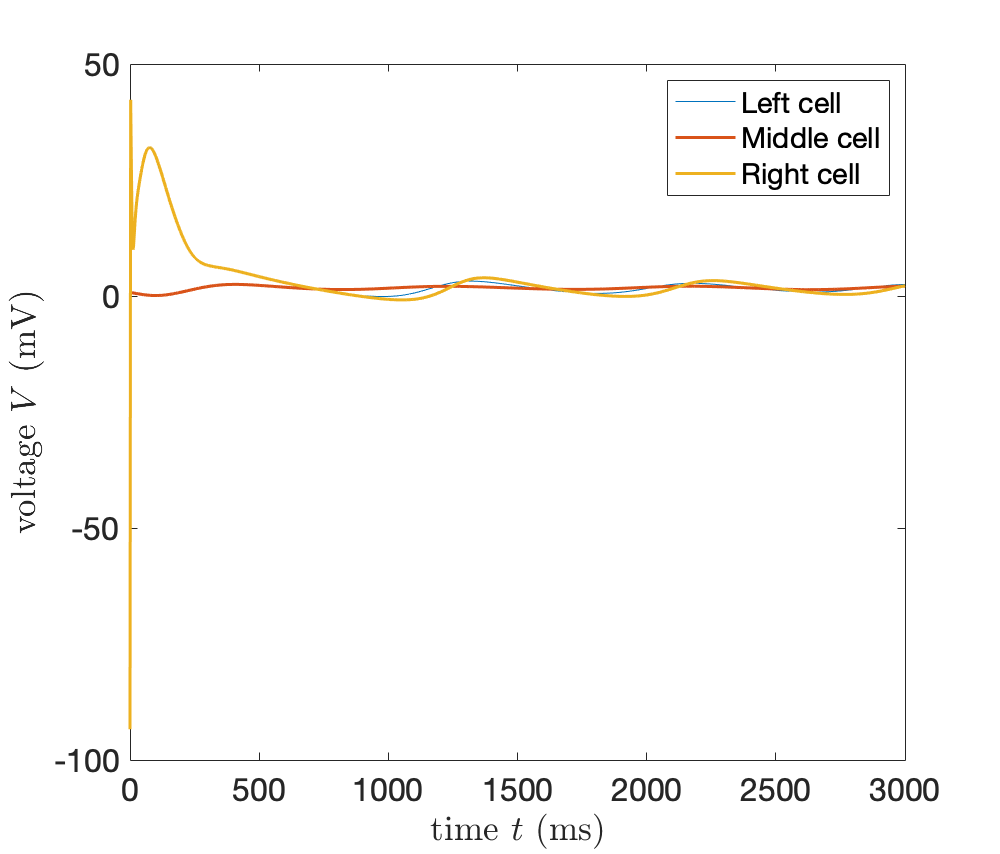}}}
\caption{Simulations of system~\eqref{monodomain_bernus}. In $\mathcal{D}$ the calcium conductance is $G_{\mathrm{Ca}}=0.0962518\frac{mS}{cm^2}$ and the stimulus to $I_\mathrm{stimulus}=0$, while in the remaining parts $G_{\mathrm{Ca}}=0.064\frac{mS}{cm^2}$ and $I_\mathrm{stimulus}=40\frac{\mu A}{cm^2}$ in (a)--(b), and $G_{\mathrm{Ca}}=0.09616\frac{mS}{cm^2}$ and $I_\mathrm{stimulus}=40\frac{\mu A}{cm^2}$ in (c). In all simulations $\mathcal{D}=[0.2, 0.7)$. Top line: the 1D cable simulations. Bottom line: the corresponding left, middle, and right cell (blue, red, and yellow). Computed with $128$ cells.} \label{fig:mono-bernus-chaos}
\end{figure}

{Finally, we briefly study the effects of initially setting the cells in the region $\mathcal{D}=[0.2, 0.7)$ to be prone to chaotic behaviour. From Figure \ref{fig:mono-bernus-chaos} we observe that chaotic behaviour does not spread in the three cases considered. In the first two simulations (Figure \ref{fig:mono-bernus-chaos})(a)--(b)) the surrounding cells are set to produce normal APs, which they indeed do for both the diffusion constant $\frac{1}{360}$\,{\it{mS}} and $0.00005$\,{\it{mS}}. However, if the surrounding cells are set to produce normal APs but are close to the EAD setting, we see that the dynamics die out ($V \to 0$).}

{
In conclusion, chaotic behaviour does not spread even if $50 \%$ of the cells are initially chaotic. However, cardiac death can occur in the model ($V \to 0$), see Figure \ref{fig:mono-bernus-chaos}(c). In Figure \ref{fig:mono-bernus-chaos}(a) the system produces normal APs (apart from the fact that the middle cells are initially set to zero). However, reducing the diffusion constant, the non-physiological behaviour that we saw in the previous section appears, see Figure \ref{fig:mono-bernus-chaos}(b). This indicates a critical lower bound on the diffusion constant in this model for the dynamics to be physiologically relevant. 

Concerning EADs, the above results indicate that the spreading of EADs on the tissue level is depending on 1. the number of cells prone to establish EADs, 2. how close to the normal setting the surrounding cells are, and 3. the diffusivity of the monodomain model~\eqref{monodomain_bernus}.}
 
\section{Summary and discussion}\label{sec:conclusion} 

In this paper, we investigated and analysed the behaviour of two mathematical models describing the action potentials of a Purkinje and a human ventricular cardiac muscle cell. To this end, we utilised bifurcation theory, numerical bifurcation analysis, and computational tools to establish an increased understanding of the dynamics of these models. This enabled us to find hidden features in both models considered. Furthermore, carrying out this analysis, we aimed at convincing the reader that 1. bifurcation analysis is very beneficial in the study of the dynamics of an ODE model and to detect hidden features of the considered system, 2. it is important to know how to interpret the corresponding bifurcation diagram, and 3. advancing cardiac cell research benefits from collaborations between mathematicians and physiologists/biologists.

First, we studied the dynamics of the Noble model~(\ref{model}) with respect to the leak current $I_\mathrm{L}$ based on the discussion in~\cite{Noble}. In~\cite{Noble} the author already varied the leak current conductance $G_\mathrm{L}$, resulting in the observation that the conductance $G_\mathrm{L}$ influences the period of the AP. Even more, if one chooses $G_\mathrm{L}$ large enough, e.g. $G_\mathrm{L}=0.4\ \frac{mS}{cm^2}$, the system converges into a stable equilibrium and no AP can appear, cf.~\cite{Noble}. This behaviour was analysed in more detail using numerical bifurcation theory. It turns out that this system changes stability via a supercritical Andronov--Hopf bifurcation from which a stable limit cycle branch bifurcates. This limit cycle branch loses and wins stability via limit point of cycle bifurcations and a period doubling bifurcation, respectively. Moreover, from the first period doubling bifurcations of the second limit cycle branch a (stable) period doubling cascade bifurcates, which is also the route to chaos. Interestingly, every limit cycle branch contains two period doubling bifurcations, which are connected via two limit cycle branches. Dependent on the initial values and the choice of $G_\mathrm{L}$, the Noble model~(\ref{model}) exhibits complex patterns and chaos. However, although mathematically interesting, the chaos detected (see Figure \ref{fig:chaos}) might be an artefact of the mathematical model and not within the physiologically relevant range. 

We used the same approach to study the more complex model~\eqref{model_bernus}, more complex in the sense that it contains more ion currents, pumps and exchangers compared to model~(\ref{model}), including the missing calcium current $I_\text{Ca}$ and the fast and slow potassium current, $I_{\mathrm{K}_\mathrm{r}}$ and $I_{\mathrm{K}_\mathrm{s}}$. On the single cell level it turned out that system~\eqref{model_bernus} exhibits both chaotic behaviour and EADs via a combination of a reduction of the fast potassium current and an enhanced calcium current. 

We would like to remark that this approach is also applicable for more modern and advanced models as the ones in \cite{TP06,Mahajan,Solittle,Ohara}, provided that the model of interest is regular and smooth enough. However, some up-to-date models, as in \cite{Mahajan} or \cite{Ohara}, might cause issues due to their complexity and lack of smoothness.


For system~\eqref{model_bernus}, we showed that for a 80\% block of the fast potassium current and an enhanced calcium current, system~\eqref{model_bernus} exhibits a subcritical Andronov--Hopf bifurcation from which an unstable limit cycle branch bifurcates, which stays unstable and contains a period doubling bifurcation. From this period doubling bifurcation a stable period doubling cascade bifurcates, which causes both EADs and deterministic chaos. Note that also for other (fast) potassium block rates system~\eqref{model_bernus} contains a subcritical Andronov--Hopf bifurcation and may exhibit EADs and/or chaos. {As experiments have shown that this is how EADs can occur, we can interpret these findings as a further validation of the model. However, it is unclear whether the occurring chaotic and self-oscillatory patterns also appear in a biological cardiac muscle cell (or are modelling artefacts). In particular, the dynamics in Figure \ref{fig:bernus_chaos}(a) is likely not appearing within a physiologically relevant range. This has to be investigated further with accurate experiments.


We would also like to highlight that further investigations using bifurcation analysis can be performed in the study of, for example, the potassium or calcium dynamics of cardiac cell models. For instance, system~\eqref{model_bernus} contains the calcium current and the equilibrium potential of the calcium current is
     $$
     E_\mathrm{Ca}=\frac{RT}{2F}\log\left(\frac{[\mathrm{Ca}]_e}{[\mathrm{Ca}]_i}\right).
     $$
     Different equilibrium potentials may change the behaviour of the considered system, see~\cite{MR4102409}, since the ion currents are depending on the equilibrium potentials, i.e.
     $$
     I_\mathrm{Ca}=G_\mathrm{Ca}\cdot d_\infty\cdot f\cdot f_\mathrm{Ca}(V-E_\mathrm{Ca}),
     $$
cf.~\cite{Bernus}. Thus, one can investigate complex calcium dynamics by, e.g., choosing the intra- or extracellular calcium concentrations, $[\mathrm{Ca}]_i$ or $[\mathrm{Ca}]_e$, as a bifurcation parameter. Indeed, more complex calcium dynamics is expected for more up-to-date models, such as the ones in \cite{TP06, Mahajan, Ohara}.

Besides comprehending the dynamics of the single cell models, it is crucial to understand how these dynamics affect the behaviour on the macro-scale ($cm$) due to the fact that multiple cardiac single cells may synchronise and cause arrhythmias. To this end, we introduced monodomain models for both systems, cf.~\eqref{monodomain} and \eqref{monodomain_bernus}. Based on the analysis of the single cell dynamics, we investigated cell synchronisation in both models. Both analyses showed that 1. the diffusivity of the model, 2. the number of cells, and 3. the placement of chaotic/EAD regions affect the global dynamics of the monodomain models. Furthermore, the non-physiological behaviour observed for single cells can transfer over to the macro-scale models, but this depends on the size of the diffusion constant. In particular, this warrants a bifurcation analysis with respect to the diffusion parameter and the size of the unstable regions for the two models considered.
An analysis finding a criterion for instabilities has been performed for the Luo--Rudy model~\cite{LR} in \cite{MR2445240}. Whether a criterion like this is obtainable for the Bernus model (and for more complex models) remains an open question.

\section{Conclusion}
Bifurcation theory is itself a powerful tool to study the behaviour of dynamical systems. It is used in many contexts and gets more and more attention also in (mathematical and computational) cardiac and neuroscience. In cardiac science, one can use this approach to establish a better understanding of cardiac arrhythmia~\cite{Kurata,Tran,Tsumoto2017,Xie,AE_MMOs,KBE,Ae_control}. In interdisciplinary research, bifurcation theory can also be an important component in successful treatment of human diseases. 

We considered two specific cardiac cell models among a multitude. Although a large number of models exist, there is still a lot of future work to do to derive a complete understanding of all cardiac dynamics. One has to deal with several issues, e.g. complexity of realistic models, and numerical and computational problems. In the very end, the extension from the cellular level to the tissue level has to be understood~\cite{computing_heart,Tri-domain,10.1371/journal.pcbi.1007042,Kaboudianeaav6019,niederer2019computational}. In particular, the physiological relevance of complex dynamics at all levels of modelling, as the ones we consider in this paper, has to be proven by experimental data.

In conclusion, this research would benefit from close interdisciplinary collaborations, since 1. a good, robust and realistic mathematical model can be developed based on experimental data, 2. an in-depth mathematical analysis can validate the accuracy or display weaknesses of the model, 3. these new findings either would help to improve the modelling or reformulation of the system to derive a most realistic model, including all expected dynamics, 4. the analysis would also increase the understanding of the occurrence of certain phenomena, and 5. the new obtained knowledge would help to develop new potential treatments for human diseases.

\vspace{0.5cm}
\noindent \small{\textbf{Conflict of interest} The authors declare no conflict of interest.}
\\[2ex]
\noindent \small{\textbf{Acknowledgements} 
A.E., supported by the \textit{Kristine Bonnevie scholarship 2020} during his research stay at Lund University in 2020, wishes to thank Erik Wahl\'en and the Centre of Mathematical Sciences, Lund University, Sweden for hosting him. 
}
\bibliography{erhardt_noble62_mybibfile}

\end{document}